\documentclass[12pt]{elsarticleNoFoot}

\usepackage[a4paper,text={16cm,24cm},centering]{geometry}
\usepackage[latin1]{inputenc}
\usepackage[pdftex]{color}
\usepackage{graphicx}
\usepackage{amssymb}
\usepackage[hidelinks]{hyperref}

\renewcommand{\arraystretch}{1.0}

\makeatletter

\def\newlog#1#2{%
	\def#1{\mathop{\operator@font #2}\nolimits}%
}
\def\Ite#1#2#3{\setbox\@tempboxa\hbox{$#1$}%
\ifdim\wd\@tempboxa=0cm#3\else#2\fi\relax}

\makeatother

\definecolor{coLabel}	{rgb}{0.50,0.00,0.00}	
\definecolor{coRemark}	{rgb}{0.99,0.00,0.00}	

\definecolor{coRef}	{rgb}{0.00,0.00,0.99}	
\definecolor{coThmName}	{rgb}{0.00,0.00,0.99}
\definecolor{coQDef}	{rgb}{0.00,0.00,0.60}	
\definecolor{coQRef}	{rgb}{0.00,0.00,0.80}	
\definecolor{coQNoref}	{rgb}{0.00,0.00,0.20}	

\definecolor{coLowlight}{rgb}{0.50,0.50,0.50}	
\definecolor{coLowlowlight}{rgb}{0.75,0.75,0.75}

\newcommand{\LABEL}[1]{\label{#1}}

\newcommand{\REF}[1]{\textcolor{coRef}{\ref{#1}}}
\newcommand{\REFF}[2]{\textcolor{coRef}{\ref{#1}.\ref{#1 #2}}}

\newcommand{\OEF}[3]{\cite[#1, p.#2]{Burghardt.2018c}}

\newcommand{\La}{\Leftarrow}
\newcommand{\Lra}{\Leftrightarrow}
\newcommand{\ra}{\rightarrow}
\newcommand{\la}{\leftarrow}
\newcommand{\lra}{\leftrightarrow}
\newcommand{\Ra}{\Rightarrow}
\newcommand{\raa}{\longrightarrow}
\renewcommand{\leq}{\leqslant}
\renewcommand{\geq}{\geqslant}

\newlog{\dom}{dom}
\newlog{\ran}{ran}

\newcommand{\tpl}[1]{\langle #1 \rangle}	
\newcommand{\set}[1]{\{ #1 \}}			
\newcommand{\disjUnion}{\mathop{\stackrel{.}{\cup}}}
\newcommand{\false}{{\it false}}
\newcommand{\true}{{\it true}}
\newcommand{\N}{I\!\!N}
\newcommand{\Z}{Z\!\!\!Z}
\newcommand{\Q}{Q\!\!\!\!Q}

\newcommand{\idiv}{/\!\!/}

\newcommand{\compose}{\circ}	

\newcommand{\filename}[1]{\texttt{#1}}

\newcommand{\lxor}{\oplus}

\newcommand{\QDef}[1]{\label{#1}\textcolor{coQDef}{\textsf{#1}}}%

\newcommand{\QRef}[1]{\hyperref[#1]{\textcolor{coQRef}{\textsf{#1}}}}%

\newcommand{\QRem}[1]{\textcolor{coQNoref}{\textsf{#1}}}%

\newcommand{\rr}[1]{#1}

\newcommand{\ro}[2]{\rr{{#1}^{#2}}}

\renewcommand{\:}[4]{%
	{%
	\renewcommand{\:}[4]{%
		{%
		\renewcommand{\:}[4]{error\error}%
		\renewcommand{\j}{{##2}}%
		{##4}%
		##1...##1%
		\renewcommand{\j}{{##3}}%
		{##4}%
		}%
	}%
	\renewcommand{\i}{{#2}}%
	{#4}%
	#1\ldots#1%
	\renewcommand{\i}{{#3}}%
	{#4}%
	}%
}

\newproof{pf}{Proof}
\newcommand{\THM}[3]{%
	\vspace{0ex plus 1ex}%
	\begin{#1}%
	\Ite{#2}{ \textcolor{coThmName}{(#2)} }{ }%
	{\rm #3}%
	\end{#1}%
	}
\newcommand{\PRF}[2]{%
	\vspace{-2ex plus 1ex}%
	\begin{pf}%
	{#2}%
	\end{pf}%
	\vspace{0ex plus 1ex}%
	}
\newcommand{\LEMMA}[2]		{\THM{lemma}		{#1}{#2}}
\newcommand{\THEOREM}[2]	{\THM{theorem}		{#1}{#2}}
\newcommand{\DEFINITION}[2]	{\THM{definition}	{#1}{#2}}
\newcommand{\EXAMPLE}[2]	{\THM{example}		{#1}{#2}}

\newcommand{\PROOF}[1]		{\PRF{\em Proof. }	{#1}}

\renewcommand{\.}[1]{\!#1\!}

\begin{document}

\newsavebox{\atpic}
\savebox{\atpic}{\includegraphics[scale=0.08]{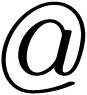}}

\begin{frontmatter}

\title{An Algebra of Properties of Binary Relations}
\author{Jochen Burghardt}
\address{jochen.burghardt\usebox{\atpic}alumni.tu-berlin.de}
\address{\rm Feb 2021}

\begin{abstract}
We consider all $16$ unary
operations that, given a homogeneous binary relation
$R$, define a new one by a boolean combination of $xRy$ and $yRx$.
Operations can be composed, and connected by pointwise-defined logical
junctors.
We consider the usual properties of relations,
and allow them to be lifted by prepending an operation.

We investigate extensional equality between lifted properties
(e.g.\ {\em a relation is connex iff its complement is asymmetric\/}),
and give a table to decide this equality.
Supported by a counter-example generator
and a resolution theorem prover,
we investigate all $3$-atom implications between lifted properties,
and give a sound and complete axiom set for them (containing e.g.\
{\em ``if $R$'s complement is left Euclidean and $R$ is right serial, then
$R$'s symmetric kernel is left serial''\/}).
\end{abstract}


\end{frontmatter}

\vfill

{

\setlength{\parskip}{0cm}

\tableofcontents


\listoffigures

}

\clearpage

\section{Introduction}
\LABEL{Introduction}

We strive to compute new laws about homogeneous binary relations.
In \cite{Burghardt.2018c},
we considered the $24$ ``best-known'' basic properties
of relations, and
used a counter-example search combined with
the Quine-McCluskey algorithm and followed by a manual
confirmation / disconfirmation process, to obtain
a complete list of laws of the form
$$\forall R. \; \:{\lor}1n{{\it prop}_\i(R)} ,$$
where ${\it prop}_i$ are negated or unnegated properties,
and $R$ ranges over all relations.
A classical example of such a law is ``{\em each irreflexive and
transitive relation is asymmetric}'';
a lesser known example is ``{\em  a left unique and right serial
relation (on a domain of $\geq 5$ elements)
cannot be incomparability-transitive}''.
Once such a law has been stated it is usually easy to prove.

We also briefly sketched in \cite[p.6-7]{Burghardt.2018c}
a more general approach based on regular tree grammars
which could obtain laws involving arbitrary
operators on relations (union, closures, \ldots);
however, this approach would require infeasible large
computation time and memory.

In the current
report, we investigate an intermediate approach which allows
one to detect laws involving certain unary operations on relations.
The considered set of operations can be motivated
by the definition of quasi-transitivity:
a relation $R$ is quasi-transitive iff ${\it op}_2(R)$ is transitive,
where the relation ${\it op}_2(R)$ is defined as
$$x \mathrel{{\it op}_2(R)} y 
\mbox{ ~ iff ~ }
x \mathrel{R} y 
\mbox{ and not } 
y \mathrel{R} x .$$
Slightly more general,
we will allow all unary operations
that can be defined by a boolean combination of 
$x \rr{R} y$ and $y \rr{R} x$.
Since there are only $16$ such operations, generating law suggestions
will be still feasible.
This allows one to search for
laws of the form
$$\forall R. \;
\:{\lor}1n{{\it prop}_\i( \:{\circ}1{m_\i}{\it op_{\i\j}}  (R))} ,$$
An example law is
``{\em the symmetric closure of a left-serial and transitive relation
is always dense''};
once it is stated, it is straight forward to prove
(Lem.~\REFF{3-Implications group 1}{qudq}).

Our set of unary operations is closed w.r.t.\ composition 
(``$\compose$'', Def.~\REF{Operation composition}),
therefore it is
sufficient to look for laws of the form
$$\forall R. \; \:{\lor}1n{{\it prop}_\i({\it op_\i}(R))} .$$
We define a ``lifted property'' to be the composition of a basic
property and a unary operation, that is, an expression of the form
$\lambda R. \; {\it prop}({\it op}(R))$ in Lambda-calculus notation.
Considering all $24 \cdot 16=384$
possible combinations of basic properties (Def.~\REF{def})
and unary operations (Def.~\REF{optn def}), it turns
out that many combinations are extensionally equal,
and there are only $81$ different
ones (Thm.~\REF{optn equiv}).

Found laws may now involve e.g.\ the converse, the complement,
the symmetric closure, and the symmetric kernel of a relation.\footnote{
	Note that we can't express e.g.\ reflexive closure,
	since our unary operations aren't concerned with information
	about $x \rr{R} x$.
	To achieve this, our approach could be extended
	to all $65536$ boolean combinations of $x \rr{R} y$, 
	$y \rr{R} x$, $x \rr{R} x$,
	and $y \rr{R} y$; 
	however, this is likely to be beyond feasibility
	again.
}
Moreover, some basic
properties are rendered redundant, as they can be expressed
by lifting other properties
(Lem.~\REF{optn redundant}).
For example, a relation is empty and co-reflexive
iff its symmetric closure is asymmetric and anti-symmetric,
respectively.

While the basic properties were obtained from mathematical
applications (order, equivalence, \ldots; cf.\ Def.~\REF{kinds})
developed over the years,
and therefore pretty much ad hoc,
considering all $81$ lifted properties is slightly more systematic;
the latter are closed w.r.t.\ our unary operations,
while the former are not.
A completely systematic approach would consider each property
definable by a predicate logic formula with a given number of
quantifiers, rather than just the $24$ of them shown in
Def.~\REF{def}; however, such an approach is, again,
computationally infeasible.
So, the unary-operations approach presented here is a good
compromise;
we will demonstrate in this report that it is in fact feasible.

In Sect.~\REF{Definitions},
we give some basic definitions prior to the formal introduction of our
unary-operation approach in
Sect.~\REF{An algebra of unary operations on relations}.
In Sect.~\REF{Equivalent lifted properties},
we compute the $81$ equivalence classes (w.r.t.\ extensional equality)
of lifted properties.
We also define a default representative of each class.
While the Quine-McCluskey approach from \cite{Burghardt.2018c} is
still infeasible for $81$ different lifted properties, we investigate,
in Sect.~\REF{Implications between lifted properties},
all laws of the form
$$\forall R. \; {\it lprop}_1(R) \land {\it lprop}_2(R)
	\ra {\it lprop}_3(R) ,$$
where ${\it lprop}_i$ are (unnegated) lifted properties.
Using an approach different from Quine-McCluskey to eliminate
redundant laws, we isolate in Sect.~\REF{Axiom set}
a total of $124$ ``axioms'' which imply all found laws 
(Thm.~\REF{$3$-implication axioms}).
Our \texttt{C} source code is provided in the ancillary files.

\clearpage

\section{Definitions}
\LABEL{Definitions}

In this section, we give some preparatory definitions.
Definition~\REF{def} defines the ``basic''
properties we consider throughout this report.
Definition~\REF{kinds} indicates their historical origins,
and at the same time names some applications for them.
Definition~\REF{optn names def} introduces some operators on
relations; each of them will later be identified with an admitted
unary operation (Lem.~\REF{optn names lem}).
Definition~\REF{Monotonicity} defines the notion of monotonic and
antitonic properties; Lem.~\REF{Monotonic properties} classifies our
basic properties by these categories.

\DEFINITION{Binary relation properties}{%
\LABEL{def}%
Let $X$ be a set.
A (homogeneous) binary relation $R$ on $X$ is a subset of $X \times X$.
The relation $R$ is called
\begin{enumerate}
\item\LABEL{def Refl}
	reflexive
	(``Refl'', ``rf'')
	~ if ~ $\forall x \in X. \;\; x \rr{R} x$;
\item\LABEL{def Irrefl}
	irreflexive
	(``Irrefl'', ``ir'')
	~ if ~ $\forall x \in X. \;\; \lnot x \rr{R} x$;
\item\LABEL{def CoRefl}
	co-reflexive
	(``CoRefl'', ``cr'')
	~ if ~ $\forall x,y \in X. \;\; x \rr{R} y \ra x=y$;
\item\LABEL{def LfQuasiRefl}
	left quasi-reflexive
	(``lq'')
	~ if ~ 
	$\forall x,y \in X. \;\; x \rr{R} y \ra x \rr{R} x$;
\item\LABEL{def RgQuasiRefl}
	right quasi-reflexive
	(``rq'')
	~ if ~ 
	$\forall x,y \in X. \;\; x \rr{R} y \ra y \rr{R} y$;
\item\LABEL{def QuasiRefl}
	quasi-reflexive
	(``QuasiRefl'')
	~ if ~ it is both left and right quasi-reflexive;
\item\LABEL{def Sym}
	symmetric
	(``Sym'', ``sy'')
	~ if ~ 
	$\forall x,y \in X. \;\; x \rr{R} y \ra y \rr{R} x$;
\item\LABEL{def ASym}
	asymmetric
	(``ASym'', ``as'')
	~ if ~ 
	$\forall x,y \in X. \;\; 
	x \rr{R} y \ra \lnot y \rr{R} x$;
\item\LABEL{def AntiSym}
	anti-symmetric
	(``AntiSym'', ``an'')
	~ if ~ 
	$\forall x,y \in X. \;\; 
	x \rr{R} y \land x \neq y \ra \lnot y \rr{R} x$;
\item\LABEL{def SemiConnex}
	semi-connex
	(``SemiConnex'', ``sc'')
	~ if ~ 
	 $\forall x,y \in X. \;\; 
	 x \rr{R} y \lor y \rr{R} x \lor x=y$;
\item\LABEL{def Connex}
	connex
	(``Connex'', ``co'')
	~ if ~ 
	$\forall x,y \in X. \;\; x \rr{R} y \lor y \rr{R} x$;
\item\LABEL{def Trans}
	transitive
	(``Trans'', ``tr'')
	~ if ~ 
	$\forall x,y,z \in X. \;\; 
	x \rr{R} y \land y \rr{R} z \ra x \rr{R} z$;
\item\LABEL{def AntiTrans}
	anti-transitive
	(``AntiTrans'', ``at'')
	~ if ~ 
	$\forall x,y,z \in X. \;\; 
	x \rr{R} y \land y \rr{R} z \ra \lnot x \rr{R} z$;
\item\LABEL{def QuasiTrans}
	quasi-transitive
	(``QuasiTrans'', ``qt'')
	~ if ~
	\\
	$\forall x,y,z \in X. \;\;
	x \rr{R} y \land \lnot y \rr{R} x 
	\land y \rr{R} z \land \lnot z \rr{R} y
	\ra x \rr{R} z \land \lnot z \rr{R} x$;
\item\LABEL{def RgEucl}
	right Euclidean
	(``RgEucl'', ``re'')
	~ if ~ 
	$\forall x,y,z \in X. \;\; 
	x \rr{R} y \land x \rr{R} z \ra y \rr{R} z$;
\item\LABEL{def LfEucl}
	left Euclidean
	(``LfEucl'', ``le'')
	~ if ~ 
	$\forall x,y,z \in X. \;\; 
	y \rr{R} x \land z \rr{R} x \ra y \rr{R} z$;
\item\LABEL{def SemiOrd1}
	semi-order property~1
	(``SemiOrd1'', ``s1'')
	~ if ~
	\\
	$\forall w,x,y,z \in X. \;\;
	w \rr{R} x 
	\land \lnot x \rr{R} y \land \lnot y \rr{R} x 
	\land y \rr{R} z 
	\ra w \rr{R} z$;
\item\LABEL{def SemiOrd2}
	semi-order property~2
	(``SemiOrd2'', ``s2'')
	~ if ~
	\\
	$\forall w,x,y,z \in X. \;\;
	x \rr{R} y \land y \rr{R} z
	\ra w \rr{R} x \lor x \rr{R} w 
	\lor w \rr{R} y \lor y \rr{R} w 
	\lor w \rr{R} z \lor z \rr{R} w$.
\item\LABEL{def RgSerial}
	right serial
	(``RgSerial'', ``rs'')
	~ if ~ $\forall x \in X \; \exists y \in X. \;\; x \rr{R} y$
\item\LABEL{def LfSerial}
	left serial
	(``LfSerial'', ``ls'')
	~ if ~ $\forall y \in X \; \exists x \in X. \;\; x \rr{R} y$
\item\LABEL{def Dense}
	dense
	(``Dense'', ``de'')
	~ if ~ $\forall x, z \in X \; \exists y \in X. \;\;
	x \rr{R} z \ra x \rr{R} y \land y \rr{R} z$.
\item\LABEL{def IncTrans}
	incomparability-transitive
	(``IncTrans'', ``it'')
	~ if ~
	\\
	$\forall x, y, z \in X. \;\;
	\lnot x \rr{R} y \land \lnot y \rr{R} x 
	\land \lnot y \rr{R} z \land \lnot z \rr{R} y
	\ra \lnot x \rr{R} z \land \lnot z \rr{R} x$.
\item\LABEL{def LfUnique}
	left unique
	(``LfUnique'', ``lu'')
	~ if ~ $\forall x_1, x_2, y \in X \;\;
	x_1 \rr{R} y \land x_2 \rr{R} y \ra x_1 = x_2$.
\item\LABEL{def RgUnique}
	right unique
	(``RgUnique'', ``ru'')
	~ if ~ $\forall x, y_1, y_2 \in X \;\;
	x \rr{R} y_1 \land x \rr{R} y_2 \ra y_1 = y_2$.
\end{enumerate}
The capitalized abbreviations in parentheses are used by our
software; the two-letter codes are used in tables and pictures when
space is scarce.

We say that $x,y$ are incomparable w.r.t.\ $R$,
if $\lnot x \rr{R} y \land \lnot y \rr{R} x$ holds.
\qed
}

\DEFINITION{Kinds of binary relations}{%
\LABEL{kinds}%
A binary relation $R$ on a set $X$ is called
\begin{enumerate}
\item an equivalence
	~ if ~
	it is reflexive, symmetric, and transitive;
\item a partial equivalence
	~ if ~
	it is symmetric and transitive;

\item a tolerance relation
	~ if ~
	it is reflexive and symmetric;
\item idempotent
	~ if ~
	it is dense and transitive;
\item trichotomous
	~ if ~
	it is irreflexive, asymmetric, and semi-connex;

\item a non-strict partial order
	~ if ~
	it is reflexive, anti-symmetric, and transitive;
\item a strict partial order
	~ if ~
	it is irreflexive, asymmetric, and transitive;
\item\LABEL{kinds semi-order}
	a semi-order
	~ if ~
	it is asymmetric and satisfies semi-order properties 1 and 2;
\item a preorder
	~ if ~
	it is reflexive and transitive;
\item a weak ordering
	~ if ~
	\\
	it is irreflexive, asymmetric, transitive, and
	incomparability-transitive;

\item a partial function
	~ if ~
	it is right unique;
\item a total function
	~ if ~
	it is right unique and right serial;
\item an injective function
	~ if ~
	it is left unique, right unique, and right serial;
\item a surjective function
	~ if ~
	it is right unique and and left and right serial;
\item a bijective function
	~ if ~
	it is left and right unique and left and right serial.
\qed
\end{enumerate}
}

\DEFINITION{Operation names}{%
\LABEL{optn names def}%
\begin{enumerate}
\item\LABEL{optn names def 1}%
	The symmetric kernel of a relation $R$ is defined as the
	largest subset of $R$ that is a symmetric relation.
\item\LABEL{optn names def 2}%
	The symmetric closure of a relation $R$ is defined as the
	smallest superset of $R$ that is a symmetric relation.
\item\LABEL{optn names def 3}%
	We define the asymmetric kernel of a relation $R$ as the
	intersection of all maximal subsets of $R$ that are
	asymmetric relations.

	Note that for an arbitrary relation $R$,
	a largest subset that is an asymmetric relation
	need not exist.
	For example, on the set $X = \set{0,1}$,
	the relation $R = \set{ \tpl{0,1}, \tpl{1,0} }$
	has three asymmetric subsets, viz.\
	$R_1 = \set{ \tpl{0,1} }$,
	$R_2 = \set{ \tpl{1,0} }$, and
	$R_3 = \set{}$.
	While $R_1$ and $R_2$ are maximal w.r.t.\ $\subseteq$,
	none of them is largest.
\qed
\end{enumerate}
}

\DEFINITION{Monotonicity}{%
\LABEL{Monotonicity}%
A property ${\it prop}$ of binary relations is called monotonic
if
$\forall R_1, R_2. \;
R_1 \subseteq R_2 \land {\it prop}(R_1) \Ra {\it prop}(R_2)$.
It is called antitonic if
$\forall R_1, R_2. \;
R_1 \supseteq R_2 \land {\it prop}(R_1) \Ra {\it prop}(R_2)$.
\qed
}

\definecolor{coMonNeeded}	{rgb}{0.00,0.80,0.00}
\definecolor{coMonDestroys}	{rgb}{0.80,0.00,0.00}

\LEMMA{Monotonic properties}{%
\LABEL{Monotonic properties}%
We use first-order formulas without constants and function symbols and
with equality and one binary relation symbol $R$ to define properties
of binary relations, as in Def.~\REF{def}.
Such a formula is called an $\land$-$\lor$-formula if it is
in prenex normal form,
contains no other binary junctors than $(\land)$ and $(\lor)$,
and contains $(\lnot)$ only applied to atoms.
\begin{enumerate}
\item\LABEL{Monotonic properties 1}%
	A property is monotonic if its definition can be written as
	a closed $\land$-$\lor$-formula
	without any negated occurrence of $R$.
\item\LABEL{Monotonic properties 2}%
	A property is antitonic if its definition can be written as
	a closed $\land$-$\lor$-formula
	without any unnegated occurrence of $R$.
\item\LABEL{Monotonic properties 3}%
	The following basic properties from Def.~\REF{def} are monotonic:
	\\
	Refl,
	SemiConnex,
	Connex,
	RgSerial,
	LfSerial.
\item\LABEL{Monotonic properties 4}%
	The following basic properties are antitonic:
	\\
	Irrefl,
	Corefl,
	ASym,
	AntiSym,
	AntiTrans,
	LfUnique,
	RgUnique.
\item\LABEL{Monotonic properties 5}%
	The following basic properties are neither monotonic nor antitonic:
	\\
	LfQuasiRefl,
	RgQuasiRefl,
	QuasiRefl,
	Sym,
	Trans,
	QuasiTrans,
	RgEucl,
	LfEucl,
	SemiOrd1,
	SemiOrd2,
	Dense,
	IncTrans.
\end{enumerate}
}
\PROOF{
\begin{enumerate}
\item
	Given an
	$\land$-$\lor$-formula $\psi$ in $n$ free variables
	$x_1,\ldots,x_n$ and a relation $R$ on a domain $X$,
	define $M(R,\psi)$ to be the set of all $n$-tuples in $X^n$
	satisfying $\psi$, w.r.t.\ $R$.
	For an $\land$-$\lor$-formula
	$\psi$ without quantifiers or negated occurrences of $R$,
	show by induction on the structure of $\psi$ that
	$R_1 \subseteq R_2$
	implies $M(R_1,\psi) \subseteq M(R_2,\psi)$:
	\begin{itemize}
	\item If $\psi$ has the form $\psi_1 \land \psi_2$,
		then
		\\
		$M(R_1,\psi)
		= M(R_1,\psi_1) \cap M(R_1,\psi_2)
		\subseteq M(R_2,\psi_1) \cap M(R_2,\psi_2)
		= M(R_2,\psi)$.
	\item The case $\psi_1 \lor \psi_2$ is similar, relying on the
		monotonicity of $\cup$, rather than of $\cap$.
	\item If $\psi$ has the form $x_i \rr{R} x_j$,
		then
		$M(R_1,\psi)
		= \set{ \tpl{ x_1, ..., x_i, ..., x_j, ..., x_n}
			\mid x_i~\rr{R_1}~x_j }
		\subseteq \set{ \tpl{ x_1, ..., x_i, ..., x_j, ..., x_n}
			\mid x_i~\rr{R_2}~x_j }
		= M(R_2,\psi)$.
	\item If $\psi$ does not contain $R$,
		then $M(R_1,\psi) = M(R_2,\psi)$.
		\\
		Note that $\psi$ needn't be an atom in this case.
	\end{itemize}
	Finally, show that $M(R_1,\psi) \subseteq M(R_2,\psi)$
	implies both
	$M(R_1,\forall x_n. \psi) \subseteq M(R_2,\forall x_n. \psi)$
	and
	$M(R_1,\exists x_n. \psi) \subseteq M(R_2,\exists x_n. \psi)$.
	By induction on $n$, this extends to arbitrarily long quantifier
	prefixes.
\item
	The proof is similar to~\REF{Monotonic properties 1}.
\item
	For each of the listed properties, its definition in
	Def.~\REF{def} is in the form required
	by~\REF{Monotonic properties 1}.
\item
	The definition of irreflexivity in
	Def.~\REF{def} is in the form required
	by~\REF{Monotonic properties 2}.
	For the remaining properties, resolving $(\ra)$ brings it into
	that form:
	\begin{itemize}
	\item Corefl: $\forall x, y \in X. \; 
		\lnot x \rr{R} y \lor x=y$
	\item ASym: $\forall x, y \in X. \; 
		\lnot x \rr{R} y \lor \lnot y \rr{R} x$
	\item AntiSym: $\forall x, y \in X. \;
		\lnot x \rr{R} y \lor \lnot y \rr{R} x \lor x=y$
	\item AntiTrans: $\forall x, y, z \in X. \;
		\lnot x \rr{R} y 
		\lor \lnot y \rr{R} z 
		\lor \lnot x \rr{R} z$.
	\item LfUnique: $\forall x_1, x_2, y \in X. \;
		\lnot x_1Ry \lor \lnot x_2Ry \lor x_1=x_2$.
	\item RgUnique: $\forall x, y_1, y_2 \in X. \;
		\lnot x \rr{R} y_1 
		\lor \lnot x \rr{R} y_2 
		\lor x_1=x_2$.
	\end{itemize}
\item
	\newcommand{\2}[1]{\textcolor{coMonNeeded}{#1}}
	\newcommand{\3}[1]{\textcolor{coMonDestroys}{#1}}
	For each property, we give relations
	$R_1 \subset R_2 \subset R_3$
	on the domain $X = \set{0,1,2,3}$
	such that $R_2$, but neither $R_1$ nor $R_3$,
	has the property.
	$R_1$ consists of all black pairs,
	$R_2$ consists of all black or green pairs, and
	$R_3$ consists of all pairs.
	Intuitively, adding the green pair establishes the property,
	and adding the red destroys it again.
	We use semi-colons to indicate the sub-relation separations in
	grey-scale renderings.
	\begin{itemize}
	\item LfQuasiRefl:
		$\set{ \tpl{0,1} \2{; \tpl{0,0}} \3{; \tpl{1,2}} }$
	\item RgQuasiRefl, QuasiRefl:
		$\set{ \tpl{0,0}, \tpl{0,1}
		\2{; \tpl{1,1}}
		\3{; \tpl{1,2}} }$
	\item Sym: $\set{ \tpl{0,1} \2{; \tpl{1,0}} \3{; \tpl{0,2}} }$
	\item Trans, QuasiTrans:
		$\set{ \tpl{0,1} , \tpl{1,2}
		\2{; \tpl{0,2}}
		\3{; \tpl{2,3}} }$
	\item RgEucl:
		$\set{ \tpl{0,1}
		\2{; \tpl{1,1}}
		\3{; \tpl{0,0}} }$
	\item LfEucl:
		$\set{ \tpl{0,1}
		\2{; \tpl{0,0}}
		\3{; \tpl{1,0}} }$
	\item SemiOrd1:
		$\set{ \tpl{0,0}, \tpl{1,1}
		\2{; \tpl{0,1}}
		\3{; \tpl{2,2}} }$
	\item SemiOrd2:
		$\set{ \tpl{0,1}, \tpl{1,2}
		\2{; \tpl{0,3}}
		\3{; \tpl{0,0}} }$
	\item Dense: $\set{ \tpl{0,1} \2{; \tpl{0,0}} \3{; \tpl{1,2}} }$
	\item IncTrans:
		$\set{ \tpl{0,1}, \tpl{1,2}
		\2{; \tpl{1,3}}
		\3{; \tpl{0,0}} }$
\qed
	\end{itemize}
\end{enumerate}
}

\clearpage

\section{An algebra of unary operations on relations}%
\LABEL{An algebra of unary operations on relations}

In this section, we introduce and investigate
our admitted unary operations on
relations.

\begin{figure}
\begin{center}
$$\begin{array}{|l|cccc|}
\hline
x \rr{R} y: & \false & \false & \true  & \true	\\
y \rr{R} x: & \false & \true  & \false & \true	\\
\hline
x \ro{R}{q} y: & q_8 & q_4 & q_2 & q_1	\\
\hline
\end{array}$$
\caption{Encoding of unary operations}
\LABEL{Encoding of unary operations}
\end{center}
\end{figure}

\DEFINITION{Unary operations}{%
\LABEL{optn def}%
We represent a unary operation on relations
as a $4$-bit number
$p \in \set{0, \ldots, 9, A, \ldots, F}$,
and denote its bits as $p_8, p_4, p_2, p_1$,
that is,
$$p = 8 \cdot p_8 + 4 \cdot p_4 + 2 \cdot p_2 + 1 \cdot p_1 .$$
Given a binary relation $R$ on a domain set $X$,
and $x, y \in X$,
we write $\ro{R}{p}$ to denote the application of $p$ to $R$, which we
define as
$$\begin{array}{cc@{\;}r@{\;}c@{\;}r@{\;}c@{\;}r@{\;}c}
     & \multicolumn{7}{l}{x \ro{R}{p} y}	\\
\Lra&(&\lnot x\rr{R}y&\land&\lnot y\rr{R}x&\land&p_8&) \\
\lor&(&\lnot x\rr{R}y&\land&      y\rr{R}x&\land&p_4&) \\
\lor&(&      x\rr{R}y&\land&\lnot y\rr{R}x&\land&p_2&) \\
\lor&(&      x\rr{R}y&\land&      y\rr{R}x&\land&p_1&) \\
\end{array}$$
tacitly identifying $0$ with $\false$ and $1$ with $\true$,
see Fig.~\REF{Encoding of unary operations}.
For example, $x \ro{R}{1} y$ is true 
iff both $x \rr{R} y$ and $y \rr{R} x$ is;
moreover,
${\it op_2}(R)$ from Sect.~\REF{Introduction} can now be written as
$R^2$.
Figure~\REF{Semantics of unary operations}
shows the semantics of each of the $16$ possible unary
operations.\footnote{
	We use $\lxor$ to denote exclusive or.
}
They allow us to express e.g.\ the converse, the complement,
the symmetric kernel, and the symmetric closure of a relation.
For unary operations $p$ and $q$,
we define $p \subseteq q$ bitwise as
$$p_8 \leq q_8
\land p_4 \leq q_4
\land p_2 \leq q_2
\land p_1 \leq q_1 .$$
Moreover,
we extend boolean connectives to unary operations in a bitwise manner,
e.g.\ $\lnot p$ is defined as bitwise complement:
$$\begin{array}{c*{20}{@{\;}c}}
(\lnot p) & = 
	& 8 & \cdot & (\lnot p_8) & + 
	& 4 & \cdot & (\lnot p_4) & + 
	& 2 & \cdot & (\lnot p_2) & + 
	& 1 \cdot & (\lnot p_1) \\
          & = & 8 & \cdot & (1 - p_8)   & + & 4 & \cdot & (1 - p_4)
	  & + & 2 & \cdot & (1 - p_2)   & + & 1 \cdot & (1 - p_1) & ,	\\
\end{array}$$
similar for the other connectives.
\qed
}

{
\newcommand{\8}{\scriptscriptstyle\mid\mid\mid\mid}
\newcommand{\4}{\la}
\newcommand{\2}{\ra}
\newcommand{\1}{\lra}
\newcommand{\0}{\color{coLowlowlight}}

It may be helpful to think of a unary operation as a set of 
graph rewriting rules.
This view is supported in column ``Rewriting'' of 
Fig.~\REF{Semantics of unary operations}.
Representing a binary relation as a directed graph, two given vertices
$x$ and $y$ can be connected by 
\begin{enumerate}
\item no edge at all (depicted $\8$),
\item just an edge from $y$ to $x$ (depicted $\4$),
\item just an edge from $x$ to $y$ (depicted $\2$), or
\item both an edge from $x$ to $y$ and a reverse edge (depicted $\1$).
\end{enumerate}
For each of the four situations, column ``Rewriting''
gives the appropriate replacement performed by an operation.
For example, operation 7 replaces all $\4$ and all $\2$
situations by $\1$,
and thus obtains the symmetric closure.
Since the column ``$\4$'' is just a mirror of ``$\2$'', it is grayed
out.
Observe that $0$ and $1$ in the most significant bit in column ``bin''
corresponds
to $\8$ and $\1$, respectively; similar for the least significant bit;
for the two middle bits, $00$, $01$, $10$, and $11$
correspond to
${\0\8}\;\8$, ${\0\4}\;\2$, ${\0\2}\;\4$, and ${\0\1}\;\1$, respectively.

\begin{figure}
\begin{center}
$$\begin{array}{|c|c@{}c@{}c@{}c|r@{\;}c@{\;}r|c@{\;}c@{\;}c@{\;}c|l|}
\hline
\multicolumn{5}{|c|}{\rm Optn}
	& \multicolumn{3}{c|}{\rm Formal}
	& \multicolumn{4}{c|}{\rm Rewriting}
	& \rm Intuitively	\\
\multicolumn{1}{|c}{\rm hx} & \multicolumn{4}{c|}{\rm bin}
	& \multicolumn{3}{c|}{}
	& \8 & \0\4 & \2 & \1
	& \\
\hline
\hline
0&0&0&0&0&\multicolumn{3}{c|}{\false}        & \8 & \0\8 & \8 & \8&\mbox{empty relation}	\\
1&0&0&0&1&      x\rr{R}y&\land&      y\rr{R}x& \8 & \0\8 & \8 & \1&\mbox{symmetric kernel}	\\
2&0&0&1&0&      x\rr{R}y&\land&\lnot y\rr{R}x& \8 & \0\4 & \2 & \8&\mbox{asymmetric kernel}	\\
3&0&0&1&1&      x\rr{R}y&     &              & \8 & \0\4 & \2 & \1&\mbox{identity}	\\
\hline
4&0&1&0&0&\lnot x\rr{R}y&\land&      y\rr{R}x& \8 & \0\2 & \4 & \8&\mbox{converse asymmetric kernel}	\\
5&0&1&0&1&              &     &      y\rr{R}x& \8 & \0\2 & \4 & \1&\mbox{converse}	\\
6&0&1&1&0&      x\rr{R}y&\lxor&      y\rr{R}x& \8 & \0\1 & \1 & \8&	\\
7&0&1&1&1&      x\rr{R}y&\lor &      y\rr{R}x& \8 & \0\1 & \1 & \1&\mbox{symmetric closure}	\\
\hline
8&1&0&0&0&\lnot x\rr{R}y&\land&\lnot y\rr{R}x& \1 & \0\8 & \8 & \8&\mbox{incomparable}	\\
9&1&0&0&1&      x\rr{R}y&\lxor&\lnot y\rr{R}x& \1 & \0\8 & \8 & \1&	\\
A&1&0&1&0&              &     &\lnot y\rr{R}x& \1 & \0\4 & \2 & \8&\mbox{complement of converse}	\\
B&1&0&1&1&      x\rr{R}y&\lor &\lnot y\rr{R}x& \1 & \0\4 & \2 & \1&\mbox{complement of converse asymmetric kernel}	\\
\hline
C&1&1&0&0&\lnot x\rr{R}y&     &              & \1 & \0\2 & \4 & \8&\mbox{complement}	\\
D&1&1&0&1&\lnot x\rr{R}y&\lor &      y\rr{R}x& \1 & \0\2 & \4 & \1&\mbox{complement of asymmetric kernel}	\\
E&1&1&1&0&\lnot x\rr{R}y&\lor &\lnot y\rr{R}x& \1 & \0\1 & \1 & \8&\mbox{complement of symmetric kernel}	\\
F&1&1&1&1&\multicolumn{3}{c|}{\true}         & \1 & \0\1 & \1 & \1&\mbox{universal relation}	\\
\hline
\end{array}$$
\caption{Semantics of unary operations}
\LABEL{Semantics of unary operations}
\end{center}
\end{figure}

}

\LEMMA{Boolean connectives on operations}{%
\LABEL{optn bool}%
\begin{enumerate}
\item\LABEL{optn bool 1}%
	$x R^{\lnot p} y$ ~ iff ~ $\lnot x \ro{R}{p} y$
\item\LABEL{optn bool 2}%
	$x \ro{R}{p} y \land x \ro{R}{q} y$ ~ iff ~ $x R^{p \land q} y$
\item\LABEL{optn bool 3}%
	Any other boolean connective
	distributes over operation application in a similar way.
\item\LABEL{optn bool 4}%
	If $p \subseteq q$,
	then $R^p \subseteq R^q$.
\end{enumerate}
}
\PROOF{
\begin{enumerate}
\item[1,2]
	We distinguish four cases:
	\begin{itemize}
	\item $      x \rr{R} y \land       y \rr{R} x$:	\\
		Then ~ $x R^{\lnot p} y$
		~ iff ~ $(\lnot p_1) = 1$
		~ iff ~ $p_1 = 0$
		~ iff ~ $\lnot x \ro{R}{p} y$.
		\\
		And ~ $x \ro{R}{p} y \land x \ro{R}{q} y$
		~ iff ~ $(p_1 = 1) \land (q_1 = 1)$
		~ iff ~ $(p_1 \land q_1) = 1$
		~ iff ~ $x R^{p \land q} y$.
	\item $      x \rr{R} y \land \lnot y \rr{R} x$:
	\item $\lnot x \rr{R} y \land       y \rr{R} x$:
	\item $\lnot x \rr{R} y \land \lnot y \rr{R} x$:
		~ ~ ~ These cases are similar,
		using $p_2$, $p_4$, and $p_8$ instead of $p_1$.
	\end{itemize}
\setcounter{enumi}{2}
\item
	Any other boolean connective $c$
	distributes over operation application,
	since
	$\lnot$ and $\land$ do, and $c$ can be obtained as an expression
	over $\lnot$ and $\land$;
	e.g.\ $x R^{p \lor q} y$
	iff $x R^{\lnot (\lnot p \land \lnot q)} y$
	iff $\lnot (\lnot x \ro{R}{p} y \land \lnot x \ro{R}{q} y)$
	iff $x \ro{R}{p} y \lor x \ro{R}{q} y$.
\item
	Follows from~\REF{optn bool 3}:
	$p \subseteq q$
	iff $(p \land q) = p$,
	iff $(R^p \cap R^q) = R^p$
	iff $R^p \subseteq R^q$ .
\qed
\end{enumerate}
}

As an example,
Sen's construction \cite[p.381]{Sen.1969}
of a transitive relation $I$ and a symmetric relation $P$ from a given
quasi-transitive relation $R$, such that $R = I \disjUnion P$,
cf.\ \OEF{Lem.17.3}{26}{quasiTrans 1 3},
can now be paraphrased as $I := R^2$ and $P := R^1$.
Using Lem.~\REF{optn bool}, the proof of disjoint union boils down to
two simple computations:
$I \cup P
= \ro{R}{2} \cup R^1
= R^3
= R$,
and
$I \cap P
= \ro{R}{2} \cap R^1
= R^0
= \set{}$.
See Lem.~\REFF{3-Implications group 5}{gtlu} for another proof using
Lem.~\REF{optn bool}.

\begin{figure}
\begin{center}
$\begin{array}{|cc|c|cc|c|c|}
\hline
x \rr{R} y & y \rr{R} x 
	&    & x \ro{R}{q} y & y \ro{R}{q} x &    & x \ro{(R^q)}{p} y	\\
\hline
0   &  0  &\raa&   q_8\.=0     &   q_8\.=0     &\raa&       p_8	\\
0   &  1  &\raa&   q_4\.=1     &   q_2\.=0     &\raa&       p_2	\\
1   &  0  &\raa&   q_2\.=0     &   q_4\.=1     &\raa&       p_4	\\
1   &  1  &\raa&   q_1\.=1     &   q_1\.=1     &\raa&       p_1	\\
\hline
\end{array}$
\caption{Composition computation example for $q=5$}
\LABEL{Composition computation example}
\end{center}
\end{figure}

\begin{figure}
\begin{center}
$\begin{array}[b]{|c || c*{3}{@{\;}c} | c*{3}{@{\;}c}
	| c*{3}{@{\;}c} | c*{3}{@{\;}c}|}
\hline
_p \backslash ^q
    & 0 & 1 & 2 & 3 & 4 & 5 & 6 & 7 & 8 & 9 & A & B & C & D & E & F \\
\hline
\hline
0   & 0 & 0 & 0 & 0 & 0 & 0 & 0 & 0 & 0 & 0 & 0 & 0 & 0 & 0 & 0 & 0 \\
1   & 0 & 1 & 0 & 1 & 0 & 1 & 6 & 7 & 8 & 9 & 8 & 9 & 8 & 9 & E & F \\
2   & 0 & 0 & 2 & 2 & 4 & 4 & 0 & 0 & 0 & 0 & 2 & 2 & 4 & 4 & 0 & 0 \\
3   & 0 & 1 & 2 & 3 & 4 & 5 & 6 & 7 & 8 & 9 & A & B & C & D & E & F \\
\hline
4   & 0 & 0 & 4 & 4 & 2 & 2 & 0 & 0 & 0 & 0 & 4 & 4 & 2 & 2 & 0 & 0 \\
5   & 0 & 1 & 4 & 5 & 2 & 3 & 6 & 7 & 8 & 9 & C & D & A & B & E & F \\
6   & 0 & 0 & 6 & 6 & 6 & 6 & 0 & 0 & 0 & 0 & 6 & 6 & 6 & 6 & 0 & 0 \\
7   & 0 & 1 & 6 & 7 & 6 & 7 & 6 & 7 & 8 & 9 & E & F & E & F & E & F \\
\hline
8   & F & E & 9 & 8 & 9 & 8 & 9 & 8 & 7 & 6 & 1 & 0 & 1 & 0 & 1 & 0 \\
9   & F & F & 9 & 9 & 9 & 9 & F & F & F & F & 9 & 9 & 9 & 9 & F & F \\
A   & F & E & B & A & D & C & 9 & 8 & 7 & 6 & 3 & 2 & 5 & 4 & 1 & 0 \\
B   & F & F & B & B & D & D & F & F & F & F & B & B & D & D & F & F \\
\hline
C   & F & E & D & C & B & A & 9 & 8 & 7 & 6 & 5 & 4 & 3 & 2 & 1 & 0 \\
D   & F & F & D & D & B & B & F & F & F & F & D & D & B & B & F & F \\
E   & F & E & F & E & F & E & 9 & 8 & 7 & 6 & 7 & 6 & 7 & 6 & 1 & 0 \\
F   & F & F & F & F & F & F & F & F & F & F & F & F & F & F & F & F \\
\hline
\end{array}$
\caption{Composition of unary operations ($p \compose q = ?$)}
\LABEL{Composition p compose q of unary operations}
\end{center}
\end{figure}

\newlength{\colsep}
\setlength{\colsep}{0.6mm}

\begin{figure}
\newcommand{\0}[1]{\begin{array}{@{}c@{}}#1\end{array}}
\small
\begin{center}
$$\begin{array}[b]{ @{} | @{\hspace*{\colsep}} c
	@{\hspace*{\colsep}} || @{\hspace*{\colsep}}
	c*{16}{ @{\hspace*{\colsep}} | @{\hspace*{\colsep}} c }
	@{\hspace*{\colsep}} | @{} }
\hline
_p \backslash ^q
 &0             &1   &2 &3&4 &5&6   &7   &8   &9   &A&B &C&D &E   &F              \\
\hline
\hline
0&\0{0123\\4567}&0246&01&0&01&0&0246&0246&0246&0246&0&08&0&08&0246&\0{0246\\8ACE} \\
\hline
1&              &1357&  &1&  &1&    &    &    &    &8&  &8&  &8ACE&               \\
\hline
2&              &    &23&2&45&4&    &    &    &    &2&2A&4&4C&    &               \\
\hline
3&              &    &  &3&  &5&    &    &    &    &A&  &C&  &    &               \\
\hline
4&              &    &45&4&23&2&    &    &    &    &4&4C&2&2A&    &               \\
\hline
5&              &    &  &5&  &3&    &    &    &    &C&  &A&  &    &               \\
\hline
6&              &    &67&6&67&6&1357&    &    &8ACE&6&6E&6&6E&    &               \\
\hline
7&              &    &  &7&  &7&    &1357&8ACE&    &E&  &E&  &    &               \\
\hline
8&              &    &  &8&  &8&    &8ACE&1357&    &1&  &1&  &    &               \\
\hline
9&              &    &89&9&89&9&8ACE&    &    &1357&9&19&9&19&    &               \\
\hline
A&              &    &  &A&  &C&    &    &    &    &3&  &5&  &    &               \\
\hline
B&              &    &AB&B&CD&D&    &    &    &    &B&3B&D&5D&    &               \\
\hline
C&              &    &  &C&  &A&    &    &    &    &5&  &3&  &    &               \\
\hline
D&              &    &CD&D&AB&B&    &    &    &    &D&5D&B&3B&    &               \\
\hline
E&              &8ACE&  &E&  &E&    &    &    &    &7&  &7&  &1357&               \\
\hline
F&\0{89AB\\CDEF}&9BDF&EF&F&EF&F&9BDF&9BDF&9BDF&9BDF&F&7F&F&7F&9BDF&\0{1357\\9BDF} \\
\hline
\end{array}$$
\caption{Left inverses w.r.t.\ composition ($? \compose q = p$)}
\LABEL{Left inverses w.r.t. composition}
\end{center}
\end{figure}

\begin{figure}
\newcommand{\0}[1]{\begin{array}{@{}c@{}}#1\end{array}}
\small
\begin{center}
$$\begin{array}[b]{ @{} | @{\hspace*{\colsep}} c
	@{\hspace*{\colsep}} || @{\hspace*{\colsep}}
	c*{16}{ @{\hspace*{\colsep}} | @{\hspace*{\colsep}} c }
	@{\hspace*{\colsep}} | @{} }
\hline
_p \backslash ^q
 &0             &1  &2   &3&4   &5&6             &7  &8  &9             &A&B   &C&D   &E  &F             \\
\hline
\hline
0&*             &   &    & &    & &              &   &   &              & &    & &    &   &              \\
\hline
1&024           &135&    & &    & &6             &7  &8AC&9BD           & &    & &    &E  &F             \\
\hline
2&\0{0167\\89EF}&   &23AB& &45CD& &              &   &   &              & &    & &    &   &              \\
\hline
3&0             &1  &2   &3&4   &5&6             &7  &8  &9             &A&B   &C&D   &E  &F             \\
\hline
4&\0{0167\\89EF}&   &45CD& &23AB& &              &   &   &              & &    & &    &   &              \\
\hline
5&0             &1  &4   &5&2   &3&6             &7  &8  &9             &C&D   &A&B   &E  &F             \\
\hline
6&\0{0167\\89EF}&   &    & &    & &\0{2345\\ABCD}&   &   &              & &    & &    &   &              \\
\hline
7&0             &1  &    & &    & &246           &357&8  &9             & &    & &    &ACE&BDF           \\
\hline
8&BDF           &ACE&    & &    & &9             &8  &357&246           & &    & &    &1  &0             \\
\hline
9&              &   &    & &    & &              &   &   &\0{2345\\ABCD}& &    & &    &   &\0{0167\\89EF}\\
\hline
A&F             &E  &B   &A&D   &C&9             &8  &7  &6             &3&2   &5&4   &1  &0             \\
\hline
B&              &   &    & &    & &              &   &   &              & &23AB& &45CD&   &\0{0167\\89EF}\\
\hline
C&F             &E  &D   &C&B   &A&9             &8  &7  &6             &5&4   &3&2   &1  &0             \\
\hline
D&              &   &    & &    & &              &   &   &              & &45CD& &23AB&   &\0{0167\\89EF}\\
\hline
E&F             &E  &    & &    & &9BD           &8AC&7  &6             & &    & &    &135&024           \\
\hline
F&              &   &    & &    & &              &   &   &              & &    & &    &   &*             \\
\hline
\end{array}$$
\caption{Right inverses w.r.t.\ composition ($p \compose ? = q$)}
\LABEL{Right inverses w.r.t. composition}
\end{center}
\end{figure}

\DEFINITION{Operation composition}{%
\LABEL{Operation composition}%
We define operation composition in the usual way by
$R^{p \compose q} = (R^q)^p$.
The set of operations is closed w.r.t.\ composition;
Fig.~\REF{Composition computation example} shows,
by way of an example ($q=5$)
how to compute the
bit representation of $p \compose q$, given that of $p$ and $q$.
\qed
}

Figure~\REF{Composition p compose q of unary operations}
shows a computer-generated composition table.
Observe that operation
$3$ is a neutral element, and the operations $3,5,A,C$
have inverses; in fact, this set is a group w.r.t.\ composition.
Figure~\REF{Left inverses w.r.t. composition} shows the sets of
left inverses w.r.t\ composition;
in line $p$, column $q$, all operations $x$ are listed that satisfy 
$x \compose q = p$.
To save space, we omitted braces and commas.
Similarly, Fig.~\REF{Right inverses w.r.t. composition}
shows the right inverses; line $p$, column $q$ contains all $x$ such
that $p \compose x = q$.

A machine-supported investigation of the algebraic structure w.r.t.\
composition showed nothing interesting.
Of the $65536$ possible operation sets, $461$ are closed w.r.t.\
composition; besides $\set{3, 5, A, C}$ the subsets $\set{ 0 }$,
$\set{F}$,
$\set{ 1, E }$, $\set{ 2, 4}$ $\set{7 ,8}$,  and $\set{B, D}$
are maximal groups, each with a different neutral element.
Each of the $296$ closed subset containing $3$ is, of course, a
monoid;
besides them, the subsets
$\set{ 0, 1, E, F }$,
$\set{ 0, 2, 4 }$,
$\set{ 0, 7, 8, F }$,
and
$\set{ B, D, F }$
are maximal monoids, again each with a different neutral element.
The ancillary file \filename{optnComposition\_GroupsMonoids.txt}
lists each
operation set that is closed w.r.t.\ composition, and indicates
whether it is a monoid, or even a group.

\DEFINITION{Lifted property}{%
If ${\it prop}$ is a property of binary relations,
and $q$ is a unary operation,
then we call $q$-${\it prop}$ a lifted property.
We define that a relation $R$ has the property $q$-${\it prop}$ if
$R^q$ has the property ${\it prop}$.
\qed
}

In this way, unary operations can be lifted from relations to
relation properties.
For example, $8$-transitivity is a synonym for
incomparability-transitivity,
and irreflexivity coincides with $C$-reflexivity.

Starting from the $24$ basic properties from Def.~\REF{def},%
\footnote{
	Unlike in \cite{Burghardt.2018c},
	we this time
	included left-, right-, and two-sided quasi-reflexivity,
	to obtain machine-generated evidence for the laws about
	them; cf.\ Fig.~\REF{Redundant properties}.
	We didn't include co-transitivity (which we only recently
	became aware of) since it obviously can be expressed as
	C-Trans.
}
we obtain
$24 \cdot 16=384$ lifted properties, some of which coincide.
We call two lifted properties equivalent if they agree on every
relation on a sufficiently large,%
\footnote{
	It will turn out that two properties agree on every relation on
	a $7$-element domain iff they agree on every relation of a
	larger domain, see Thm.~\REF{optn equiv}.
}
finite or infinite, domain set.
For example, 3-transitivity is equivalent to 5-transitivity,
due to the self-duality of the definition.
We first show a few simple laws about lifted properties that are needed
later on.
In the next section, we compute equivalence classes of lifted properties.

\LEMMA{Operations and properties}{%
\LABEL{optn univ}%
\begin{enumerate}
\item\LABEL{optn univ 1}%
	Operation 0 always yields the empty relation.
\item\LABEL{optn univ 2}%
	Operation F always yields the universal relation.
\item\LABEL{optn univ 3}%
	Operations B, D, F always yield a connex relation.
\item\LABEL{optn univ 4}%
	Operations 0, 2, 4, 6 always yield an irreflexive relation.
\item\LABEL{optn univ 5}%
	Operations 1, 3, 5, 7 preserve reflexivity and irreflexivity.
\item\LABEL{optn univ 7}%
	Operations 8, A, C, E invert   reflexivity and irreflexivity.
\item\LABEL{optn univ 6}%
	Operations 0, 2, 4 always yield an asymmetric relation.
\item\LABEL{optn univ 14}%
	Operations 9, B, D, F always yield a reflexive relation.
\item\LABEL{optn univ 9}%
	Operations 0, 9, B, D, F always yield a dense relation.
\item\LABEL{optn univ 11}%
	Operations 9, B, D, F always yield a left serial relation.
\item\LABEL{optn univ 13}%
	Operations 0, 9, B, D, F always yield a left quasi-reflexive
	relation.
\item\LABEL{optn univ 15}%
	Operations 0, B, D, F always yield a relation that
	satisfies semi-order property~1.
\item\LABEL{optn univ 16}%
	Operations 0, B, D, F always yield a relation that
	satisfies semi-order property~2.
\item\LABEL{optn univ 17}%
	Operations 0, 1, 6, 7, 8, 9, E, F always yield a symmetric
	relation.
\end{enumerate}
}
\PROOF{
\begin{enumerate}
\item
	Obvious from Def.~\REF{optn def}.
	See \OEF{Exm.74}{48}{} for the
	properties satisfied by the empty relation,
\item
	Obvious from Def.~\REF{optn def}.
	See \OEF{Exm.75}{48}{} for the properties satisfied by the
	universal relation,
\item
	Each listed operation $q$ satisfies
	$7 \compose q = F$,
	hence the symmetric closure of its result relation
	is the universal one.
\item
	$x \ro{R}{q} x$ boils down to $\false$ 
	for $q \in \set{ 0, 2, 4, 6}$.
\item
	For $q \in \set{ 1, 3, 5, 7}$,
	the formula $x \ro{R}{q} x$ boils down to $x \rr{R} x$;
	hence $R^q$ is (ir)reflexive iff $R$ is.
\item
	For $q \in \set{ 8, A, C, E }$,
	the formula $x \ro{R}{q} x$ boils down to $\lnot x \rr{R} x$;
	hence $R^q$ is reflexive iff $R$ is irreflexive, and
	$R^q$ is irreflexive iff $R$ is reflexive.
\item
	Each listed operation $q$ satisfies
	$q \cap (5 \compose q) = 0$,
	hence its result relation is disjoint from its own converse.
\item
	$x \ro{R}{9} x$ iff 
	$x \rr{R} x \lxor \lnot x \rr{R} x$ iff $\true$.
	Similarly, $x \ro{R}{q} x$ boils down to $\true$
	for $q \in \set{ B, D, F}$.
\item
	The empty relation $R^0$ is dense by
	\OEF{Exm.74.11}{48}{empty 9}.
	Operations 9, B, D, F always yield a reflexive relation
	by~\REF{optn univ 14},
	which is dense by \OEF{Lem.48.1}{38}{den 1 1}.
\item
	Operations 9, B, D, F always yield a reflexive relation
	by~\REF{optn univ 14},
	which is left-serial by \OEF{Lem.54}{40}{ser 1}.
\item
	The empty relation $R^0$ is quasi-reflexive by
	\OEF{Exm.74.1}{48}{empty 1}, hence left quasi-reflexive in
	particular.
	Operations 9, B, D, F always yield a reflexive relation
	by~\REF{optn univ 14},
	which is quasi-reflexive by \OEF{Lem.9}{23}{quasiRefl 2},
	and hence left quasi-reflexive.
\item
	The empty relation $R^0$ satisfies semi-order property~1
	by \OEF{Exm.74.9}{48}{empty 8}.
	Operations B, D, F always yield a connex relation by
	Lem.~\REFF{optn univ}{3},
	which satisfies semi-order property~1 by
	\OEF{Lem.66}{45}{semiOrd1 1a}.
\item
	The empty relation $R^0$ satisfies semi-order property~2
	by \OEF{Exm.74.9}{48}{empty 8}.
	Operations B, D, F always yield a connex relation by
	Lem.~\REFF{optn univ}{3},
	which satisfies semi-order property~2 by
	\OEF{Lem.66}{45}{semiOrd1 1a}.
\item
	Each listed operation $q$ satisfies $q = 5 \compose q$,
	hence its result relation agrees with its own converse.
\qed
\end{enumerate}
}

\LEMMA{Unsatisfiable lifted properties}{%
\LABEL{optn unsat}%
On a domain set $X$ with $\geq 1$ elements,
no relation has any of the following properties:
\begin{enumerate}
\item\LABEL{optn unsat 1}%
	9-, B-, D-, F-ASym
\item\LABEL{optn unsat 3}%
	9-, B-, D-, F-AntiTrans
\item\LABEL{optn unsat 8}%
	0-LfSerial
\item\LABEL{optn unsat 10}%
	0-, 2-, 4-, 6-Refl
\end{enumerate}

On a domain set with $\geq2$ elements,
no relation has one of the following properties:
\begin{enumerate}
\setcounter{enumi}{4}%
\item\LABEL{optn unsat 2}%
	F-AntiSym
\end{enumerate}

On a domain set with $\geq 4$ elements,
no relation has one of the following properties:
\begin{enumerate}
\setcounter{enumi}{5}%
\item\LABEL{optn unsat 9}%
	B-, D-, F-LfUnique
\end{enumerate}
}
\PROOF{
On a nonempty domain set:
\begin{enumerate}
\item
	By Lem.~\REFF{optn univ}{14},
	operations 9, B, D, F yield a reflexive relation,
	which cannot be asymmetric by \OEF{Lem.10}{23}{irrefl 2}.
\item
	By Lem.~\REFF{optn univ}{14},
	operations 9, B, D, F yield a reflexive relation,
	which cannot be anti-transitive
	by \OEF{Lem.10}{23}{irrefl 2}.
\item
	Operation 0 always yields the empty relation which isn't
	left-serial on a non-empty domain.
\item
	By Lem.~\REFF{optn univ}{4},
	operations 0, 2, 4, 6 yield an irreflexive relation,
	which cannot be reflexive on a nonempty domain.
\end{enumerate}

On a domain set with $\geq 2$ elements:
\begin{enumerate}
\setcounter{enumi}{4}%
\item
	For $x \neq y$,
	we have $x \ro{R}{F} y \land y \ro{R}{F} x$,
	contradicting anti-symmetry.
\end{enumerate}

On a domain set with $\geq 4$ elements:
\begin{enumerate}
\setcounter{enumi}{5}%
\item
	By Lem.~\REFF{optn univ}{3},
	operations B, D, F yield a connex relation,
	which cannot be left-unique
	on a set of $4$ or more elements
	by \OEF{Lem.51}{39}{conn 2}.
\qed
\end{enumerate}
}

\LEMMA{Operation names}{%
\LABEL{optn names lem}%
\begin{enumerate}
\item\LABEL{optn names lem 1}%
	The symmetric kernel of a relation $R$ is obtained as $R^1$.
\item\LABEL{optn names lem 2}%
	The symmetric closure of a relation $R$ is obtained as $R^7$.
\item\LABEL{optn names lem 3}%
	The asymmetric kernel of a relation $R$ is obtained as $R^2$.
	As discussed in Def.~\REFF{optn names def}{3},
	it is an asymmetric subset of $R$,
	but need not be a maximal one.
\end{enumerate}
}
\PROOF{
\begin{enumerate}
\item
	We show that $\ro{R}{1}$ is the largest symmetric sub-relation
	of $R$.
	First, $R^1$ is symmetric by Lem.~\REFF{optn univ}{17},
	and a subset of $R$ by Lem.~\REFF{optn bool}{4}.
	If $R' \subseteq R$ is a symmetric relation,
	then $x \rr{R'} y$ implies $y \rr{R'} x$ by symmetry,
	hence $x \rr{R} y \land y \rr{R} x$ 
	by the subset property,
	hence $x \ro{R}{1} y$ by definition.
\item
	We show that $\ro{R}{7}$ is the smallest symmetric
	super-relation of $R$.
	First, $R^7$ is symmetric by Lem.~\REFF{optn univ}{17},
	and a superset of $R$ by Lem.~\REFF{optn bool}{4}.
	Let $R' \supseteq R$ be a symmetric relation.
	If $x \ro{R}{7} y$, then by definition
	$x \rr{R} y$ or $y \rr{R} x$.
	In the former case, we have $x \rr{R'} y$ per superset,
	in the latter, we additionally use the symmetry of $R'$.
\item
	We show that $\ro{R}{2}$ equals the intersection of all
	maximal asymmetric sub-relations of $R$.
	\begin{itemize}
	\item[``$\subseteq$'':]
		Let $R'$ be a maximal asymmetric subset of $R$, we show
		$\ro{R}{2} \subseteq R'$:
		Let $x \ro{R}{2} y$,
		then $x \rr{R} y \land \lnot y \rr{R} x$,
		hence $\lnot y \rr{R'} x$.
		%
		If $x \rr{R'} y$ did not hold, 
		then $R' \cup \set{ \tpl{x,y} }$
		was a larger, but still asymmetric subset of $R$.
	\item[``$\supseteq$'':]
		First, $\set{}$ is an asymmetric subset of $R$,
		hence%
		\footnote{
			We use Zorn's lemma here.
		}
		a maximal one can also be found; we call it $R_0$.
		Now let $x \rr{R'} y$
		for every maximal asymmetric subset $R'$
		of $R$, we show $x \ro{R}{2} y$.
		We have in particular $x R_0 y$,
		hence $x \rr{R} y$,
		Assume for contradiction $y \rr{R} x$ does also hold.
		Then
		$R_1 = 
		(R_0 \setminus \set{ \tpl{x,y} }) \cup \set{ \tpl{y,x} }$
		is another maximal asymmetric subset of $R$, 
		which, however,
		does not satisfy $x R_1 y$, contrary to our assumption.
	\end{itemize}
	Note that $R^2$ is asymmetric by Lem.~\REFF{optn univ}{6},
	and a subset of $R$ by Lem.~\REFF{optn bool}{4}.
	In contrast, $\ro{R}{4}$ is different from the asymmetric kernel
	since it is not a subset of $R$.
\qed
\end{enumerate}
}

\clearpage

\section{Equivalent lifted properties}
\LABEL{Equivalent lifted properties}

In this section, we investigate extensional equality of lifted
properties.
We show that, starting from the basic property set from Def.~\REF{def},
a of total $81$ different equivalence classes exist
(Thm.~\REF{optn equiv}).

\definecolor{fgPartSet}	{rgb}{0.00,0.00,0.70}
\definecolor{bgPartSet}	{rgb}{0.90,0.90,0.99}
\definecolor{fgPart}	{rgb}{0.00,0.70,0.00}
\definecolor{bgPart}	{rgb}{0.90,0.99,0.90}
\definecolor{fgRel}	{rgb}{0.70,0.00,0.00}
\definecolor{bgRel}	{rgb}{0.99,0.90,0.90}

\begin{figure}
\begin{center}
\begin{picture}(150,45)
\thicklines
\textcolor{bgPartSet}{\put(0,0){\makebox(0,0)[bl]{\rule{60mm}{40mm}}}}%
\textcolor{fgPartSet}{\put(0,0){\line(1,0){60}}}%
\textcolor{fgPartSet}{\put(0,0){\line(0,1){40}}}%
\textcolor{fgPartSet}{\put(0,40){\line(1,0){60}}}%
\textcolor{fgPartSet}{\put(60,0){\line(0,1){40}}}%
\textcolor{bgPart}{\put(5,7){\framebox(9,4){\rule{9mm}{4mm}}}}%
\textcolor{bgPart}{\put(15,20){\framebox(9,3){\rule{9mm}{3mm}}}}%
\textcolor{bgPart}{\put(25,13){\framebox(9,5){\rule{9mm}{5mm}}}}%
\textcolor{bgPart}{\put(33,29){\framebox(9,6){\rule{9mm}{6mm}}}}%
\textcolor{bgPart}{\put(7,25){\framebox(9,8){\rule{9mm}{8mm}}}}%
\textcolor{bgPart}{\put(49,22.5){\framebox(9,5){\rule{9mm}{5mm}}}}%
\textcolor{bgPart}{\put(49,7){\framebox(9,2){\rule{9mm}{2mm}}}}%
\textcolor{bgPart}{\put(49,3){\framebox(9,2){\rule{9mm}{2mm}}}}%
\textcolor{fgPart}{\put(5,7){\framebox(9,4){}}}%
\textcolor{fgPart}{\put(15,20){\framebox(9,3){}}}%
\textcolor{fgPart}{\put(25,13){\framebox(9,5){}}}%
\textcolor{fgPart}{\put(33,29){\framebox(9,6){}}}%
\textcolor{fgPart}{\put(7,25){\framebox(9,8){}}}%
\textcolor{fgPart}{\put(49,22.5){\framebox(9,5){}}}%
\textcolor{fgPart}{\put(49,7){\framebox(9,2){}}}%
\textcolor{fgPart}{\put(49,3){\framebox(9,2){}}}%
\textcolor{bgPart}{\put(88,18){\makebox(0,0)[bl]{\rule{22mm}{12mm}}}}%
\textcolor{fgPart}{\put(88,18){\line(1,0){22}}}%
\textcolor{fgPart}{\put(88,24){\line(1,0){22}}}%
\textcolor{fgPart}{\put(88,18){\line(0,1){12}}}%
\textcolor{fgPart}{\put(88,30){\line(1,0){22}}}%
\textcolor{fgPart}{\put(110,18){\line(0,1){12}}}%
\textcolor{fgPart}{\put(89,19){\makebox(0,0)[bl]{\small C-LfSerial}}}%
\textcolor{fgPart}{\put(89,25){\makebox(0,0)[bl]{\small 3-LfUnique}}}%
\textcolor{bgRel}{\put(128,5){\makebox(0,0)[bl]{\rule{22mm}{19mm}}}}%
\textcolor{fgRel}{\put(135,5){\line(1,0){15}}}%
\textcolor{fgRel}{\put(135,10){\line(1,0){15}}}%
\textcolor{fgRel}{\put(135,15){\line(1,0){15}}}%
\textcolor{fgRel}{\put(128,20){\line(1,0){22}}}%
\textcolor{fgRel}{\put(135,5){\line(0,1){19}}}%
\textcolor{fgRel}{\put(140,5){\line(0,1){15}}}%
\textcolor{fgRel}{\put(145,5){\line(0,1){15}}}%
\textcolor{fgRel}{\put(150,5){\line(0,1){15}}}%
\textcolor{fgRel}{\put(137.5,17.5){\makebox(0,0){\sf X}}}%
\textcolor{fgRel}{\put(142.5,17.5){\makebox(0,0){\sf X}}}%
\textcolor{fgRel}{\put(150,21){\makebox(0,0)[br]{$\scriptstyle x$}}}%
\textcolor{fgRel}{\put(134,5){\makebox(0,0)[br]{$\scriptstyle y$}}}%
\textcolor{fgRel}{\put(134,21){\makebox(0,0)[br]{$\scriptstyle x\rr{R}y$}}}%
\thinlines
\put(60,25){\vector(1,0){28}}
\put(110,15){\oval(6,16)[r]}
\put(110,15){\oval(10,20)[r]}
\put(110,7){\vector(-1,0){50}}
\put(110,5){\vector(-1,0){50}}
\put(110,8){\makebox(0,0)[br]{$\true$}}
\put(110,4){\makebox(0,0)[tr]{$\false$}}
\put(30,45){\makebox(0,0)[t]{Set of partitions}}
\put(99,45){\makebox(0,0)[t]{Current partition}}
\put(139,45){\makebox(0,0)[t]{Current relation}}
\end{picture}
\caption{Partition split algorithm}
\LABEL{Partition split algorithm}
\end{center}
\end{figure}
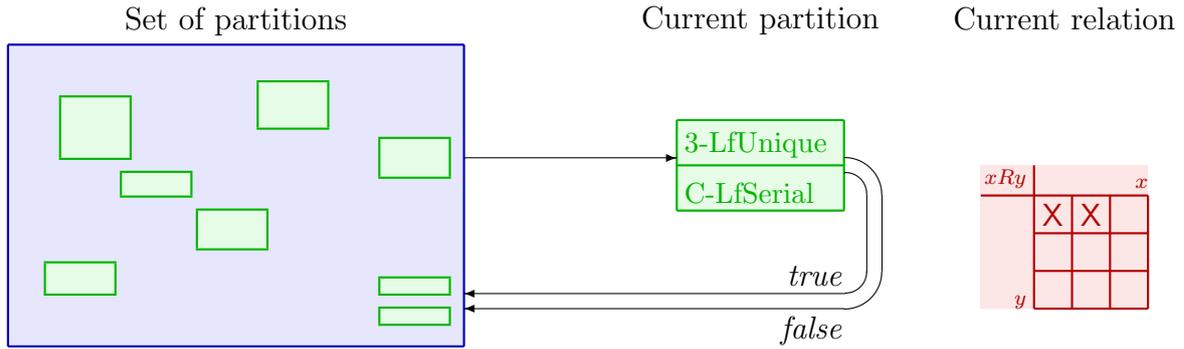

To obtain an approximation of equivalence,
we implemented a partition-refinement
routine: initially, all $384$ lifted properties are in one partition;
cycling though each relation $R$ on a small finite set, we split each
partition according to the property behavior on $R$.
Figure~\REF{Partition split algorithm} shows a snapshot of
the algorithm, when it is about to split
the partition $\set{ \mbox{3-LfUnique}, \mbox{C-LfSerial} }$
into two singleton sets, since the current relation is not left unique,
but its complement is left serial.
Running the routine on the relations over a $5$-element set, we ended
up with $80$ partitions;
for a $6$-element set, the very same partitions were obtained;
checking on a $7$-element set would take far too long.

The raw output about the computed partitions
is available in the ancillary file
\filename{lpEqns.txt};
each partition is represented as an equation chain like
\texttt{5-Trans = 3-Trans}.
Most of them are easily seen to be in fact equivalences,
we give the formal proofs in the following.
Some equivalences allow us to define one property in terms of another,
see Fig.~\REF{Redundant properties} and
Lem.~\REF{optn redundant} below.

Figure~\REF{Computed partitions of lifted properties}
shows the partitions, omitting the redundant properties for brevity.
Each box in the left part denotes one partition.
The two topmost boxes denote the trivial partitions, viz.\ that of
properties that hold for no relation at all (left) and for every
relation (right); only two of the $46$ and $17$ members, are
shown, respectively.
These two partitions together
contain all lifted properties starting with 0- or F-,%
\footnote{
	The property 0-{\it prop} applies to a relation $R$
	iff {\it prop} applies to the empty relation, independent of
	$R$.
	Hence, 0-{\it prop} is either true on all relations or false
	on all.
	Similar for F-{\it prop}.
}
along with many others, like 2-Refl and 7-Sym.

Partitions containing only different operations prepended to a single
basic property are called pure partitions.
They are not shown as a box, but listed in the right part
of Fig.~\REF{Computed partitions of lifted properties}.
For example, the first line is short notation for the raw output
chains
\texttt{1-ASym = 3-ASym = 5-ASym} and
\texttt{8-ASym = A-ASym = C-ASym}.
Singleton partitions are a subclass of pure partitions; they are
listed separately.
For example, the first line is short notation for the trivial
raw output chains \texttt{7-ASym} and \texttt{E-ASym}.%
\footnote{
	Note that the actual chains in file
	\filename{lpEqns.txt} are longer than shown here,
	since they still contain redundant properties.
}

By construction, the computed partitioning is a coarsening of the
proper extensional equality partitioning.
To obtain the latter, it sufficient to split the largest box partition
along the dotted line, as we will show
in Sect.~\REF{Proof of equivalence classes},
in particular in Thm.~\REF{optn equiv}.
We thus arrive at a total of $81$ equivalence classes of lifted properties.
We chose a default representation for each class, they are shown in
Figure~\REF{Default representations of lifted properties}.%
\footnote{
	The figure shows, for each unary operation $q$ and each basic
	property
	${\it prop}$, the default representation to which $q$-{\it prop}
	is extensionally equal.
	If $q$-{\it prop} is its own default representation, 
	it is highlighted in
	blue and cyan, for a non-singleton and a singleton equivalence
	class, respectively.
	The trivial partitions are denoted by $+$ and $-$.
}
Considering all of them
in our naive Quine-McCluskey implementation from
\cite[Sect.6]{Burghardt.2018c}
would require more
than $2^{81}$ bits, i.e.\
$274877$ millions of Tera bytes.
This still renders that approach infeasible.

However, it could be used to investigate small subsets of properties,
looking for law suggestions including complement, converse, etc.
For example, to investigate the connections between quasi-transitivity
and semi-orders, one may restrict oneself to QuasiTrans, SemiOrd1,
SemiOrd2, and ASym.
From Fig.~\REF{Default representations of lifted properties},
it can seen that $21$ distinct properties would be needed.%
\footnote{
	 +,    -,     as, 7-as, C-as, E-as,   s1, 1-s1, 2-s1, 7-s1,
	C-s1,   s2, 2-s2, C-s2,   sy, 1-tr, 2-tr, 7-tr, 8-tr, 9-tr,
	E-tr.
}
However, we didn't yet perform such an investigation.

Instead, we checked all possible implications between default
representations up
to a length of $3$, i.e.\ with at most $2$ antecedents.
This is presented in
Sect.~\REF{Implications between lifted properties}.

\begin{figure}
\begin{center}
\begin{tabular}{rcl}
RgSerial	& $\lra$ & 5-LfSerial	\\
RgEucl		& $\lra$ & 5-LfEucl	\\
RgUnique	& $\lra$ & 5-LfUnique	\\
Univ		& $\lra$ & E-ASym	\\
Empty		& $\lra$ & 7-ASym	\\
Irrefl		& $\lra$ & C-Refl	\\
Connex		& $\lra$ & C-ASym	\\
SemiConnex	& $\lra$ & C-AntiSym	\\
IncTrans	& $\lra$ & 8-Trans	\\
CoRefl		& $\lra$ & 7-AntiSym	\\
QuasiTrans	& $\lra$ & 2-Trans	\\
RgQuasiRefl	& $\lra$ & 5-LfQuasiRefl\\
QuasiRefl	& $\lra$ & 7-LfQuasiRefl\\
\end{tabular}
\caption{Redundant properties}
\LABEL{Redundant properties}
\end{center}
\end{figure}

\newcommand{\oEqc}[1]{%
	\begin{tabular}[t]{|l|}
	\hline
	#1	\\
	\hline
	\end{tabular}
	\\
}

\begin{figure}
\begin{tabular}[t]{@{}l@{}}
\oEqc{F-AntiSym\\0-Refl \ldots}
\oEqc{D-Sym\\D-LfEucl\\C-Sym\\B-Sym\\B-LfEucl\\A-Sym\\9-SemiOrd2\\9-SemiOrd1\\6-AntiSym\\6-LfQuasiRefl\\6-SemiOrd1\\6-Trans\\6-ASym\\6-LfEucl\\5-Sym\\4-LfQuasiRefl\\4-Sym\\4-LfEucl\\3-Sym\\2-LfQuasiRefl\\2-Sym\\2-LfEucl\\[-1ex]..................\\4-Dense\\2-Dense}
\end{tabular}
\hfill
\begin{tabular}[t]{@{}l@{}}
\oEqc{F-LfSerial\\0-Sym \ldots}
\oEqc{D-AntiSym\\B-AntiSym\\9-AntiSym\\9-LfUnique}
\\
\\
\oEqc{1-LfEucl\\E-SemiOrd2\\1-Trans}
\oEqc{7-LfEucl\\8-SemiOrd2\\7-Trans}
\oEqc{8-LfEucl\\7-SemiOrd2\\8-Trans}
\oEqc{9-LfEucl\\6-SemiOrd2\\9-Trans}
\oEqc{E-LfEucl\\1-SemiOrd2\\E-Trans}
\end{tabular}
\hfill
\begin{tabular}[t]{@{}l@{}}
\\
\bf Pure Partitions
\\
\\
\begin{tabular}[t]{@{}l*{5}{@{$\;$}l}@{}}
ASym:      & \multicolumn{2}{@{}l}{$1 \.\lra 3 \.\lra 5$,}
	& \multicolumn{2}{@{}l}{$8 \.\lra A \.\lra C$}	\\
AntiSym:   & \multicolumn{2}{@{}l}{$1 \.\lra 3 \.\lra 5$,}
	& \multicolumn{2}{@{}l}{$8 \.\lra A \.\lra C$}	\\
Refl:      & \multicolumn{2}{@{}l}{$1 \.\lra 3 \.\lra 5 \.\lra 7$,}
	& \multicolumn{2}{@{}l}{$8 \.\lra A \.\lra C \.\lra E$}	\\
\\
SemiOrd1: &$1\.\lra E$,&$2\.\lra 4$,&$3\.\lra 5$,&$7\.\lra 8$,&$A\.\lra C$ \\
SemiOrd2: &            &$2\.\lra 4$,&$3\.\lra 5$,&            &$A\.\lra C$ \\
AntiTrans:&            &$2\.\lra 4$,&$3\.\lra 5$,&            &$A\.\lra C$ \\
Trans:    &            &$2\.\lra 4$,&$3\.\lra 5$,&$B\.\lra D$,&$A\.\lra C$ \\
Dense:    &            &            &$3\.\lra 5$,&            &$A\.\lra C$ \\
\end{tabular}
\\
\\
\\
\\
\\
\\
\bf Singleton Partitions
\\
\\
\begin{tabular}[t]{@{}l*{11}{@{$\;\;$}l}@{}}
ASym:        &    &    &    &    &    &    & 7, &    &    &    & E	\\
AntiSym:     &    &    &    &    &    &    & 7, &    &    &    & E	\\
AntiTrans:   & 1, &    &    &    &    & 6, & 7, & 8, &    &    & E	\\
Dense:       & 1, &    &    &    &    & 6, & 7, & 8, &    &    & E	\\
LfEucl:      &    &    & 3, &    & 5, &    &    &    & A, & C  &  	\\
LfSerial:    & 1, & 2, & 3, & 4, & 5, & 6, & 7, & 8, & A, & C, & E	\\
LfUnique:    & 1, & 2, & 3, & 4, & 5, & 6, & 7, & 8, & A, & C, & E	\\
LfQuasiRefl: & 1, &    & 3, &    & 5, &    & 7, & 8, & A, & C, & E	\\
\end{tabular}
\\
\end{tabular}
\begin{center}
\caption{Computed partitions of lifted properties}
\LABEL{Computed partitions of lifted properties}
\end{center}
\end{figure}

\setlength{\colsep}{0.22mm}

\definecolor{coNfNone}		{rgb}{0.00,0.00,0.00}
\definecolor{coNfAll}		{rgb}{0.00,0.00,0.00}
\definecolor{coNfSelf}		{rgb}{0.00,0.80,0.99}
\definecolor{coNfSELF}		{rgb}{0.00,0.00,0.99}
\definecolor{coNfRedundant}	{rgb}{0.50,0.50,0.50}

\begin{figure}
\newcommand{\3}{}
\newcommand{\0}{\textcolor{coNfNone}{$-$}}
\newcommand{\1}{\textcolor{coNfAll}{$+$}}
\newcommand{\2}{\3sy}
\newcommand{\4}{\color{coNfRedundant}}
\newcommand{\8}{\color{coNfSELF}}
\newcommand{\9}{\color{coNfSelf}}
\begin{center}
\begin{tabular}[t]{@{}|@{\hspace*{\colsep}}c@{\hspace*{\colsep}}|@{\hspace*{\colsep}}l*{16}{@{\hspace*{\colsep}}|@{\hspace*{\colsep}}c}@{\hspace*{\colsep}}|@{}}
\hline
  &               &0 &  1   &  2   &  3   &  4   &  5   &  6   &  7   &  8   &  9   &  A   &  B   &  C   &  D   &  E   &F \\\hline

em&\4 Empty       &\1&  \3as&  \3sy&  7-as&  \3sy&  7-as&  \3sy&  7-as&  C-as&  \0  &  E-as&  \0  &  E-as&  \0  &  E-as&\0\\\hline
un&\4 Univ        &\0&  E-as&  \0  &  E-as&  \0  &  E-as&  \0  &  C-as&  7-as& \3sy &  7-as&  \3sy&  7-as&  \3sy&  3-as&\1\\\hline

rf&   Refl        &\0&  \3rf&  \0  &\8\3rf&  \0  &  \3rf&  \0  &  \3rf&  C-rf&  \1  &  C-rf&  \1  &\8C-rf&  \1  &  C-rf&\1\\\hline
  &\4 Irrefl      &\1&  C-rf&  \1  &  C-rf&  \1  &  C-rf&  \1  &  C-rf&  \3rf&  \0  &  \3rf&  \0  &  \3rf&  \0  &  \3rf&\0\\\hline
  &\4 CoRefl      &\1&  \3an&  \2  &  7-an&  \2  &  7-an&  \2  &  7-an&  C-an&  9-an&  E-an&  \0  &  E-an&  \0  &  E-an&\0\\\hline
lq&   LfQuasiRefl &\1&\91-lq&  \2  &\9\3lq&  \2  &\95-lq&  \2  &\97-lq&\98-lq&  \1  &\9A-lq&  \1  &\9C-lq&  \1  &\9E-lq&\1\\\hline
  &\4 RgQuasiRefl &\1&  1-lq&  \2  &  5-lq&  \2  &  \3lq&  \2  &  7-lq&  8-lq&  \1  &  C-lq&  \1  &  A-lq&  \1  &  E-lq&\1\\\hline
  &\4 QuasiRefl   &\1&  1-lq&  \2  &  7-lq&  \2  &  7-lq&  \2  &  7-lq&  8-lq&  \1  &  E-lq&  \1  &  E-lq&  \1  &  E-lq&\1\\\hline

sy&   Sym         &\1&  \1  &  \2  &\8\2  &  \2  &  \2  &  \1  &  \1  &  \1  &  \1  &  \2  &  \2  &  \2  &  \2  &  \1  &\1\\\hline
as&   ASym        &\1&  \3as&  \1  &\8\3as&  \1  &  \3as&  \2  &\97-as&  C-as&  \0  &  C-as&  \0  &\8C-as&  \0  &\9E-as&\0\\\hline
an&   AntiSym     &\1&  \3an&  \1  &\8\3an&  \1  &  \3an&  \2  &\97-an&  C-an&\89-an&  C-an&  9-an&\8C-an&  9-an&\9E-an&\0\\\hline

  &\4 SemiConnex  &\0&  E-an&  9-an&  C-an&  9-an&  C-an&  9-an&  C-an&  7-an&  \2  &  \3an&  \1  &  \3an&  \1  &  \3an&\1\\\hline
  &\4 Connex      &\0&  E-as&  \0  &  C-as&  \0  &  C-as&  \0  &  C-as&  7-as&  \2  &  \3as&  \1  &  \3as&  \1  &  \3as&\1\\\hline

tr&   Trans       &\1&\81-tr&\82-tr&\8\3tr&  2-tr&  \3tr&  \2  &\87-tr&\88-tr&\89-tr&  C-tr&  D-tr&\8C-tr&\8D-tr&\8E-tr&\1\\\hline
at&   AntiTrans   &\1&\91-at&\82-at&\8\3at&  2-at&  \3at&\96-at&\97-at&\98-at&  \0  &  C-at&  \0  &\8C-at&  \0  &\9E-at&\0\\\hline
  &\4 QuasiTrans  &\1&  \1  &  2-tr&  2-tr&  2-tr&  2-tr&  \1  &  \1  &  \1  &  \1  &  2-tr&  2-tr&  2-tr&  2-tr&  \1  &\1\\\hline
  &\4 RgEucl      &\1&  1-tr&  \2  &  5-le&  \2  &  \3le&  \2  &  7-tr&  8-tr&  9-tr&  C-le&  \2  &  A-le&  \2  &  E-tr&\1\\\hline
le&   LfEucl      &\1&  1-tr&  \2  &\9\3le&  \2  &\95-le&  \2  &  7-tr&  8-tr&  9-tr&\9A-le&  \2  &\9C-le&  \2  &  E-tr&\1\\\hline

s1&   SemiOrd1    &\1&\81-s1&\82-s1&\8\3s1&  2-s1&  \3s1&  \2  &\87-s1&  7-s1&  \2  &  C-s1&  \1  &\8C-s1&  \1  &  1-s1&\1\\\hline
s2&   SemiOrd2    &\1&  E-tr&\82-s2&\8\3s2&  2-s2&  \3s2&  9-tr&  8-tr&  7-tr&  \2  &  C-s2&  \1  &\8C-s2&  \1  &  1-tr&\1\\\hline

  &\4 RgSerial    &\0&  1-ls&  4-ls&  5-ls&  2-ls&  \3ls&  6-ls&  7-ls&  8-ls&  \1  &  C-ls&  \1  &  A-ls&  \1  &  E-ls&\1\\\hline
ls&   LfSerial    &\0&\91-ls&\92-ls&\9\3ls&\94-ls&\95-ls&\96-ls&\97-ls&\98-ls&  \1  &\9A-ls&  \1  &\9C-ls&  \1  &\9E-ls&\1\\\hline

de&   Dense       &\1&\91-de&\82-de&\8\3de&  2-de&  \3de&\96-de&\97-de&\98-de&  \1  &  C-de&  \1  &\8C-de&  \1  &\9E-de&\1\\\hline
  &\4 IncTrans    &\1&  E-tr&  9-tr&  8-tr&  9-tr&  8-tr&  9-tr&  8-tr&  7-tr&  \2  &  1-tr&  \1  &  1-tr&  \1  &  1-tr&\1\\\hline

lu&   LfUnique    &\1&\91-lu&\92-lu&\9\3lu&\94-lu&\95-lu&\96-lu&\97-lu&\98-lu&  9-an&\9A-lu&  \0  &\9C-lu&  \0  &\9E-lu&\0\\\hline
  &\4 RgUnique    &\1&  1-lu&  4-lu&  5-lu&  2-lu&  \3lu&  6-lu&  7-lu&  8-lu&  9-an&  C-lu&  \0  &  A-lu&  \0  &  E-lu&\0\\\hline

\end{tabular}
\caption{Default representations of lifted properties}
\LABEL{Default representations of lifted properties}
\end{center}
\end{figure}

\clearpage

\subsection{Proof of equivalence classes}
\LABEL{Proof of equivalence classes}

\LEMMA{}{%
\LABEL{optn 2}%
\begin{enumerate}
\item\LABEL{optn 2 1}%
	If $R$ is symmetric, 
	then $R^0 = \ro{R}{2} = \ro{R}{4} = R^6$ are the empty
	relation, 
	and $R^9 = \ro{R}{B} = \ro{R}{D} = R^F$ are the universal relation.
\item\LABEL{optn 2 3}%
	If $R$ is asymmetric, then $R^0 = R^1$ are the empty relation,
	and $R^E = R^F$ are the universal relation.
\end{enumerate}
}
\PROOF{
\begin{enumerate}
\item
	Symmetry implies $R = R^7$,
	hence $R^2 = R^{(2 \compose 7)} = R^0$,
	similar for the other operations.
\item
	Asymmetry implies $R = R^2$,
	hence $R^1 = R^{(1 \compose 2)} = R^0$,
	similar for $E$.
\qed
\end{enumerate}
}

\LEMMA{Redundant properties}{%
\LABEL{optn redundant}%
\begin{enumerate}
\item\LABEL{optn redundant 1}%
	RgSerial	$\Lra$ 5-LfSerial
\item\LABEL{optn redundant 2}%
	RgEucl	$\Lra$ 5-LfEucl
\item\LABEL{optn redundant 3}%
	RgUnique	$\Lra$ 5-LfUnique
\item\LABEL{optn redundant 12}%
	RgQuasiRefl	$\Lra$ 5-LfQuasiRefl
\item\LABEL{optn redundant 13}%
	QuasiRefl	$\Lra$ 7-LfQuasiRefl
\item\LABEL{optn redundant 4}%
	Univ	$\Lra$ E-ASym
\item\LABEL{optn redundant 5}%
	Empty	$\Lra$ 7-ASym
\item\LABEL{optn redundant 6}%
	Irrefl	$\Lra$ C-Refl
\item\LABEL{optn redundant 7}%
	Connex	$\Lra$ C-ASym
\item\LABEL{optn redundant 8}%
	SemiConnex	$\Lra$ C-AntiSym
\item\LABEL{optn redundant 9}%
	IncTrans	$\Lra$ 8-Trans
\item\LABEL{optn redundant 10}%
	CoRefl	$\Lra$ 7-AntiSym
\item\LABEL{optn redundant 11}%
	QuasiTrans	$\Lra$ 2-Trans
\end{enumerate}
In Fig.~\REF{Redundant properties},
we summarize these redundancies.
}
\PROOF{
\begin{enumerate}
\item[\REF{optn redundant 1}--%
	\REF{optn redundant 12}]
	From the definitions in Def.~\REF{def} is is obvious that each
	``Rg'' property is the converse of the corresponding ``Lf''
	property.
\setcounter{enumi}{4}%
\item
	First, observe that $x \rr{R} x$ 
	iff $x \ro{R}{7} x$ by definition.
	\\
	``$\Ra$'':
	If $x \ro{R}{7} y$,
	then $x \rr{R} y \lor y \rr{R} x$;
	both cases imply $x \rr{R} x$ by quasi-reflexivity.
	\\
	``$\La$'':
	If $x \rr{R} y$,
	then $x \rr{R} y \lor y \rr{R} x$ by weakening,
	hence $x \ro{R}{7} y$, and similarly $y \ro{R}{7} x$,
	hence $x \ro{R}{7} x$ and $y \ro{R}{7} y$.
\item
	``$\Ra$'':
	Let $x \ro{R}{E} y$,
	then $\lnot x \rr{R} y \lor \lnot y \rr{R} x$;
	both cases contradict universality.
	\\
	``$\La$'':
	Let $R^E$ be asymmetric;
	it is also symmetric
	by Lem.~\REFF{optn univ}{17},
	hence empty by \OEF{Lem.16}{25}{sym 2},
	hence 
	$\lnot x \rr{R} y \lor \lnot y \rr{R} x$ never holds,
	hence $R$ is universal.
\item
	``$\Ra$'':
	$x \ro{R}{7} y$ would imply $x \rr{R} y \lor y \rr{R} x$,
	contradicting $R$'s emptiness.
	\\
	``$\La$'':
	$x \rr{R} y$ 
	would imply $x \rr{R} y \lor y \rr{R} x$ by weakening,
	hence $x \ro{R}{7} y$ and $y \ro{R}{7} x$,
	contradicting asymmetry.
\item
	Since $\lnot x \rr{R} x$ iff $x \ro{R}{C} x$, for all $x$.
\item
	``$\Ra$'':
	Let $x \ro{R}{C} y$,
	then $\lnot x \rr{R} y$ by definition,
	hence $y \rr{R} x$ by connexity,
	hence $\lnot y \ro{R}{C} x$.
	\\
	``$\La$'':
	If $\lnot x \rr{R} y$,
	then $x \ro{R}{C} y$,
	hence $\lnot y \ro{R}{C} x$ by asymmetry,
	hence $y \rr{R} x$.
\item
	Similar to~\REF{optn redundant 7}:
	\\
	``$\Ra$'':
	Let $x \ro{R}{C} y$,
	then $\lnot x \rr{R} y$ by definition,
	then $x=y$ or hence $y \rr{R} x$ by semi-connexity,
	hence $x=y$ or $\lnot y \ro{R}{C} x$.
	\\
	``$\La$'':
	If $\lnot x \rr{R} y$,
	then $x \ro{R}{C} y$,
	hence $x=y$ or $\lnot y \ro{R}{C} x$ by anti-symmetry,
	hence $x=y$ or $y \rr{R} x$.
\item
	By definition.
\item
	``$\Ra$'':
	If $x \ro{R}{7} y$ and $y \ro{R}{7} x$,
	then $x \rr{R} y \lor y \rr{R} x$;
	both cases imply $x=y$ by co-reflexivity.
	\\
	``$\La$'':
	Let $x \rr{R} y$,
	then $x \rr{R} y \lor y \rr{R} x$ by weakening;
	this implies both $x \ro{R}{7} y$ and $y \ro{R}{7} x$,
	hence $x=y$ by anti-symmetry.
\item
	By definition.
\qed
\end{enumerate}
}

\LEMMA{}{%
\LABEL{optn semiord2 1}%
\begin{enumerate}
\item\LABEL{optn semiord2 1 1}%
	If $R$ is transitive, then its complement $R^C$ satisfies
	semi-order property~2.
\item\LABEL{optn semiord2 1 2}%
	If $R$ is symmetric and satisfies semi-order property~2,
	then its complement $R^C$ is transitive.
\end{enumerate}
}
\PROOF{
\begin{enumerate}
\item
	Let $x \ro{R}{C} y$ and $y \ro{R}{C} z$ hold,
	assume for contradiction that neither
	$w \ro{R}{C} x$ nor $x \ro{R}{C} w$ nor
	$w \ro{R}{C} y$ nor $y \ro{R}{C} w$ nor
	$w \ro{R}{C} z$ nor $z \ro{R}{C} w$ holds.
	From $x \rr{R} w$ and $w \rr{R} y$, 
	we obtain $x \rr{R} y$ by transitivity,
	contradicting our assumption.
\item
	Let $x \ro{R}{C} y$ and $y \ro{R}{C} z$ hold,
	assume for contradiction $x \rr{R} z$.
	Then $z \rr{R} x$ by symmetry, 
	and
	$x \rr{R} y \lor y \rr{R} x 
	\lor y \rr{R} z \lor z \rr{R} y$
	by semi-order property~2.
	Due to symmetry, the latter disjunction boils down to
	$x \rr{R} y \lor y \rr{R} z$, 
	contradicting our assumption.
\qed
\end{enumerate}
}

\LEMMA{}{%
\LABEL{optn trans eucl}%
For $q \in \set{ 0, 1, 6, 7, 8, 9, E, F}$,
a relation $R$ is $q$-left-Euclidean iff it is $q$-transitive.
}
\PROOF{
By Lem.~\REFF{optn univ}{17}, $R^q$ is symmetric.
Hence $q$-left-Euclideanness and $q$-transitivity coincide by
\OEF{Lem.36}{33}{eucl 1}.
\qed
}

\LEMMA{Characterization of symmetry}{%
\LABEL{optn sym}%
All of the following properties are equivalent:
\begin{enumerate}
\item\LABEL{optn sym 2}%
	6-AntiSym
\item\LABEL{optn sym 1}%
	6-ASym
\item\LABEL{optn sym 4}%
	2-, 4-, 6-LfQuasiRefl
\item\LABEL{optn sym 3}%
	2-, 4-, 6-, B-, D-LfEucl
\item\LABEL{optn sym 5}%
	6-, 9-SemiOrd1
\item\LABEL{optn sym 6}%
	9-SemiOrd2
\item\LABEL{optn sym 7}%
	2-, 3-, 4-, 5-, A-, B-, C-, D-Sym
\item\LABEL{optn sym 8}%
	6-Trans
\end{enumerate}
}
\PROOF{
We show for each property that is applies to a relation $R$ iff $R$
is symmetric.

For the ``if'' part, observe that for a symmetric $R$,
we have $R^0 = \ro{R}{2} = \ro{R}{4} = R^6$ empty 
by Lem.~\REFF{optn 2}{1},
and hence asymmetric and anti-symmetric 
by \OEF{Exm.74.6}{48}{empty 6},
left Euclidean by \OEF{Exm.74.2}{48}{empty 2},
quasi-reflexive by \OEF{Exm.74.1}{48}{empty 1},
symmetric by \OEF{Exm.74.3}{48}{empty 4},
transitive by \OEF{Exm.74.7}{48}{empty 7},
and satisfy semi-order property~1 and~2 
by \OEF{Exm.74.9}{48}{empty 8}.
Moreover, for a symmetric $R$, 
we have that $R^9 = \ro{R}{B} = \ro{R}{D} = R^F$
are universal by Lem.~\REFF{optn 2}{1},
and hence left Euclidean by \OEF{Exm.75.1}{48}{univ 1},
symmetric by \OEF{Exm.75.2}{48}{univ 2}
and satisfy semi-order property~1 and~2 
by \OEF{Exm.75.5}{48}{univ 5}.
And for symmetric $R$ we have $R^5$ symmetric by self-duality,
$R^C$ symmetric by contraposition, and $R^A$ symmetric by duality to
$R^C$.

For the ``only if'' part:
\begin{enumerate}
\item
	Let $\ro{R}{6}$ be anti-symmetric, we show that $R$ is symmetric.
	Let $x \rr{R} y$, 
	assume for contradiction $\lnot y \rr{R} x$;
	then $x \neq y$.
	Moreover, $x \ro{R}{6} y$ by definition,
	hence $\lnot y \ro{R}{6} x$ by anti-symmetry,
	contradicting Lem.~\REFF{optn univ}{17}.
\item
	Follows from~\REF{optn sym 2} and
	\OEF{Lem.13.2}{25}{asym 1b 2}.
\item
	\begin{itemize}
	\item Case $2$:
		Let $x \rr{R} y$, 
		assume for contradiction $\lnot y \rr{R} x$.
		Then $x \ro{R}{2} y$,
		hence $x \ro{R}{2} x$ by left quasi-reflexivity;
		this contradicts Lem.~\REFF{optn univ}{4}.
	\item Cases $4$ and $6$ are shown similar.
	\end{itemize}
\item
	\begin{itemize}
	\item Cases $2$, $4$, and $6$ follow
		from~\REF{optn sym 4},
		using \OEF{Lem.46}{37}{eucl 11}.
	\item Case $B$:
		Let $x \rr{R} y$, then $x \ro{R}{B} y$ by definition.
		Moreover $y \ro{R}{B} y$ 
		by Lem.~\REFF{optn univ}{14},
		hence $y \ro{R}{B} x$ by Euclideaness,
		hence $y \rr{R} x \lor \lnot x \rr{R} y$ 
		by definition,
		hence $y \rr{R} x$ by our assumption.
	\item Case $D$ is similar.
	\end{itemize}
\item
	\begin{itemize}
	\item Case $6$:
		Let $x \rr{R} y$, 
		assume for contradiction $\lnot y \rr{R} x$.
		Then $x \ro{R}{6} y$ by definition,
		hence $y \ro{R}{6} x$ by Lem.~\REFF{optn univ}{17}.
		Moreover, $\lnot x \ro{R}{6} x$
		by Lem.~\REFF{optn univ}{14}.
		Applying semi-order property~1 to
		$y \ro{R}{6} x$, $\lnot x \ro{R}{6} x$, and $x \ro{R}{6} y$
		infers $y \ro{R}{6} y$,
		contradicting Lem.~\REFF{optn univ}{4}.
	\item Case $9$:
		Let $x \rr{R} y$, 
		assume for contradiction $\lnot y \rr{R} x$.
		Then $\lnot x \ro{R}{9} y$ and also $\lnot y \ro{R}{9} x$
		by definition.
		Moreover, $x \ro{R}{6} x$ and $y \ro{R}{9} y$
		by Lem.~\REFF{optn univ}{14}.
		Applying semi-order property~1 to $x \ro{R}{9} x$,
		$x,y$ incomparable w.r.t.\ $R^9$, and $y \ro{R}{9} y$
		infers $x \ro{R}{9} y$,
		contradicting the incomparability of $x,y$.
	\end{itemize}
\item
	By Lem.~\REFF{optn univ}{14},
	$R^9$ is reflexive,
	hence also connex by \OEF{Lem.66}{45}{semiOrd1 1a}.
	Since $R^9$ is also symmetric by Lem.~\REFF{optn univ}{17},
	it is universal by \OEF{Lem.53.2}{40}{conn 4 2}.
	By definition, 
	this means that $x \rr{R} y$ and $y \rr{R} x$ agree
	everywhere, i.e.\ that $R$ is symmetric.

\item
	\begin{itemize}
	\item Case $2$:
		Let $x \rr{R} y$, 
		assume for contradiction $\lnot y \rr{R} x$.
		Then $x \ro{R}{2} y$ by definition,
		hence $y R^2$ by symmetry,
		hence $y \rr{R} x$ by definition, 
		contradicting our assumption.
	\item Case $B$:
		Let $x \rr{R} y$,
		then $x \ro{R}{B} y$ by definition,
		hence $y \ro{R}{B} x$ by symmetry,
		hence $y \rr{R} x$ by definition,
		since $\lnot x \rr{R} y$ cannot hold.
	\item Cases $4$ and $D$ follow by duality.
	\item Case $3$ is trivial.
	\item Cases $5$, $A$, and $C$
		follow by duality and contraposition.
	\end{itemize}
\item
	Let $x \rr{R} y$, 
	assume for contradiction $\lnot y \rr{R} x$.
	Then, $x \ro{R}{6} y$ and also $y \ro{R}{6} x$ by definition,
	hence $x \ro{R}{6} x$ by transitivity,
	contradicting Lem.~\REFF{optn univ}{4}.
\qed
\end{enumerate}
}

\LEMMA{Symmetry and semi-order property~1}{%
\LABEL{optn semiord1}%
If $R$ is symmetric and satisfies semi-order property~1,
then $R^C$ does, too.
}
\PROOF{
First, we have
	$R$ symmetric
	~ iff ~ $x \rr{R} y \lra y \rr{R} x$
	~ iff ~ $\lnot x \rr{R} y \lra \lnot y \rr{R} x$
	~ iff ~ $R^C$ symmetric.

Now let $R$ be symmetric and satisfy semi-order property~1,
let $w \ro{R}{C} x$, $\lnot x \ro{R}{C} y$, $\lnot y \ro{R}{C} x$, 
and $y \ro{R}{C} z$ hold.
Assume for contradiction $\lnot w \ro{R}{C} z$.
Then by definition
$x \rr{R} y$, $\lnot y \rr{R} z$, and $w \rr{R} z$,
hence by symmetry
$\lnot z \rr{R} y$, and $z \rr{R} w$.
Applying semi-order property~1 to these facts implies $x \rr{R} w$,
hence $w \rr{R} x$, contradicting $w \ro{R}{C} x$.
\qed
}

\LEMMA{}{%
\LABEL{optn dense}%
\begin{enumerate}
\item\LABEL{optn dense 1}%
	A relation is 2-dense iff it is 4-dense.
\item\LABEL{optn dense 2}%
	If $R$ is symmetric, then $R^2$ is dense.
\item\LABEL{optn dense 3}%
	The converse direction does not hold if the universe set
	has $\geq 7$ elements.
\end{enumerate}
}
\PROOF{
\begin{enumerate}
\item
	$R^4$ is the converse relation of $R^2$, since $4 = 5 \compose 2$.
	Since the definition of density is self-dual, we are done.
\item
	By Lem.~\REFF{optn 2}{1}, $R^2$ is empty, hence dense by
	\OEF{Exm.74.11}{48}{empty 9}.
\item
	\OEF{Exm.76}{49}{exm 186}
	and Fig.~\REF{Counter-example relation for aiir and qij0} here
	show an example relation on a $7$-element
	domain that is 2-dense and not symmetric.
\qed
\end{enumerate}
}

\LEMMA{}{%
\LABEL{optn antisym 1}%
The following are equivalent:
\begin{enumerate}
\item\LABEL{optn antisym 1 6}%
	$R$ is anti-symmetric and semi-connex.
\item\LABEL{optn antisym 1 1}%
	$R^9$ is the identity relation.
\item\LABEL{optn antisym 1 5}%
	$R^9$ is left-unique.
\item\LABEL{optn antisym 1 2}%
	$R^9$ is anti-symmetric.
\item\LABEL{optn antisym 1 3}%
	$R^B$ is anti-symmetric
\item\LABEL{optn antisym 1 4}%
	$R^D$ is anti-symmetric.
\end{enumerate}
}
\PROOF{
By Lem.~\REFF{optn univ}{17} and Lem.~\REFF{optn univ}{14},
$R^9$ is always symmetric and reflexive.
\begin{itemize}
\item \REF{optn antisym 1 6} $\Ra$ \REF{optn antisym 1 1}:
	Let $x \ro{R}{9} y$ hold, by definition, we have two cases:
	\begin{itemize}
	\item If $x \rr{R} y \land y \rr{R} x$,
		then $x=y$ by anti-symmetry of $R$.
	\item If $\lnot x \rr{R} y \land \lnot y \rr{R} x$,
		then $x=y$ since $R$ is semi-connex.
	\end{itemize}
	This shows that $R^9$ is co-reflexive.
	Since it is also reflexive, it must be the identity,
	by \OEF{Lem.5.1}{21}{}.
\item \REF{optn antisym 1 1} $\Ra$ \REF{optn antisym 1 6}:
	\begin{itemize}
	\item Anti-symmetry:
		Let $x \rr{R} y \land y \rr{R} x$,
		then $x \ro{R}{9} y$ by definition, hence $x=y$.
	\item Semi-connex:
		Let $\lnot x \rr{R} y \land \lnot y \rr{R} x$,
		then again $x \ro{R}{9} y$ by definition, hence $x=y$.
	\end{itemize}
\item \REF{optn antisym 1 1}
	$\Ra$ \REF{optn antisym 1 5}
	$\land$ \REF{optn antisym 1 2}:
	trivial.
\item \REF{optn antisym 1 2}
	$\Ra$ \REF{optn antisym 1 1}:
	\\
	If $R^9$ is anti-symmetric,
	then it is co-reflexive by \OEF{Lem.7.7}{22}{corefl 6 5},
	hence the identity by \OEF{Lem.5.1}{21}{corefl 2 1}.
\item \REF{optn antisym 1 5} $\Ra$ \REF{optn antisym 1 1}:
	$R^9$ is co-reflexive by \OEF{Lem.7.1}{22}{corefl 6 1},
	hence the identity by \OEF{Lem.5.1}{21}{corefl 2 1}.
\item \REF{optn antisym 1 3}
	$\lor$ \REF{optn antisym 1 4}
	$\Ra$ \REF{optn antisym 1 2}:
	\\
	By Lem.~\REFF{optn bool}{4},
	both $R^B$ and $R^D$ are supersets of $R^9$,
	by Lem.~\REFF{Monotonic properties}{4},
	anti-symmetry is antitonic.
\item \REF{optn antisym 1 6}
	$\Ra$ \REF{optn antisym 1 3}
	$\land$ \REF{optn antisym 1 4}:
	Let $x \ro{R}{B} y$ and $y \ro{R}{B} x$,
	then by definition
	$x \rr{R} y \lor \lnot y \rr{R} x$ and
	$y \rr{R} x \lor \lnot x \rr{R} y$.
	This is logically equivalent to
	$(x \rr{R} y \land y \rr{R} x) 
	\lor (\lnot x \rr{R} y \land \lnot y \rr{R} x)$,
	that is, to $x \ro{R}{9} y$.
	Since $R^9$ is the identity as shown above,
	we have $x=y$.
	The proof for \REF{optn antisym 1 4} is dual.
\qed
\end{itemize}
}

\LEMMA{Derived equivalences}{%
\LABEL{Derived equivalences}%
Let ${\it prop}$, ${\it prop}_i$ denote properties of binary relations.
\begin{enumerate}
\item\LABEL{Derived equivalences 1}%
	If
	$\forall R. \; {\it prop}_1(R^{p_1}) \ra {\it prop}_2(R^{p_2})$,
	then
	$\forall R. \; {\it prop}_1(R^{p_1 \compose r})
	\ra {\it prop}_2(R^{p_2 \compose r})$
	for every operation $r$.
\item\LABEL{Derived equivalences 2}%
	If 3-{\it prop} $\lra$ 5-{\it prop},
	then
	2-{\it prop} $\lra$ 4-{\it prop},
	A-{\it prop} $\lra$ C-{\it prop},
	and
	B-{\it prop} $\lra$ D-{\it prop}.
	\\
	This applies to the properties
	Refl, QuasiRefl, Sym, ASym, AntiSym, Trans, AntiTrans,
	SemiOrd1, SemiOrd2, and Dense.
\item\LABEL{Derived equivalences 3}%
	If 1-{\it prop} $\lra$ 3-{\it prop},
	then 0-{\it prop} $\lra$ 4-{\it prop},
	1-{\it prop} $\lra$ 5-{\it prop},
	8-{\it prop} $\lra$ A-{\it prop},
	9-{\it prop} $\lra$ B-{\it prop},
	and
	8-{\it prop} $\lra$ C-{\it prop}.
	\\
	This applies to the properties
	Refl, ASym, and AntiSym.
\item\LABEL{Derived equivalences 4}%
	If 3-{\it prop} $\lra$ 7-{\it prop},
	then
	2-{\it prop} $\lra$ 6-{\it prop},
	4-{\it prop} $\lra$ 6-{\it prop},
	5-{\it prop} $\lra$ 7-{\it prop},
	A-{\it prop} $\lra$ E-{\it prop},
	B-{\it prop} $\lra$ F-{\it prop},
	C-{\it prop} $\lra$ E-{\it prop},
	and
	D-{\it prop} $\lra$ F-{\it prop}.
	\\
	This applies to the properties
	Refl, and QuasiRefl.
\item\LABEL{Derived equivalences 5}%
	If 1-{\it prop} $\lra$ E-{\it prop},
	then
	0-{\it prop} $\lra$ F-{\it prop},
	6-{\it prop} $\lra$ 9-{\it prop},
	and
	7-{\it prop} $\lra$ 8-{\it prop}.
	\\
	This applies to the properties Sym, and SemiOrd1.
\end{enumerate}
}
\PROOF{
\begin{enumerate}
\item
	Since $R^{p_i \compose r} = (R^r)^{p_i}$,
	the consequent is an instance of the antecedent.
\item
	Apply~\REF{Derived equivalences 1} to $r=2, A, B$.
	The listed properties are self-dual by Def.~\REF{def}.
\item
	Apply~\REF{Derived equivalences 1}
	to $r=4, 5, A, B, C$.
	To prove the list:
	\begin{itemize}
	\item If $R$ is reflexive, then $R^1$ is, too,
		by Lem.~\REFF{optn univ}{5}.
		If $R^1$ is reflexive, then its superset $R^3$ is, too,
		by Lem.~\REFF{Monotonic properties}{3}.
	\item If $R$ is asymmetric, then its subset $R^1$ is, too,
		by Lem.~\REFF{Monotonic properties}{4}.
		If $R^1$ is asymmetric, it is empty
		by \OEF{Lem.16}{25}{sym 2}
		using Lem.~\REFF{optn univ}{17},
		that is, 
		$x \rr{R} y \land y \rr{R} x$ can never happen,
		that is, $R$ is asymmetric.
	\item If $R$ is anti-symmetric, then its subset $R^1$ is, too,
		by Lem.~\REFF{Monotonic properties}{4}.
		If $R^1$ is is anti-symmetric,
		then it is co-reflexive
		by \OEF{Lem.7.7}{22}{corefl 6 5}
		using Lem.~\REFF{optn univ}{17},
		that is, 
		$x \rr{R} y \land y \rr{R} x$ 
		can happen only for $x=y$,
		that is, $R$ is anti-symmetric.
	\end{itemize}
\item
	Apply~\REF{Derived equivalences 1}
	to $r=2, 4, 5, A, B, C, D$.
	To prove the list:
	\begin{itemize}
	\item If $R$ is reflexive, then $R^7$ is, too,
		by Lem.~\REFF{optn univ}{5}.
		If $R^7$ is reflexive,
		then $x \rr{R} x \lor x \rr{R} x$ for all $x$,
		hence $R$ is reflexive.
	\item Let $R$ be quasi-reflexive, let $x \ro{R}{7} y$ hold.
		Then $x \rr{R} y \lor y \rr{R} x$ by definition.
		Each alternative implies 
		$x \rr{R} x \land y \rr{R} y$,
		and hence $x \ro{R}{7} x \land y \ro{R}{7} y$.
		\\
		Conversely, let $R^7$ be quasi-reflexive,
		let $x \rr{R} y$ hold.
		Then $x \ro{R}{7} y$,
		hence $x \ro{R}{7} x \land y \ro{R}{7} y$,
		which boils down to 
		$x \rr{R} x \land y \rr{R} y$.
	\end{itemize}
\item
	Apply~\REF{Derived equivalences 1}
	to $r=0, 6, 7$.
	To prove the list:
	\begin{itemize}
	\item Both $R^1$ and $R^E$ are always symmetric
		by Lem.~\REFF{optn univ}{17}.
	\item For semi-order property~1, the equivalence follows from
		Lem.~\REF{optn semiord1},
		since $R^E = (R^1)^C$ and $R^1 = (R^E)^C$.
\qed
	\end{itemize}
\end{enumerate}
}

\THEOREM{Equivalent lifted properties}{%
\LABEL{optn equiv}%
On a set with $\geq 7$ elements,
the equivalence classes of lifted properties coincide with the
partition shown in Fig.~\REF{Computed partitions of lifted properties},
except that an own class $\set{ \mbox{2-Dense}, \mbox{4-Dense} }$
must be split off from the partition of 3-Sym (indicated by the dotted
line).
}
\PROOF{
First, 2-Dense is not equivalent to 3-Sym:
for the usual $<$ relation on the rational numbers,
3-Dense and 2-Dense coincide since $<$ is asymmetric;
therefore $<$ satisfies 2-Dense, but not 3-Sym;
see also \OEF{Exm.76}{49}{exm 186}.
This shows the split is necessary.

By construction of the partition split algorithm from
Fig.~\REF{Partition split algorithm},
no two lifted properties from different partitions
in Fig.~\REF{Computed partitions of lifted properties}
can agree on all relations.
Therefore, it is sufficient to show that each two members of a
partition are equivalent, with the above exception.

\begin{itemize}
\item Boxes:
	\\
	All properties in the topmost left and right box are equivalent by
	Lem.~\REF{optn unsat} and~\REF{optn univ},
	respectively.
	The two members of the split-off class are equivalent by
	Lem.~\REFF{optn dense}{1};
	all members of the remainder of that
	partition are equivalent by Lem.~\REF{optn sym}.

	The bottommost $5$ right boxes are covered by
	Lem.~\REF{optn semiord2 1}
	(using symmetry of $R^1$, $R^7$, $R^8$, $R^9$, $R^E$
	by Lem.~\REFF{optn univ}{17})
	and~\REF{optn trans eucl}.

	The right box second from top is covered
	by Lem.~\REF{optn antisym 1}.

\item Pure partitions:
	\\
	The ASym and the AntiSym partitions
	are covered
	by Lem.~\REFF{Derived equivalences}{3};
	the Refl partitions
	by this
	and Lem.~\REFF{Derived equivalences}{4}.

	The SemiOrd1 partition are covered
	by Lem.~\REFF{Derived equivalences}{2}
	and Lem.~\REFF{Derived equivalences}{5};
	the SemiOrd2, AntiTrans, Trans, and Dense partitions by
	Lem.~\REFF{Derived equivalences}{2}.

	Note that not all equivalences implied
	by Lem.~\REF{Derived equivalences}
	lead to pure partitions,
	e.g.\ B-Dense $\lra$ D-Dense holds trivially
	(by Lem.~\REFF{optn univ}{14}
	and~\OEF{48.1}{38}{den 1 1})
	and is therefore reflected in the universal (``$+$'') partition.

\item Singleton Partitions:
	\\
	Nothing to show.
\qed
\end{itemize}
}

\LEMMA{Default representations}{%
\begin{enumerate}
\item For every ``$+$'' in the matrix
	in Fig.~\REF{Default representations of lifted properties},
	the operation of its column always yields a relation satisfying the
	property of its row.
	For example, operation $9$ always yields a reflexive and
	symmetric relation.
\item For every ``$-$'', the column operation never yields a relation
	satisfying the row property.
	For example, operation $B$ never yields a anti-transitive or
	left-unique relation.
\item Whenever in a row {\it prop} the default representations
	coincide for column $p$ and $q$, and $q = p \compose r$,
	then $r$ preserves $p$-{\it prop}.
	For example, operations 2, 3, 4, 5, A, B, C, D all preserve
	symmetry.
\end{enumerate}
This generalizes Lem.~\REF{optn univ}.
}
\PROOF{
By Thm.~\REF{optn equiv},
two properties are equivalent iff
Fig.~\REF{Default representations of lifted properties}
shows the same default representation for them.
\begin{enumerate}
\item
	If operation $p$-{\it prop} has default representation $+$,
	then $R^p$ always satisfies {\it prop}.
\item
	If operation $p$-{\it prop} has default representation $-$,
	then $R^p$ never satisfies {\it prop}.
\item
	Let $R$ satisfy $p$-{\it prop},
	then $R$ equivalently satisfies $q$-{\it prop}.
	That is, $\ro{(\ro{R}{r})}{p} = \ro{R}{q}$ satisfies {\it prop},
	hence $\ro{R}{r}$ satisfies $p$-{\it prop}.
\qed
\end{enumerate}
}

\LEMMA{}{
A relation $R$ is E-AntiSym iff its complement is CoRefl.
}
\PROOF{
$\ro{R}{E}$ is AntiSym
iff $\ro{R}{7 \compose C}$ is AntiSym
iff $\ro{R}{C}$ is 7-AntiSym
iff $\ro{R}{C}$ is CoRefl.
\qed
}

\clearpage

\section{Implications between lifted properties}
\LABEL{Implications between lifted properties}

The approach of Sect.~\REF{Equivalent lifted properties} could produce
only ``equivalence laws'', of the form
$$\forall R. \;
{\it prop}_1({\it op}_1(R)) \lra {\it prop}_2({\it op}_2(R)) .$$
On the other hand, investigating all ``disjunction laws'', of the form
$$\forall R. \; \:{\lor}1n{{\it prop}_\i({\it op}_\i(R))} ,$$
for ${\it prop}_i$ a negated or unnegated property,
is computationally infeasible for large $n$, as discussed in
Sect.~\REF{Equivalent lifted properties}.
As a compromise,
we investigate in this section ``$3$-implication laws'', of the form
$$\forall R. \;
{\it prop}_1({\it op}_1(R)) \land
{\it prop}_2({\it op}_2(R)) \ra
{\it prop}_3({\it op}_3(R)) ,$$
that is, implications of unnegated lifted properties with up to $2$
antecedents and one conclusion.\footnote{
	In \cite{Burghardt.2018c}, $58$ of the found ``disjunction laws''
	(of basic properties only)
	had $\geq 4$ literals;
	another $10$ involved $3$ literals, but no negation;
	the remaining $206$ laws (= $75\%$) could be written as
	``$3$-implication laws''.
	This supports our speculation that for lifted properties,
	too, ``$3$-implication laws'' will cover the majority of
	``disjunction laws''.
}
Note that we don't cover the negation of a lifted property;
in particular, it cannot be obtained by
composing it with the $C$ operation; for example, ``$R$ is not reflexive''
is usually different from ``$R$ is $C$-reflexive, a.k.a.\ irreflexive''.
However, we can express $\true$ and $\false$ as 0-Sym and 0-Refl,
respectively;\footnote{
	The latter presupposes a nonempty relation domain.
}
hence we can express incompatibility between two lifted properties.

Since the Quine-McCluskey procedure would require too much
computational resources, we used a different approach.
We cycled though all $81^3 = 531441$ triples
$\langle
{\it op}_1\mbox{-}{\it prop}_1,
{\it op}_2\mbox{-}{\it prop}_2,
{\it op}_3\mbox{-}{\it prop}_3
\rangle$
of default representation of lifted properties, and for each of them
searched a counter-example relation
$R$ such that
$
{\it prop}_1({\it op}_1(R)) \land
{\it prop}_2({\it op}_2(R)) \land \lnot
{\it prop}_3({\it op}_3(R))
$
holds.
We recorded the search result in an array, indexed by the rank of the
three default representations.\footnote{
	The equivalence classes of lifted properties shown in
	Fig.~\REF{Computed partitions of lifted properties}
	are numbered in the following order, cf.\ the list
	\texttt{lpnfList[]} in file
	\filename{lpImplicationsTables.c}:
	boxed partitions starting at rank 0:
	0-rf, 0-sy, 3-sy, 2-de, 9-an, 1-tr, 7-tr, 8-tr, 9-tr, E-tr;
	pure partitions starting at rank 10, top to bottom, left to
	right: 3-as, C-as, 3-an, C-an, \ldots, C-de;
	singleton partitions starting at rank 33, in the same order:
	7-as, E-as, 7-an, E-an, \ldots, E-lq.
}
In case no counter-example was found, we asked an external resolution
prover\footnote{
	EProver version {\em E 2.3 Gielle\/},
	from \url{http://www.eprover.org}
}
to prove
$
{\it prop}_1({\it op}_1(R)) \land
{\it prop}_2({\it op}_2(R)) \ra
{\it prop}_3({\it op}_3(R))
$,
which might, or might not, succeed.
In the former case, we recorded this result in the array.
The remaining array cells were handled manually, unless they were
trivial (Sect.~\REF{Trivial inferences}).

\subsection{Referencing an implication}
\LABEL{Referencing an implication}

To reference a $3$-implication under consideration,
we use its array indices, encoded in base $27$, with digits
represented as \texttt{0}, \texttt{a}, \texttt{b}, \texttt{c}, \ldots,
\texttt{z}.
For example, the implication
``SemiOrd1 $\land$ Refl $\ra$ Connex''
corresponds to array indices \texttt{[18][14][11]},
that is, to array
cell $(18 \cdot 81+14) \cdot 81+11 = 119243 = ((6 \cdot 27+1) \cdot
27+15) \cdot 27+11$,
and is therefore encoded as \QRef{faok}.
This encoding ensures an implication is referenced by a fixed code,
unique across multiple program runs, and independent of the order in
which laws are found.
The shell script \filename{resolveImplCode.sh} can be used to
convert base-$27$ codes into implications;
the call ``\texttt{nonprominentProperties -lpImpl resolve}''
can be used for the reverse direction.

\subsection{Trivial inferences}
\LABEL{Trivial inferences}

We implemented routines to draw trivial conclusions from a found
negative (counter-example) or positive (proof) result.
An inference is considered trivial if it doesn't require any knowledge
about basic property definitions (as given in Def.~\REF{def}), except for
monotonicity information as provided in
Lem.~\REFF{Monotonic properties}{3}
and~\REFF{Monotonic properties}{4}.

The trivial inference mechanism is intended to reduce the load on both
the counter-example search and the external theorem prover.
For example, when ``LfEucl $\ra$ LfQuasiRefl'' is already established,
we needn't attempt to prove ``RgEucl $\ra$ RgQuasiRefl''
since this follows trivially by taking the converse relation.
When a trivial proof or disproof of an implication
is already known, we run neither
the counter-example search nor the external theorem prover for it.
All trivial proofs are printed to a log file during a program run;
that output may subsequently be used to extract the trivial part of
a proof tree of a
given implication, by running program \filename{justify}.
As an example, Fig.~\REF{Proof tree of ``tr land ir ra as''}
shows the proof tree
for ``Trans $\land$ Irrefl $\ra$ ASym'' (base-$27$ code \QRem{ijrj});
it also demonstrates that trivial proofs needn't be intuitively obvious.
Figure~\REF{Disproof tree of ``as land tr not ra rf''}
shows the disproof tree of ``ASym $\land$ Trans $\not\ra$ Refl''
(\QRem{clcn}); right to each base-27 code, the corresponding index triple
is shown.

\DEFINITION{Trivial inference rules}{%
\LABEL{Trivial inference rules}%
Using $I, J, K, I(\cdot), \ldots$ to denote default representations of 
lifted properties,
our trivial inference rules are:
\begin{enumerate}
\item
	if $I \land J \ra K$, then $J \land I \ra K$;
\item
	if $J \land I \not\ra K$, then $I \land J \not\ra K$;
\item\LABEL{Trivial inference rules drawConclusionsTrueTrans}
	if $I \land J \ra K$, and $H \land I \ra J$,
	then $H \land I \ra K$,
	and variants thereof;
\item 
	if $H \land I \not\ra K$, and $H \land I \ra J$,
	then $I \land J \not\ra K$,
	and variants thereof;
\item\LABEL{Trivial inference rules drawConclusionsTrueComp}
	if $\forall R. \; I(R) \land J(R) \ra K(R)$,
	then
	$\forall R, p. \;
	I(\ro{R}{p}) \land J(\ro{R}{p}) \ra K(\ro{R}{p})$;
\item 
	if $\exists R, p. \;
	I(\ro{R}{p}) \land J(\ro{R}{p}) \not\ra K(\ro{R}{p})$
	then
	$\exists R. \; I(R) \land J(R) \not\ra K(R)$.
\qed
\end{enumerate}
}

The above example trivial inference from ``LfEucl $\ra$ LfQuasiRefl''
to ``RgEucl $\ra$ RgQuasiRefl''
is achieved
by applying rule~\REFF{Trivial inference rules}{drawConclusionsTrueComp}
with ${\it op} = 5$.

In the last two rules, it does make sense to consider all lifted
properties equivalent to a given one during the search for a common
operation ${\it op}$.
For example, from the implication
C-tr $\land$ 3-rf $\ra$ C-as (\QRem{isok}),
we may infer
A-tr $\land$ 1-rf $\ra$ 8-as by equivalences,
and then
3-tr $\land$ 8-rf $\ra$ 1-as
from rule~\REFF{Trivial inference rules}{drawConclusionsTrueComp} 
with ${\it op} = A$,
which normalizes to the default representation
3-tr $\land$ C-rf $\ra$ 3-as (\QRem{ijrj});
cf.\ the two topmost lines in the proof tree in
Fig.~\REF{Proof tree of ``tr land ir ra as''}.

Initially, we fill all array cells corresponding to trivial
implications, including monotonicity information:

\DEFINITION{Trivial initialization rules}{%
\LABEL{Trivial initialization rules}%
Our trivial initialization rules are:
\begin{enumerate}
\item $x \land y \ra \true$,
\item\LABEL{Trivial initialization rules init false lf}
	$x \land \false \ra z$,
\item $\false \land y \ra z$,
\item\LABEL{Trivial initialization rules init mono lf}
	$x \land y \ra x'$
	if $x \ra x'$ is known by monotonicity or antitonicity, and
\item
	$x \land y \ra y'$ if $y \ra y'$ is known, similarly.
\qed
\end{enumerate}
}

For example, since ASym is antitonic, we initialize array cell
\texttt{[33][3][10]}, 
corresponding to ``7-ASym $\land$ 2-Dense $\ra$ 3-ASym''
(\QRem{k0ij}), to $\true$, 
by rule~\REFF{Trivial initialization rules}{init mono lf}.

It turned out that no other than the above $6$ inference rules are
needed.
For example,
we don't need an extra ``explosion'' rule:
if $I \land J \ra \false$ is known to hold,
we can infer $I \land J \ra K$ for an arbitrary lifted property $K$
using an appropriate variant of 
rule~\REFF{Trivial inference rules}{drawConclusionsTrueTrans}
and the trivial implication $J \land \false \ra K$
from~\REFF{Trivial initialization rules}{init false lf}.

\begin{figure}
\begin{center}
\begin{picture}(140,85)
\put(5,80){\oval(10,10)}%
\put(0,80){\line(0,-1){10}}%
\put(10,80){\line(0,-1){10}}%
\put(5,70){\oval(10,10)[b]}%
\put(5,70){\makebox(0,0){\texttt{in}}}%
\put(11,70){\vector(4,-1){38}}%
\put(115,80){\oval(10,10)}%
\put(110,80){\line(0,-1){10}}%
\put(120,80){\line(0,-1){10}}%
\put(115,70){\oval(10,10)[b]}%
\put(115,70){\makebox(0,0){\texttt{out}}}%
\put(71,60.5){\vector(4,1){38}}%
\textcolor{bgPartSet}{\put(50,0){\makebox(20,60){\rule{20mm}{60mm}}}}%
\textcolor{fgPartSet}{\put(50,0){\framebox(20,60){}}}%
\textcolor{fgPartSet}{\put(49,59){\makebox(0,0)[tr]{\texttt{[00][00][00]}}}}%
\textcolor{fgPartSet}{\put(37, 8){\makebox(0,0)[r]{$\vdots$}}}%
\textcolor{fgPartSet}{\put(49, 1){\makebox(0,0)[br]{\texttt{[80][80][80]}}}}%
\textcolor{fgPartSet}{\put(55,0){\line(0,1){60}}}%
\textcolor{fgPartSet}{\put(53,59){\makebox(0,0)[t]{$+$}}}%
\textcolor{fgPartSet}{\put(53, 1){\makebox(0,0)[b]{$+$}}}%
\textcolor{fgPartSet}{\put(53,61.8){\makebox(0,0)[b]{val}}}%
\textcolor{fgPartSet}{\put(63,61){\makebox(0,0)[b]{wg}}}%
\textcolor{fgPartSet}{\put(49,53){\makebox(0,0)[tr]{\texttt{[00][00][80]}}}}%
\textcolor{fgPartSet}{\put(49,40){\makebox(0,0)[br]{\texttt{[10][28][14]}}}}%
\textcolor{fgPartSet}{\put(49,35){\makebox(0,0)[br]{\texttt{[14][14][76]}}}}%
\textcolor{fgPartSet}{\put(49,25){\makebox(0,0)[br]{\texttt{[28][15][18]}}}}%
\textcolor{fgPartSet}{\put(53,53){\makebox(0,0)[t]{$+$}}}%
\textcolor{fgPartSet}{\put(53,40){\makebox(0,0)[b]{$-$}}}%
\textcolor{fgPartSet}{\put(53,35){\makebox(0,0)[b]{$+$}}}%
\textcolor{fgPartSet}{\put(53,25){\makebox(0,0)[b]{?}}}%
\textcolor{bgPart}{\put(80,55){\makebox(50,5){\rule{50mm}{5mm}}}}%
\textcolor{fgPart}{\put(80,57){\vector(-1,0){10}}}%
\textcolor{fgPart}{\put(80,55){\framebox(50,5){init}}}%
\textcolor{bgPart}{\put(80,39){\makebox(50,5){\rule{50mm}{5mm}}}}%
\textcolor{fgPart}{\put(70,42){\vector(1,0){10}}}%
\textcolor{fgPart}{\put(80,39){\framebox(50,5){trivial inferences}}}%
\textcolor{fgPart}{\put(130,42){\line(1,0){5}}}%
\textcolor{fgPart}{\put(135,39){\oval(6,6)[r]}}%
\textcolor{fgPart}{\put(135,36){\vector(-1,0){65}}}%
\textcolor{bgPart}{\put(80,22){\makebox(50,5){\rule{50mm}{5mm}}}}%
\textcolor{fgPart}{\put(70,25){\vector(1,0){10}}}%
\textcolor{fgPart}{\put(80,22){\framebox(50,5){counter-example generator}}}%
\textcolor{fgPart}{\put(130,25){\line(1,0){5}}}%
\textcolor{fgPart}{\put(135,22){\oval(6,6)[r]}}%
\textcolor{fgPart}{\put(135,19){\vector(-1,0){65}}}%
\textcolor{bgPart}{\put(80,5){\makebox(50,5){\rule{50mm}{5mm}}}}%
\textcolor{fgPart}{\put(70,8){\vector(1,0){10}}}%
\textcolor{fgPart}{\put(80,5){\framebox(50,5){EProver}}}%
\textcolor{fgPart}{\put(130,8){\line(1,0){5}}}%
\textcolor{fgPart}{\put(135,5){\oval(6,6)[r]}}%
\textcolor{fgPart}{\put(135,2){\vector(-1,0){65}}}%
\put(30,77){\vector(-1,0){19}}%
\thicklines%
\textcolor{fgRel}{\put(35,83){\circle{4}}}%
\textcolor{fgRel}{\put(31,78){\line(1,0){8}}}%
\textcolor{fgRel}{\put(35,81){\line(0,-1){8}}}%
\textcolor{fgRel}{\put(35,73){\line(-1,-1){4}}}%
\textcolor{fgRel}{\put(35,73){\line(1,-1){4}}}%
\put(20,77){\makebox(0,0)[b]{\shortstack{manual\\proofs}}}%
\end{picture}
\caption{Approach to find implications}
\LABEL{Approach to find implications}
\end{center}
\end{figure}
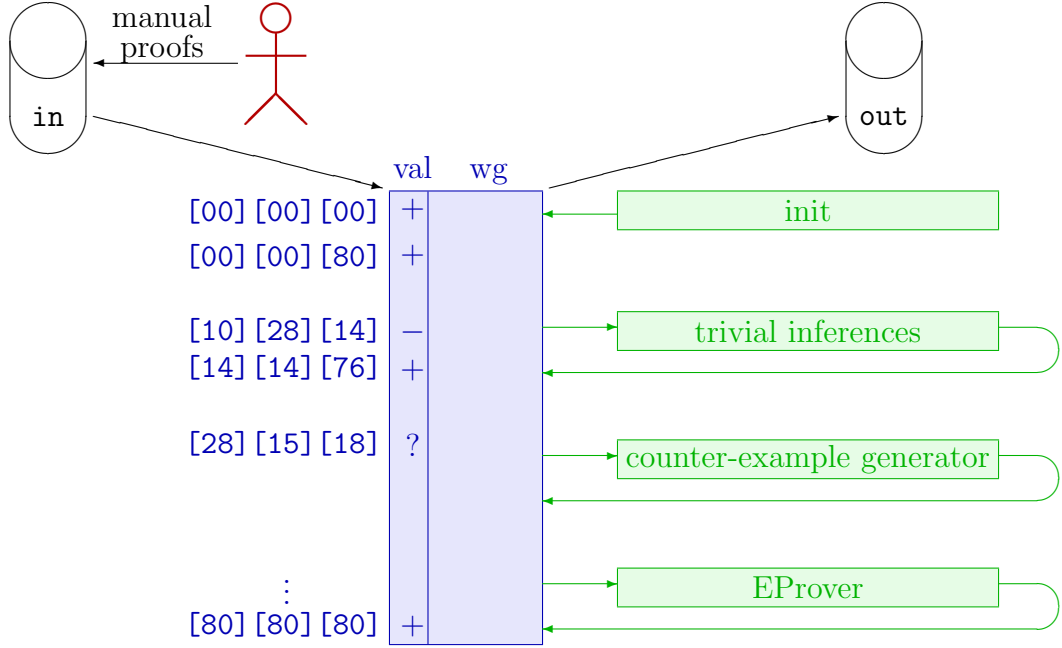

\subsection{Finding implications}
\LABEL{Finding implications}

Figure~\REF{Approach to find implications}
shows an overview of our architecture for searching for $3$-implications.
The array mentioned above is our central data structure.
Various actions on it can be triggered by command-line options.
We first perform the trivial initializations, then usually load
results from previous runs, infer all possible trivial conclusions,
and perform a disproof/proof search for
still undecided array cells.
The latter action cycles through all undecided array cells, and runs the
counter-example generator and (in case of not finding a
counter-example) the external prover; once a decision has been found,
trivial conclusions are drawn from it.

During the early runs, when most cells were yet undecided, we
drew from each new information
as much conclusions as possible,
in order to prune generator and prover calls.
Since this amounted to a depth-first search, it often lead to extremely
long proofs.\footnote{
	In Fig.~\REF{Proof tree of ``tr land ir ra as''}
	and~\REF{Disproof tree of ``as land tr not ra rf''}
	length-minimized proof trees of
	\QRem{ijrj} and \QRem{clcn} are shown.
	The corresponding depth-first versions have $13148$ and
	$124989$ nodes, even when their $3083$ and $34310$ repeated
	subtrees are counted as one node each, respectively.
	%
}
In later runs, when we tried to optimize trivial proofs
w.r.t.\ human readability, we did a breadth-first search
to obtain short proofs.

It took us a couple of program runs to prove or disprove every possible
$3$-implication.
Moreover, we modified details of our approach between successive
runs.\footnote{
	For example, we originally recorded an implication as
	``presumably true'' when no counter-example was found and no
	external proof was available.
	This turned out to be tedious to implement, and inferences
	from presumably true implications weren't of great use;
	so we dropped this kind of truth value.
}
Using a save / load mechanism, we could reuse results from earlier
runs.\footnote{
	Originally, we implemented saving and loading the complete
	array in binary format.
	However, we abandoned this approach in favor of a print / scan
	mechanism for lists of base-$27$ codes, since such files are
	easier to maintain.
}
Consistency over varying implementations is ensured as the notions of
counter-example and external proof remained unchanged.

We kept the log files of all external proof attempts.
If an implication couldn't be proven true immediately, we retried it
with increased minimum domain cardinality, up to $8$.
A domain cardinality of $\geq n$
was expressed by a first-order formula
introducing $n$ constants $\:,1n{c_\i}$, and requiring their pairwise
distinctness; for $n=1$ we used the formula $\lnot \forall x. \false$.
In Fig.~\REF{Minimal axiom set}, we noted any nontrivial minimum domain
cardinality.

Whenever the external prover noted that an implication
$I \land J \ra K$
could be proven without using its conclusion $K$
(EProver outcome ``contradictory axioms''), we
validated that by manually calling the prover for
$I \land J \ra \false$.
Altogether, we collected log files of successful proof attempts of
$4422$ implications of the latter form, and $34727$ other
implications.

A few implications that couldn't be proven this way needed to be proven
manually, viz.\ 
\QRef{g0rw},
\QRef{uivk},
\QRef{uklu},
\QRef{uydr},
\QRef{voye},
\QRef{voyr},
\QRef{wfyp},
see Sect.~\REF{Some proofs of axioms}.
All remaining implications were invalid;
we provided counter examples
on finite (Sect.~\REF{Finite counter-examples}) or countably infinite
(Sect.~\REF{Infinite counter-examples}) domains manually.
Knowing about
the validity / invalidity of every $3$-implication is the
basis to achieve completeness of our axiom set in 
Sect.~\REF{Axiom set} below.

\subsection{Finite counter-examples}
\LABEL{Finite counter-examples}

\EXAMPLE{Finite counter-examples}{
\LABEL{Finite counter-examples exm}%
The following implications are not universally valid:
\begin{enumerate}
\item (\QDef{aiir}) 9-AntiSym    $\land$ 2-Dense     $\not\ra$ 3-SemiOrd1
\item (\QDef{qij0}) 2-LfSerial   $\land$ 2-Dense     $\not\ra$ 2-Trans
\item (\QDef{eiit}) 1-SemiOrd1   $\land$ 2-Dense     $\not\ra$ C-SemiOrd1
\item (\QDef{g0ih}) 2-SemiOrd2   $\land$ 2-Dense     $\not\ra$ 9-Trans
\item (\QDef{jijo}) 3-Dense      $\land$ 2-Dense     $\not\ra$ 1-Dense
\item (\QDef{qbvk}) 1-LfSerial   $\land$ 3-AntiTrans $\not\ra$ 6-AntiTrans
\item (\QDef{qiiu}) 2-LfSerial   $\land$ 2-Dense     $\not\ra$ 2-SemiOrd2
\item (\QDef{r0jq}) 4-LfSerial   $\land$ 4-Dense     $\not\ra$ 7-Dense
\item (\QDef{reuu}) 4-LfSerial   $\land$ 2-LfSerial  $\not\ra$ 2-SemiOrd2
\item (\QDef{qism}) 2-LfSerial   $\land$ 7-Trans     $\not\ra$ 8-AntiTrans
\item (\QDef{qisr}) 2-LfSerial   $\land$ 7-Trans     $\not\ra$ 8-Dense
\item (\QDef{rg0u}) 4-LfSerial   $\land$ 2-LfUnique  $\not\ra$ 2-SemiOrd2
\item (\QDef{rgak}) 4-LfSerial   $\land$ 2-LfUnique  $\not\ra$ 6-AntiTrans
\item (\QDef{rylq}) 6-LfSerial   $\land$ 6-LfUnique  $\not\ra$ 2-SemiOrd1
\item (\QDef{rymo}) 6-LfSerial   $\land$ 6-LfUnique  $\not\ra$ 1-Dense
\item (\QDef{xrkv}) 3-LfQuasiRefl $\land$ 4-Dense   $\not\ra$ 7-LfQuasiRefl
\end{enumerate}
}
\PROOF{
We show 
in Fig.~\REF{Counter-example relation for aiir and qij0}
to~\REF{Counter-example relation for xrkv}
a counter-example relation on a finite domain\footnote{
	The shown counter-examples might not be the simplest possible.
}
for each of the implications.
In each picture,
an arrow from $x$ (light blunt end)
to $y$ (dark peaked end) indicates $x \rr{R} y$, 
colors have only didactic purpose.
All properties have been machine-checked by our implementation.
Additionally,
we give a witness for each dissatisfied implication conclusion:
\begin{enumerate}
\item
	In Fig.~\REF{Counter-example relation for aiir and qij0},
	we have $b \rr{R} a$, $a$ incomparable to itself, 
	$a \rr{R} d$, but not $b \rr{R} d$.
	Hence, $R$ is not SemiOrd1.

	Note that this relation is the same as in
	\OEF{Exm.76}{49}{};
	variations of this relation appear in the counter-examples for
	\QRef{qij0},
	\QRef{eiit},
	\QRef{g0ih},
	\QRef{qiiu},
	\QRef{r0jq},
	\QRef{reuu}, and
	\QRef{xrkv}.
\item
	The same elements in this relation
	demonstrate that $R$ is not Trans,
	and hence not QuasiTrans,
	since $R$ coincides with $\ro{R}{2}$.
\item
	In Fig.~\REF{Counter-example relation for eiit},
	we have $b \ro{R}{C} m$, $m \ro{R}{1} n$, and $n \ro{R}{C} c$,
	but not $b \ro{R}{C} c$.
\item
	In Fig.~\REF{Counter-example relation for g0ih},
	we have $a \ro{R}{9} h$ and $h \ro{R}{9} b$,
	but not $a \ro{R}{9} b$.
\item
	The counter-example relation is obtained from
	Fig.~\REF{Counter-example relation for eiit}
	by removing the loops $m \rr{R} m$ and $n \rr{R} n$.
	We then have no intermediate element for $m \ro{R}{1} n$.
\item
	In Fig.~\REF{Counter-example relation for qbvk},
	we have $a \ro{R}{6} b$ and $b \ro{R}{6} c$,
	but $a \ro{R}{6} c$.
\item
	In Fig.~\REF{Counter-example relation for qiiu, r0jq, reuu},
	we have $b \ro{R}{2} a$ and $a \ro{R}{2} c$,
	while $h$ is incomparable to each of $a,b,c$.
\item
	In the same figure,
	$R$ is not 7-Dense,
	since $a \ro{R}{7} n$ has no
	intermediate element.
\item
	In the same figure,
	$R$ is not 2-SemiOrd2,
	since $b \ro{R}{2} a$ and $a \ro{R}{2} c$,
	but $h$ is incomparable to all of them.
\item
	In Fig.~\REF{Counter-example relation for qism},
	we have
	$x \ro{R}{8} y$, and
	$y \ro{R}{8} z$, but
	$x \ro{R}{8} z$.
\item
	In Fig.~\REF{Counter-example relation for qisr},
	$x \ro{R}{8} y$ has no intermediate element
	w.r.t.\ $\ro{R}{8}$:
	each element is comparable w.r.t.\ $R$
	to $x$ or to $y$.
\item
	In Fig.~\REF{Counter-example relation for rg0u and rgak},
	we have
	$x_1 \ro{R}{2} y_1$ and $y_1 \ro{R}{2} z_1$,
	but $x_2$ is comparable with none of them.
\item
	In the same figure,
	we have $x_1 \ro{R}{6} y_1$ and $y_1 \ro{R}{6} z_1$,
	but $x_1 \ro{R}{6} z_1$,
	i.e.\ $\ro{R}{6}$ is not anti-transitive.
\item
	In Fig.~\REF{Counter-example relation for rylq and rymo},
	$a \ro{R}{2} b$, $b,c$ incomparable w.r.t.\ $\ro{R}{2}$,
	and $d\ro{R}{2} d$,
	but not $a \ro{R}{2} d$.
\item
	In the same figure,
	$b \ro{R}{1} c$ has no intermediate element.
\item
	In Fig.~\REF{Counter-example relation for xrkv},
	we have
	$h \ro{R}{7} c$, but not
	$h \ro{R}{7} h$.
\qed
\end{enumerate}
}

\newcommand{\EXMscale}{0.50}

\begin{figure}
\begin{center}
\includegraphics[scale=\EXMscale]{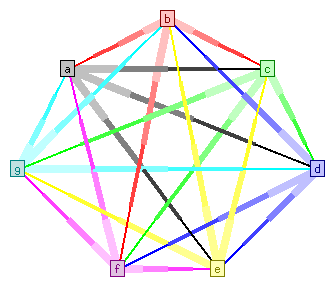}
\caption{Counter-example relation for \QRef{aiir} and \QRef{qij0}}
\LABEL{Counter-example relation for aiir and qij0}
\end{center}
\end{figure}

\begin{figure}
\begin{center}
\includegraphics[scale=\EXMscale]{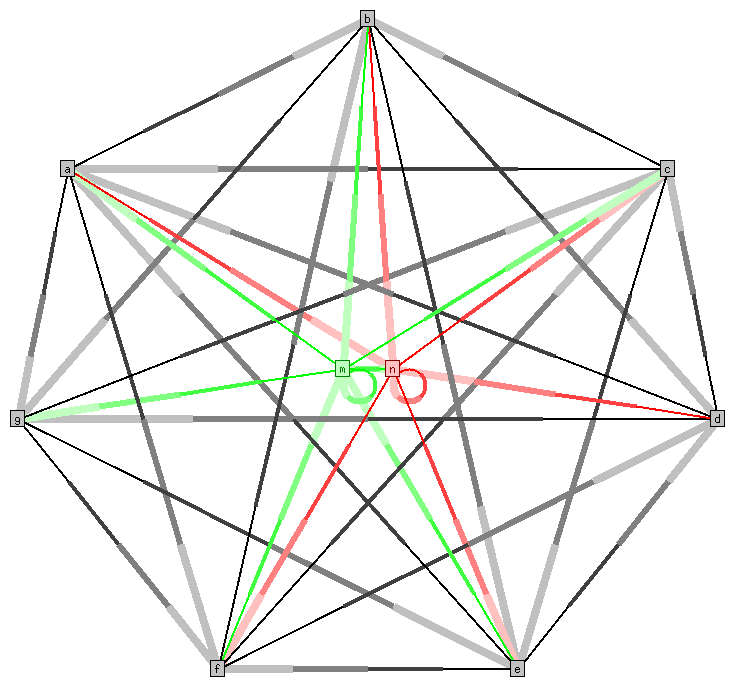}
\caption{Counter-example relation for \QRef{eiit}}
\LABEL{Counter-example relation for eiit}
\end{center}
\end{figure}

\begin{figure}
\begin{center}
\includegraphics[scale=\EXMscale]{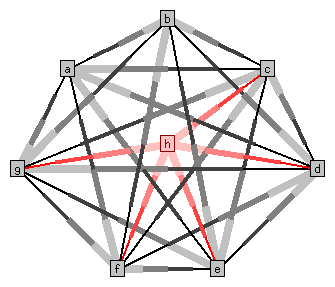}
\caption{Counter-example relation for \QRef{g0ih}}
\LABEL{Counter-example relation for g0ih}
\end{center}
\end{figure}

\begin{figure}
\begin{center}
\includegraphics[scale=\EXMscale]{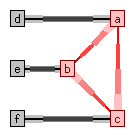}
\caption{Counter-example relation for \QRef{qbvk}}
\LABEL{Counter-example relation for qbvk}
\end{center}
\end{figure}

\begin{figure}
\begin{center}
\includegraphics[scale=\EXMscale]{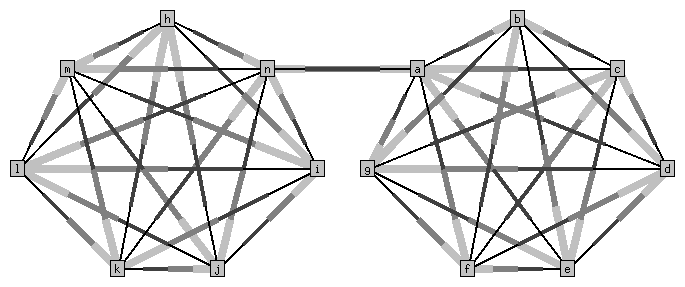}
\caption{Counter-example relation for 
	\QRef{qiiu}, \QRef{r0jq}, and \QRef{reuu}}
\LABEL{Counter-example relation for qiiu, r0jq, reuu}
\end{center}
\end{figure}

\begin{figure}
\begin{center}
\includegraphics[scale=\EXMscale]{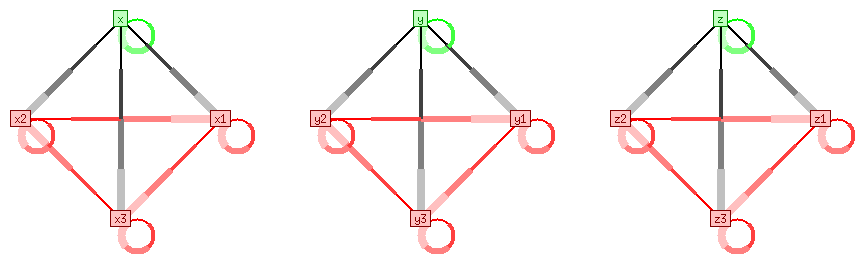}
\caption{Counter-example relation for \QRef{qism}}
\LABEL{Counter-example relation for qism}
\end{center}
\end{figure}

\begin{figure}
\begin{center}
\includegraphics[scale=\EXMscale]{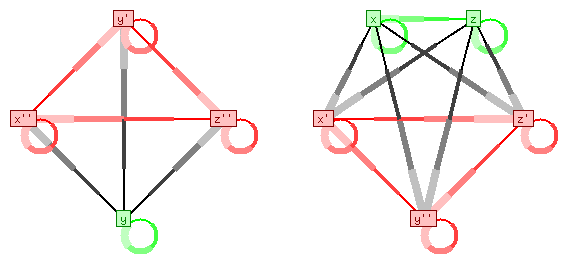}
\caption{Counter-example relation for \QRef{qisr}}
\LABEL{Counter-example relation for qisr}
\end{center}
\end{figure}

\begin{figure}
\begin{center}
\includegraphics[scale=\EXMscale]{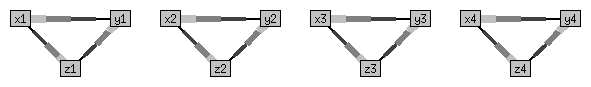}
\caption{Counter-example relation for \QRef{rg0u} and \QRef{rgak}}
\LABEL{Counter-example relation for rg0u and rgak}
\end{center}
\end{figure}

\begin{figure}
\begin{center}
\includegraphics[scale=\EXMscale]{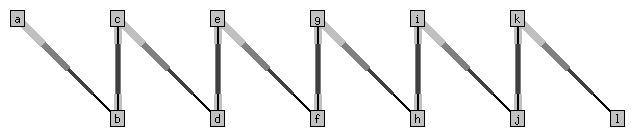}
\caption{Counter-example relation for \QRef{rylq} and \QRef{rymo}}
\LABEL{Counter-example relation for rylq and rymo}
\end{center}
\end{figure}

\begin{figure}
\begin{center}
\includegraphics[scale=\EXMscale]{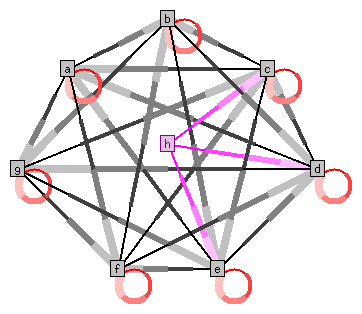}
\caption{Counter-example relation for \QRef{xrkv}}
\LABEL{Counter-example relation for xrkv}
\end{center}
\end{figure}

\clearpage

\subsection{Infinite counter-examples}
\LABEL{Infinite counter-examples}

\EXAMPLE{Non-negative rational numbers}{%
\LABEL{Non-negative rational numbers}
The set $\Q_+ = \set{q \in \Q \mid q \geq 0}$ 
of non-negative rational numbers with the usual order $(<)$
can be used to disprove the following implications:
\begin{enumerate}
\item (\QDef{a0mz}) 4-Dense $\land$ 9-AntiSym $\not\ra$ 3-LfSerial,
\item (\QDef{acib}) 2-Dense $\land$ D-Trans $\not\ra$ 3-Sym,
\item (\QDef{afay}) 2-Dense $\land$ 4-LfSerial $\not\ra$ 2-LfSerial.
\end{enumerate}
}
\PROOF{
We define $R = (<)$ to establish the usual notation.
Since $R$ is asymmetric, $\ro{R}{2i}$ coincides with $\ro{R}{2i+1}$,
for $i=\:,07\i$.
We therefore have
$R^2 = R^3 = (<)$,
$R^4 = (>)$,
$R^9 = (=)$, and
$R^D = (\geq)$.
It is well known that
$(<)$ and $(>)$ are dense,
$(>)$ is left-serial,
$(=)$ is trivially anti-symmetric, and
$(\geq)$ is transitive.
However, $(<)$ is neither left-serial (consider $x=0$) nor symmetric.
\qed
}

\EXAMPLE{Integer numbers}{%
The set $\Z$ of integer numbers with the usual order $(<)$
can be used to disprove the following implications:
\begin{enumerate}
\item (\QDef{jeuc}) D-Trans $\land$ 2-LfSerial $\not\ra$ 2-Dense,
\item (\QDef{jeux}) D-Trans $\land$ 2-LfSerial $\not\ra$ 2-AntiTrans.
\end{enumerate}
}
\PROOF{
As in Exm.~\REF{Non-negative rational numbers},
we write $R$ for $(<)$, and observe $R^2 = (<)$ and $R^D = (\geq)$.
The former is left-serial, and the latter is transitive.
However, $(<)$ is neither dense nor anti-transitive.
\qed
}

\EXAMPLE{aaje}{%
2-Dense and 3-AntiSym doesn't imply C-Dense.
}
\PROOF{
Consider $\Q \cup \set{\infty}$ with $R$ being the usual order
$(\leq)$ on $\Q$,
extended by $x R \infty$
for all $x \in \Q \cup \set{ \infty } \setminus \set{ 0 }$.
We have that $R$ is anti-symmetric, since it is so on $\Q$,
and $\infty R x$ doesn't hold, except for $x = \infty$.
$R^2$ is $(<) \cup \set{ \tpl{x,\infty} \mid 0 \neq x \neq \infty}$,
which is dense:
observe that $\Q \setminus \set{-1} \ni x \ro{R}{2} \infty$ has e.g.\ $x+1$
as intermediate element, similar for $x=-1$.
But $R^C$ is not dense,
since $0 \ro{R}{C} \infty$ has no intermediate element,
as $x \ro{R}{C} \infty$ implies $x = 0$,
but neither $0 \ro{R}{C} 0$ nor $\infty \ro{R}{C} \infty$ holds.
\qed
}

\EXAMPLE{aaxu}{%
2-Dense and 2-SemiOrd1 doesn't imply 2-SemiOrd2.
}
\PROOF{%
Consider $\Q \cup \set{\infty}$ with $R$ being
the usual order $(<)$,
and with $\infty$ incomparable to all elements.
Since this ordering is asymmetric, $\ro{R}{2}$ coincides with $R$.
Hence, it is dense, since it is dense on $\Q$,
and $\infty$ is not involved in any relation pair.
Moreover, $R$ satisfies semi-order property 1,
since $x \ro{R}{8} y \land y \rr{R} z$ 
can hold only if $x=y \neq \infty$;
in this case, $w \rr{R} x$ implies $w \rr{R} z$ by transitivity.
But $R$ doesn't satisfy semi-order property 2,
since $0 < 1 < 2$,
but $\infty$ is incomparable to all of them.
\qed
}

\EXAMPLE{i0iq}{%
2-Trans and 2-Dense doesn't imply 2-SemiOrd1.
}
\PROOF{%
Consider $\set{a,b} \times \Q$ with $R$ defined by
$\tpl{x,y} \rr{R} \tpl{u,v}$
iff $x=u$ and $y<v$, employing the usual order $<$.
This results in two independent copies of the rational numbers $\Q$.
Since this ordering is asymmetric, 2-Trans coincides with 3-Trans;
both properties are well-known to be satisfied by $<$.
Since no elements of different copies of $\Q$ are related,
$R$ itself is also both 2- and 3-Transitive, and moreover 2-Dense.
However, $R^2$ doesn't satisfy semi-order property 1,
since $\tpl{a,0} \ro{R}{2} \tpl{a,1}$,
and $\tpl{a,1}$ is incomparable w.r.t.\ $r^2$ to $\tpl{b,0}$,
and $\tpl{b,0} \ro{R}{2} \tpl{b,1}$,
but not $\tpl{a,0} \ro{R}{2} \tpl{b,1}$.
\qed
}

\EXAMPLE{jevo, jew0}{%
\LABEL{jevo jew0}%
D-Trans and 2-LfSerial implies neither 
1-Dense (\QDef{jevo}) nor 4-LfSerial (\QDef{jew0}).
}
\PROOF{
Consider the set $\N$ with
$x \rr{R} y$ defined as $x > y$ or $x=0 \land y=1$.
Then $\ro{R}{2}$ is left serial,
since e.g.\ $x+1 \ro{R}{2} x$ for $x>0$ and $2 \ro{R}{2} 0$.
The relation $\ro{R}{D}$ can be obtained as
$(\ro{R}{5}) \cup (\ro{R}{8})$;
with $x \ro{R}{5} y$ iff $x < y$ or $x=1 \land y=0$,
and $ x \ro{R}{8} y$ iff $x = y$;
hence $x \ro{R}{D} y$ iff $x \leq y$ or $x=1 \land y=0$.
Therefore, $\ro{R}{D}$ is transitive:
Let $x \ro{R}{D} y$ and $y \ro{R}{D} z$.
\begin{itemize}
\item If $x \leq y$ and $y \leq z$,
	then $x \leq z$, hence $x \ro{R}{D} z$.
\item If $x \leq y$ and $y=1 \land z=0$,
	then $x = 0$ or $x = 1$,
	in both cases, $x \ro{R}{D} 0$ holds.
\item If $x=1 \land y=0$ and $y \leq z$ and $z=0$,
	then $1 \ro{R}{D} 0$ holds.
\item If $x=1 \land y=0$ and $y \leq z$ and $z>0$,
	then $1 \leq z$, hence $x \ro{R}{D} z$.
\item The case $x=1 \land y=0$ and $y=1 \land z=0$ is a contradiction.
\end{itemize}
However, $\ro{R}{1} = \set{ \tpl{0,1}, \tpl{1,0} }$ is not dense.
Moreover, $\ro{R}{4}$ is not left serial,
i.e.\ $\ro{R}{2}$ is not right serial,
since $0 \ro{R}{2} x$ applies to no $x$.
\qed
}

\begin{figure}
\begin{center}
\begin{picture}(85,15)
\put( 2,8){\makebox(0,0)[bl]{$\cdots$}}
\put( 8,9){\vector(1,0){13}}
\put(22,9){\circle*{2}}
\put(22,11){\makebox(0,0)[b]{$5$}}
\put(23,9){\vector(1,0){13}}
\put(37,9){\circle*{2}}
\put(37,11){\makebox(0,0)[b]{$4$}}
\put(38,9){\vector(1,0){13}}
\put(52,9){\circle*{2}}
\put(52,11){\makebox(0,0)[b]{$3$}}
\put(53,9){\vector(1,0){13}}
\put(67,9){\circle*{2}}
\put(67,11){\makebox(0,0)[b]{$2$}}
\put(68,9){\vector(1,0){13}}
\put(82,9){\circle*{2}}
\put(82,11){\makebox(0,0)[b]{$1$}}
\put(68,1){\circle*{2}}
\put(70,1){\makebox(0,0)[l]{$0$}}
\put(37,9){\vector(4,-1){30}}
\end{picture}
\caption{Relation in Exm.~\REF{qjxu}}
\LABEL{Relation in Exm. qjxu}
\end{center}
\end{figure}
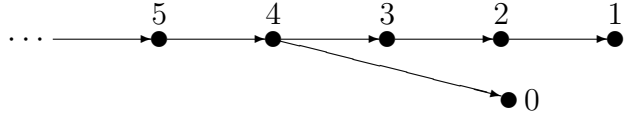

\EXAMPLE{qjxu}{%
\LABEL{qjxu}%
2-LfSerial and 2-SemiOrd1 doesn't imply 2-SemiOrd2.
}
\PROOF{%
Consider the universe set $\N$ and its usual order $(>)$,
define $R = (>) \setminus \set{ \tpl{1,0}, \tpl{2,0}, \tpl{3,0} }$,
cf.\ Fig.~\REF{Relation in Exm. qjxu}.
Since $(>)$, and therefore $R$, is asymmetric,
it coincides with $\ro{R}{2}$.
The latter is left serial since e.g.\
$i+1 \ro{R}{2} i$ for each $i \geq 1$,
and $4 \ro{R}{2} 0$.

$\ro{R}{2}$ also satisfies semi-order property 1:
Let $w \ro{R}{2} x$, $x,y$ incomparable w.r.t.\ $\ro{R}{2}$, and
$y \ro{R}{2} z$.
\begin{itemize}
\item If $x=y$ and $z > 0$,
	\\
	then $w > x = y > z > 0$,
	i.e.\ $0$ is not involved,
	hence $w \ro{R}{2} z$ by transitivity of $(>)$.
\item If $x=y$ and $z=0$,
	then $w > x = y \geq 4 > 0 = z$,
	hence again $w \ro{R}{2} z$ by transitivity.
\item The case $x \in \set{1,2,3}$ and $y=0$
	is impossible due to $y \ro{R}{2} z$.
\item If $x = 0$ and $y \in \set{1,2,3}$,
	\\
	then $y \ro{R}{2} z$ implies $3 \geq y > z > 0$,
	and $w \ro{R}{2} 0$ implies $w \geq 4$,
	hence $w \ro{R}{2} z$.
\end{itemize}
However, $\ro{R}{2}$ doesn't satisfy semi-order property 2,
since
$0$ is incomparable to all elements involved in
$3 \ro{R}{2} 2 \ro{R}{2} 1$.
\qed
}

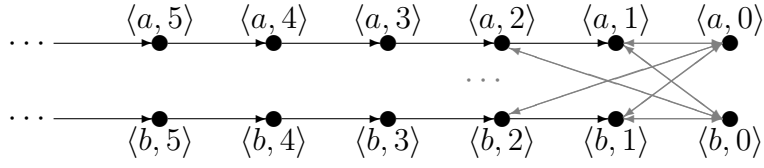
\begin{figure}
\begin{center}
\begin{picture}(100,25)
\put( 0, 6){\makebox(0,0)[bl]{$\cdots$}}
\put( 6, 7){\vector(1,0){13}}
\put(20, 7){\circle*{2}}
\put(20, 6){\makebox(0,0)[t]{$\tpl{b,5}$}}
\put(21, 7){\vector(1,0){13}}
\put(35, 7){\circle*{2}}
\put(35, 6){\makebox(0,0)[t]{$\tpl{b,4}$}}
\put(36, 7){\vector(1,0){13}}
\put(50, 7){\circle*{2}}
\put(50, 6){\makebox(0,0)[t]{$\tpl{b,3}$}}
\put(51, 7){\vector(1,0){13}}
\put(65, 7){\circle*{2}}
\put(65, 6){\makebox(0,0)[t]{$\tpl{b,2}$}}
\put(66, 7){\vector(1,0){13}}
\put(80, 7){\circle*{2}}
\put(80, 6){\makebox(0,0)[t]{$\tpl{b,1}$}}
\textcolor{coLowlight}{\put(81, 7){\vector(1,0){13}}}%
\put(95, 7){\circle*{2}}
\put(95, 6){\makebox(0,0)[t]{$\tpl{b,0}$}}
\textcolor{coLowlight}{\put(94,7){\vector(-1,0){13}}}%
\textcolor{coLowlight}{\put(94.000,7.000){\vector(-4,3){13.000}}}%
\textcolor{coLowlight}{\put(81.000,16.750){\vector(4,-3){13.000}}}%
\textcolor{coLowlight}{\put(94.000,7.000){\vector(-3,1){28.000}}}%
\textcolor{coLowlight}{\put(66.000,16.333){\vector(3,-1){28.000}}}%
\put( 0,16){\makebox(0,0)[bl]{$\cdots$}}
\put( 6,17){\vector(1,0){13}}
\put(20,17){\circle*{2}}
\put(20,18){\makebox(0,0)[b]{$\tpl{a,5}$}}
\put(21,17){\vector(1,0){13}}
\put(35,17){\circle*{2}}
\put(35,18){\makebox(0,0)[b]{$\tpl{a,4}$}}
\put(36,17){\vector(1,0){13}}
\put(50,17){\circle*{2}}
\put(50,18){\makebox(0,0)[b]{$\tpl{a,3}$}}
\put(51,17){\vector(1,0){13}}
\put(65,17){\circle*{2}}
\put(65,18){\makebox(0,0)[b]{$\tpl{a,2}$}}
\put(66,17){\vector(1,0){13}}
\put(80,17){\circle*{2}}
\put(80,18){\makebox(0,0)[b]{$\tpl{a,1}$}}
\textcolor{coLowlight}{\put(81,17){\vector(1,0){13}}}%
\put(95,17){\circle*{2}}
\put(95,18){\makebox(0,0)[b]{$\tpl{a,0}$}}
\textcolor{coLowlight}{\put(94,17){\vector(-1,0){13}}}%
\textcolor{coLowlight}{\put(94,17.000){\vector(-4,-3){13.000}}}%
\textcolor{coLowlight}{\put(81.000,7.250){\vector(4,3){13.000}}}%
\textcolor{coLowlight}{\put(94,17.000){\vector(-3,-1){28.000}}}%
\textcolor{coLowlight}{\put(66.000,7.667){\vector(3,1){28.000}}}%
\textcolor{coLowlight}{\put(65,12){\makebox(0,0)[r]{$\cdots$}}}%
\end{picture}
\caption{Relation in Exm.~\REF{Lem qlas}}
\LABEL{Relation in Exm. qlas}
\end{center}
\end{figure}

\EXAMPLE{qlas}{%
\LABEL{Lem qlas}%
2-LfSerial and 2-Trans doesn't imply E-Dense.
}
\PROOF{
Consider the universe set $\set{a,b} \times \N$,
and the relation $R$ defined by
$\tpl{x,y} \rr{R} \tpl{u,v}$
iff $x=u \land y > v$
or $y=0 \land v \neq 0$
or $y \neq 0 \land v=0$.
Figure~\REF{Relation in Exm. qlas}
sketches the relation:
it consists of two copies of $(>)$ on $\N$,
but additionally relates each zero to
each non-zero universe member, and vice versa.
We have $\tpl{x,y} \ro{R}{2} \tpl{u,v}$
iff
$x=u \land y > v > 0$;
this relation is ls and tr.
$R$ is not E-de since $\tpl{a,0} \ro{R}{E} \tpl{b,0}$
has no intermediate element by construction.
\qed
}

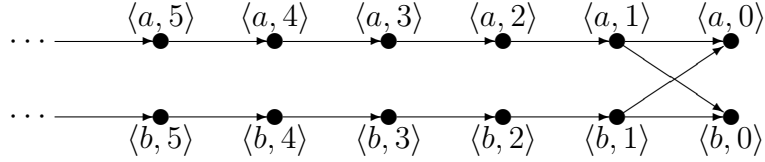
\begin{figure}
\begin{center}
\begin{picture}(110,25)%
\put( 0, 6){\makebox(0,0)[bl]{$\cdots$}}%
\put( 6, 7){\vector(1,0){13}}%
\put(20, 7){\circle*{2}}%
\put(20, 6){\makebox(0,0)[t]{$\tpl{b,5}$}}%
\put(21, 7){\vector(1,0){13}}%
\put(35, 7){\circle*{2}}%
\put(35, 6){\makebox(0,0)[t]{$\tpl{b,4}$}}%
\put(36, 7){\vector(1,0){13}}%
\put(50, 7){\circle*{2}}%
\put(50, 6){\makebox(0,0)[t]{$\tpl{b,3}$}}%
\put(51, 7){\vector(1,0){13}}%
\put(65, 7){\circle*{2}}%
\put(65, 6){\makebox(0,0)[t]{$\tpl{b,2}$}}%
\put(66, 7){\vector(1,0){13}}%
\put(80, 7){\circle*{2}}%
\put(80, 6){\makebox(0,0)[t]{$\tpl{b,1}$}}%
\put(81, 7){\vector(1,0){13}}%
\put(80, 7){\vector(3,2){14}}%
\put(95, 7){\circle*{2}}%
\put(95, 6){\makebox(0,0)[t]{$\tpl{b,0}$}}%
\put( 0,16){\makebox(0,0)[bl]{$\cdots$}}%
\put( 6,17){\vector(1,0){13}}%
\put(20,17){\circle*{2}}%
\put(20,18){\makebox(0,0)[b]{$\tpl{a,5}$}}%
\put(21,17){\vector(1,0){13}}%
\put(35,17){\circle*{2}}%
\put(35,18){\makebox(0,0)[b]{$\tpl{a,4}$}}%
\put(36,17){\vector(1,0){13}}%
\put(50,17){\circle*{2}}%
\put(50,18){\makebox(0,0)[b]{$\tpl{a,3}$}}%
\put(51,17){\vector(1,0){13}}%
\put(65,17){\circle*{2}}%
\put(65,18){\makebox(0,0)[b]{$\tpl{a,2}$}}%
\put(66,17){\vector(1,0){13}}%
\put(80,17){\circle*{2}}%
\put(80,18){\makebox(0,0)[b]{$\tpl{a,1}$}}%
\put(81,17){\vector(1,0){13}}%
\put(80,17){\vector(3,-2){14}}%
\put(95,17){\circle*{2}}%
\put(95,18){\makebox(0,0)[b]{$\tpl{a,0}$}}%
\end{picture}
\caption{Relation in Exm.~\REF{qldr}}
\LABEL{Relation in Exm. qldr}
\end{center}
\end{figure}

\EXAMPLE{qldr}{%
\LABEL{qldr}%
2-LfSerial and 3-Trans doesn't imply 8-Dense.
}
\PROOF{
On the set $\set{a,b} \times \N$,
define $R'$ such that for all $x,y,u$:
\begin{enumerate}
\item $\tpl{x,y\.+1} \rr{R'} \tpl{x,y}$,
\item $\tpl{x,1} \rr{R'} \tpl{u,0}$,
\item $\tpl{x,0} \rr{R'} \tpl{x,0}$,
\item $R'$ applies to no other pairs than given above.
\end{enumerate}
Figure~\REF{Relation in Exm. qldr} sketches that relation;
it is left serial due to rule 1.
Let $R$ be the transitive closure of $R'$,
then $R$ is transitive by construction.
$\ro{R}{2}$ is obtained from $R$ by removing the pairs
$\tpl{ \tpl{a,0}, \tpl{a,0} }$ and
$\tpl{ \tpl{b,0}, \tpl{b,0} }$;
it is left serial, since $R'$ is.
However, $\ro{R}{8}$ is not dense:
$\tpl{a,0} \ro{R}{8} \tpl{b,0}$ has no intermediate element,
since $\tpl{a,0}$ is comparable w.r.t.\ $R$ to each
element except $\tpl{b,0}$, and similar for $\tpl{b,0}$.
\qed
}

\begin{figure}
\begin{center}
\begin{picture}(110,30)
\put( 0, 1){\makebox(0,0)[bl]{$\cdots$}}
\put( 6, 2){\vector(1,0){13}}
\put(20, 2){\circle*{2}}
\put(20, 3){\makebox(0,0)[b]{$\tpl{c,5}$}}
\put(21, 2){\vector(1,0){13}}
\put(35, 2){\circle*{2}}
\put(35, 3){\makebox(0,0)[b]{$\tpl{c,4}$}}
\put(36, 2){\vector(1,0){13}}
\put(50, 2){\circle*{2}}
\put(50, 3){\makebox(0,0)[b]{$\tpl{c,3}$}}
\put(51, 2){\vector(1,0){13}}
\put(65, 2){\circle*{2}}
\put(65, 3){\makebox(0,0)[b]{$\tpl{c,2}$}}
\put(66, 2){\vector(1,0){13}}
\put(80, 2){\circle*{2}}
\put(80, 3){\makebox(0,0)[b]{$\tpl{c,1}$}}
\put(81, 2){\vector(1,0){13}}
\put(95, 2){\circle*{2}}
\put(96, 2){\makebox(0,0)[l]{$\tpl{c,0}$}}
\put( 0,11){\makebox(0,0)[bl]{$\cdots$}}
\put( 6,12){\vector(1,0){13}}
\put(20,12){\circle*{2}}
\put(20,13){\makebox(0,0)[b]{$\tpl{b,4}$}}
\put(21,12){\vector(1,0){13}}
\put(35,12){\circle*{2}}
\put(35,13){\makebox(0,0)[b]{$\tpl{b,3}$}}
\put(36,12){\vector(1,0){13}}
\put(50,12){\circle*{2}}
\put(50,13){\makebox(0,0)[b]{$\tpl{b,2}$}}
\put(51,12){\vector(1,0){13}}
\put(65,12){\circle*{2}}
\put(65,13){\makebox(0,0)[b]{$\tpl{b,1}$}}
\put(66,12){\vector(1,0){13}}
\put(80,12){\circle*{2}}
\put(83,12){\makebox(0,0)[l]{$\tpl{b,0}$}}
\put( 0,21){\makebox(0,0)[bl]{$\cdots$}}
\put( 6,22){\vector(1,0){13}}
\put(20,22){\circle*{2}}
\put(20,23){\makebox(0,0)[b]{$\tpl{a,5}$}}
\put(21,22){\vector(1,0){13}}
\put(35,22){\circle*{2}}
\put(35,23){\makebox(0,0)[b]{$\tpl{a,4}$}}
\put(36,22){\vector(1,0){13}}
\put(50,22){\circle*{2}}
\put(50,23){\makebox(0,0)[b]{$\tpl{a,3}$}}
\put(51,22){\vector(1,0){13}}
\put(65,22){\circle*{2}}
\put(65,23){\makebox(0,0)[b]{$\tpl{a,2}$}}
\put(66,22){\vector(1,0){13}}
\put(80,22){\circle*{2}}
\put(80,23){\makebox(0,0)[b]{$\tpl{a,1}$}}
\put(81,22){\vector(1,0){13}}
\put(95,22){\circle*{2}}
\put(96,22){\makebox(0,0)[l]{$\tpl{a,0}$}}
\put(81,13){\vector(3,2){13}}
\put(81,11){\vector(3,-2){13}}
\put(94,21.667){\vector(-3,-2){13}}
\put(94,2.333){\vector(-3,2){13}}
\end{picture}
\caption{Relation in Exm.~\REF{Lem qmkh}}
\LABEL{Relation in Exm. qmkh}
\end{center}
\end{figure}
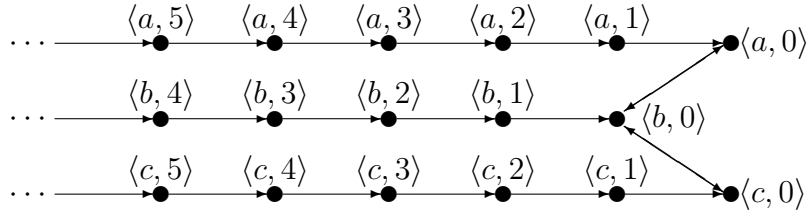

\EXAMPLE{qmkh}{%
\LABEL{Lem qmkh}%
2-LfSerial and 7-AntiTrans doesn't imply 1-LfUnique.
}
\PROOF{%
Consider the universe set $\set{a,b,c} \times \N$,
and the relation $R$ defined as follows:
$\tpl{x,y} \rr{R} \tpl{u,v}$
iff $x=u \land y = v+1$
or $y=v=0 \land x \neq u \land b \in \set{x,u}$.
The universe consists of three copies of $\N$, connected only at $0$;
see Fig.~\REF{Relation in Exm. qmkh}.
The relation $R$ is 2-LfSerial, 
since $\tpl{x,i+1} \rr{R} \tpl{x,i}$,
but not $\tpl{x,i} \rr{R} \tpl{x,i+1}$
for each $x \in \set{a,b,c}$ and $i \in \N$.
It is 7-AntiTrans, since it doesn't contain
an undirected cycle of length $3$.
In fact, the only cycles are $\tpl{b,0} \rr{R} \tpl{a,0}
\rr{R} \tpl{b,0}$ and
$\tpl{b,0} \rr{R} \tpl{c,0} \rr{R} \tpl{b,0}$, 
and compositions thereof, all of which have an even length.
However, it is not 1-LfUnique, since $\tpl{b,0} \ro{R}{1} \tpl{a,0}$
and $\tpl{b,0} \ro{R}{1} \tpl{c,0}$.
\qed
}

\EXAMPLE{qpkw}{%
\LABEL{qpkw}
2-LfSerial and 5-LfUnique doesn't imply 8-LfQuasiRefl.
}
\PROOF{
Consider the set $\N$ with
$x \rr{R} y$ defined as $x = y+1$ or $x=y=0$.
Then $\ro{R}{2}$ is left serial, since $x+1 \ro{R}{2} x$
for all $x$.
Moreover, $R$ is right-unique,
since a successor number $y+1$ is related only to $y$,
and $0$ is only related to $0$.
However, $\ro{R}{8}$ is not left quasi-reflexive, since e.g.\
$0 \ro{R}{8} 2$, but not $0 \ro{R}{8} 0$.
\qed
}

\begin{figure}
\begin{center}
\begin{picture}(55,25)%
\put(2,20){\makebox(0,0){$1$}}
\put(3,20){\vector(1,0){10}}
\put(3,20){\vector(1,-1){10}}
\put(14,20){\makebox(0,0){$2$}}
\put(14,10){\makebox(0,0){$3$}}
\put(15,20){\vector(1,0){10}}
\put(15,10){\vector(1,0){10}}
\put(15,20){\vector(2,-1){10}}
\put(15,10){\vector(2,-1){10}}
\put(26,20){\makebox(0,0){$4$}}
\put(26,15){\makebox(0,0){$5$}}
\put(26,10){\makebox(0,0){$6$}}
\put(26,5){\makebox(0,0){$7$}}
\put(27,20){\vector(1,0){10}}
\put(27,5){\vector(4,-1){10}}
\put(37,12){\makebox(0,0)[r]{$\ldots$}}
\end{picture}
\caption{Relation in Exm.~\REF{ufbk}}
\LABEL{Relation in Exm. ufbk}
\end{center}
\end{figure}
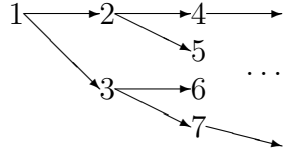

\EXAMPLE{ufbk}{%
\LABEL{ufbk}%
2-LfUnique and 4-LfSerial does not imply 4-LfUnique
}
\PROOF{
Consider $\N \setminus \set{0}$ and the relation
$xRy :\Lra x = y \idiv 2$, where $\idiv$ denotes truncating
integer division, see Fig.~\REF{Relation in Exm. ufbk}.
It is 2-LfUnique, since $x_1 = y \idiv 2 = x_2$ 
implies $x_1 = x_2$.
It is 4-LfSerial, since for a given $y$, 
choosing $x = y \cdot 2$ will satisfy
$x \neq y \idiv 2$ and $y = x \idiv 2$,
that is, $x \ro{R}{4} y$.
However, it is not 4-LfUnique, since e.g.\
$1 = 2 \idiv 2$ and $1 = 3 \idiv 2$,
that is, $2 R^4 1$ and $3 R^4 1$.
\qed
}

\begin{figure}
\begin{center}
\begin{picture}(110,25)%
\put(100.000,6.000){\makebox(0.000,0.000)[bl]{$\cdots$}}%
\put(81.000,7.000){\vector(1,0){13.000}}%
\put(80.000,7.000){\circle*{2.000}}%
\put(80.000,6.000){\makebox(0.000,0.000)[t]{$\tpl{b,5}$}}%
\put(66.000,7.000){\vector(1,0){13.000}}%
\put(65.000,7.000){\circle*{2.000}}%
\put(65.000,6.000){\makebox(0.000,0.000)[t]{$\tpl{b,4}$}}%
\put(51.000,7.000){\vector(1,0){13.000}}%
\put(50.000,7.000){\circle*{2.000}}%
\put(50.000,6.000){\makebox(0.000,0.000)[t]{$\tpl{b,3}$}}%
\put(36.000,7.000){\vector(1,0){13.000}}%
\put(35.000,7.000){\circle*{2.000}}%
\put(35.000,6.000){\makebox(0.000,0.000)[t]{$\tpl{b,2}$}}%
\put(21.000,7.000){\vector(1,0){13.000}}%
\put(20.000,7.000){\circle*{2.000}}%
\put(20.000,6.000){\makebox(0.000,0.000)[t]{$\tpl{b,1}$}}%
\put(6.000,7.000){\vector(1,0){13.000}}%
\put(5.000,7.000){\circle*{2.000}}%
\put(5.000,6.000){\makebox(0.000,0.000)[t]{$\tpl{b,0}$}}%
\put(100.000,16.000){\makebox(0.000,0.000)[bl]{$\cdots$}}%
\put(81.000,17.000){\vector(1,0){13.000}}%
\put(80.000,17.000){\circle*{2.000}}%
\put(80.000,18.000){\makebox(0.000,0.000)[b]{$\tpl{a,5}$}}%
\put(66.000,17.000){\vector(1,0){13.000}}%
\put(65.000,17.000){\circle*{2.000}}%
\put(65.000,18.000){\makebox(0.000,0.000)[b]{$\tpl{a,4}$}}%
\put(51.000,17.000){\vector(1,0){13.000}}%
\put(50.000,17.000){\circle*{2.000}}%
\put(50.000,18.000){\makebox(0.000,0.000)[b]{$\tpl{a,3}$}}%
\put(36.000,17.000){\vector(1,0){13.000}}%
\put(35.000,17.000){\circle*{2.000}}%
\put(35.000,18.000){\makebox(0.000,0.000)[b]{$\tpl{a,2}$}}%
\put(21.000,17.000){\vector(1,0){13.000}}%
\put(20.000,17.000){\circle*{2.000}}%
\put(20.000,18.000){\makebox(0.000,0.000)[b]{$\tpl{a,1}$}}%
\put(6.000,17.000){\vector(1,0){13.000}}%
\put(5.000,17.000){\circle*{2.000}}%
\put(5.000,18.000){\makebox(0.000,0.000)[b]{$\tpl{a,0}$}}%
\put(5.000,7.000){\vector(0,1){9.000}}%
\put(5.000,17.000){\vector(0,-1){9.000}}%
\end{picture}
\caption{Relation in Exm.~\REF{uo0e}}
\LABEL{Relation in Exm. uo0e}
\end{center}
\end{figure}
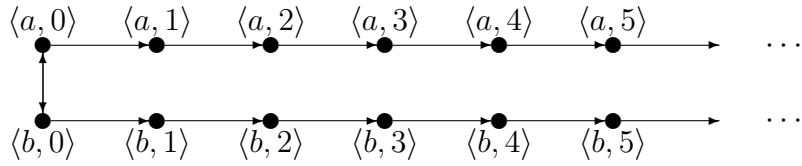

\EXAMPLE{uo0e}{%
\LABEL{uo0e}%
3-LfUnique and 4-LfSerial doesn't imply 1-Trans.
}
\PROOF{
On the set $\set{a,b} \times \N$,
define $R$ such that:
\begin{enumerate}
\item $\tpl{x,y} \rr{R} \tpl{x,y\.+1}$,
\item $\tpl{x,0} \rr{R} \tpl{u,0}$ if $x \neq u$,
\item $R$ applies to no other pairs than given above.
\end{enumerate}
Figure~\REF{Relation in Exm. uo0e} sketches $R$
and allows us to verify that
$R$ is left-unique (no two arrows end at the same point)
and $\ro{R}{4}$ is left serial (a one-sided arrow starts from each point).
However, $\ro{R}{1}$ is not transitive, since
$\tpl{a,0} \ro{R}{1} \tpl{b,0}$, and
$\tpl{b,0} \ro{R}{1} \tpl{a,0}$, but not
$\tpl{a,0} \ro{R}{1} \tpl{a,0}$.
\qed
}

\subsection{Reported laws}
\LABEL{Reported laws}

After several attempts, we eventually achieved to decide for every
$3$-implication its truth value.
We found $156 384$ valid and $375 057$ invalid implications; the
former include $42 657$ from trivial initialization.

In the ancillary file \filename{lpImpl.log}
each implication can be found together with its truth value and a
justification.
The executable \filename{justify} can be used to obtain the truth
values (and justifications) of given implications.
Figure~\REF{Proof tree of ``tr land ir ra as''} 
and~\REF{Disproof tree of ``as land tr not ra rf''}
show examples;
note that we don't print proof parts contributed by the external 
prover or the counter-example generator.

Implications with one antecedent\footnote{
	in our approach,
	they appear 
	as $3$-implications of the form e.g.\ $I \land I \ra K$
}
can be visualized as a directed
graph (Hasse diagram),
with lifted properties as vertices and implications as edges.
This has been done for basic properties in
\OEF{Fig.18}{20}{}.
For lifted properties, the graph
is still too large to be depicted in one piece.
We have split it by basic property and show the pieces in
Fig.~\REF{Implications about anti-symmetry (an)}
to~\REF{Implications about symmetry (sy), asymmetry (as), and reflexivity (rf)}.
In each figure, vertices and edges for the depicted basic
property are shown
in blue, while others, for implied or implying properties, are shown in
red.
We have two exceptions: 
edges contained in our minimal
axiom set (Sect.~\REF{Axiom set} below) are colored green,
and 
edges following from monotonicity or
antitonicity are colored cyan (for the depicted basic property)
or magenta (for others).
For brevity,
the two-character codes from Def.~\REF{def} are used, and the
separating hyphen between operation and property name is omitted.
For example, the bottommost vertex in
Fig.~\REF{Implications about anti-symmetry (an)}
represents the lifted property 9-AntiSym.
Since the vertices represent equivalence classes,
all implications are strict.

Figure~\REF{Implications about basic properties}
summarizes the part about all basic properties;
co-transitivity
is included (as ``Ctr'') although we didn't define it as a basic
property.
In this figure, in order to obtain a feasible complexity,
a red vertex is shown only if it is both implied by a blue vertex
and implying another one.

Figure~\REF{Inconsistencies between lifted properties}
shows the inconsistencies\footnote{
	i.e.\ $3$-implications of the form $I \land J \ra$ 0-Refl
}
between lifted properties.
Clusters of incompatibilities are presented using disjunctions of
lifted properties.
For example, the topmost horizontal red edge, from ``Alu,Clu,Cat'' to
``1at,1lu'', indicates that each property from the set
$\set{ \mbox{A-LfUnique}, \mbox{C-LfUnique}, \mbox{C-AntiTrans} }$
is inconsistent
with each property from
$\set{ \mbox{1-AntiTrans}, \mbox{1-LfUnique} }$.
Nested boxes indicate implications between properties, that is,
subset relations between their extensions.
For example, the ``as'' box is inside the ``1at'' box since
each ASym relation is 1-AntiTrans (\QRef{cjdj}).
The red lines' end points should be considered carefully;
for example, the topmost vertical line connects ``1ls'' with ``as''
but not with the surrounding ``1at'' box,
since 1-LfSerial is inconsistent with ASym, but not
necessarily with 1-AntiTrans.
All found inconsistencies are covered in
Fig.~\REF{Inconsistencies between lifted properties};
however, some of them require application of some implications, viz.\
$$\begin{tabular}{r@{$\;$}c@{$\;$}r@{$\;$}c@{$\;$}r@{$\;$}c@{$\;$}rl}
6-at &$\Ra$ & 2-at &      &     &     &      &(\QDef{lvfx}),	\\
E-an &$\Ra$ & C-an &$\Lra$& sc  &$\Ra$& 8-lu &(\QDef{ld0m}, \QDef{djno}),\\
E-at &$\Ra$ & C-at &      &     &     &      &(\QDef{mvoz}),	\\
cr   &$\Lra$& 7-an &$\Ra$ & an  &$\Ra$& 1-lu &(\QDef{kuxl}, \QDef{dakh}),\\
cr   &$\Lra$& 7-an &$\Ra$ & sy  &     &      &(\QDef{kuxb}),	\\
cr   &$\Lra$& 7-an &$\Ra$ & 7-lu&     &      &(\QDef{kuzn}),	\\
E-an &$\Ra$ & E-lu &      &     &     &      &(\QDef{ldbr}), and	\\
E-an &$\Ra$ & sy   &      &     &     &      &(\QDef{ld0b}).	\\
\end{tabular}$$
For example, the inconsistency of SemiConnex and CoRefl known from
\OEF{Lem.8.8}{23}{} is tacitly understood as a
consequence of the inconsistency of 1-LfUnique and 7-LfUnique shown
at the bottom of Fig.~\REF{Inconsistencies between lifted properties}.

Together, Fig.~\REF{Implications about anti-symmetry (an)}
to~\REF{Inconsistencies between lifted properties}
show a refinement of the graph from
\OEF{Fig.18}{20}{}.

\newcommand{\lpImplIScale}{0.75}

\begin{figure}
\begin{center}
\includegraphics[scale=\lpImplIScale]{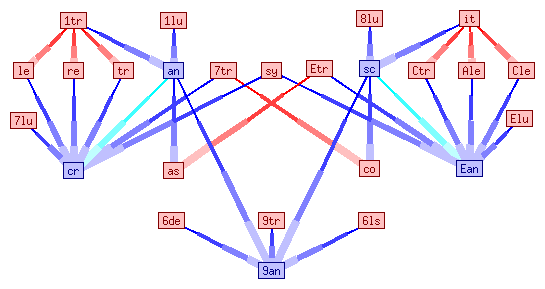}
\caption{Implications about anti-symmetry (an)}
\LABEL{Implications about anti-symmetry (an)}
\end{center}
\end{figure}

\begin{figure}
\begin{center}
\includegraphics[scale=\lpImplIScale]{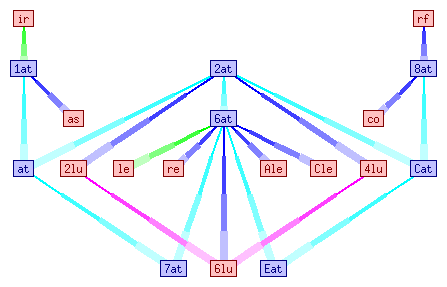}
\caption{Implications about anti-transitivity (at)}
\LABEL{Implications about transitivity-symmetry (at)}
\end{center}
\end{figure}

\begin{figure}
\begin{center}
\includegraphics[scale=\lpImplIScale]{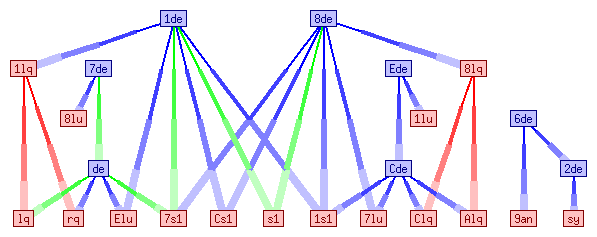}
\caption{Implications about density (de)}
\LABEL{Implications about density (de)}
\end{center}
\end{figure}

\begin{figure}
\begin{center}
\includegraphics[scale=\lpImplIScale]{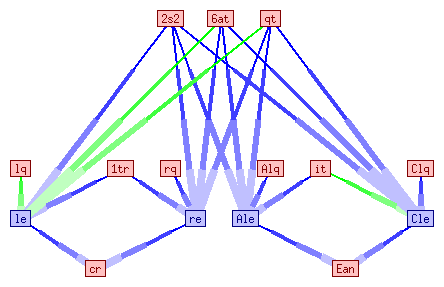}
\caption{Implications about left Euclideanness (le)}
\LABEL{Implications about left Euclideanness (le)}
\end{center}
\end{figure}

\begin{figure}
\begin{center}
\includegraphics[scale=\lpImplIScale]{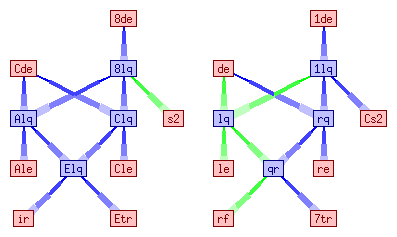}
\caption{Implications about left quasi-reflexivity (lq)}
\LABEL{Implications about left quasi-reflexivity (lq)}
\end{center}
\end{figure}

\begin{figure}
\begin{center}
\includegraphics[scale=\lpImplIScale]{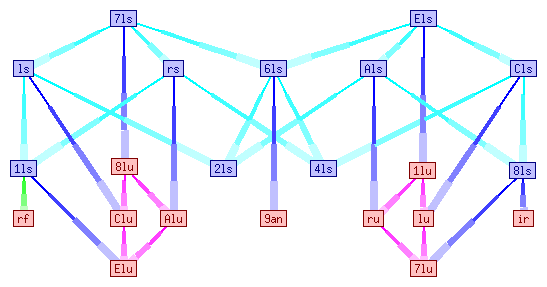}
\caption{Implications about left seriality (ls)}
\LABEL{Implications about left seriality (ls)}
\end{center}
\end{figure}

\begin{figure}
\begin{center}
\includegraphics[scale=\lpImplIScale]{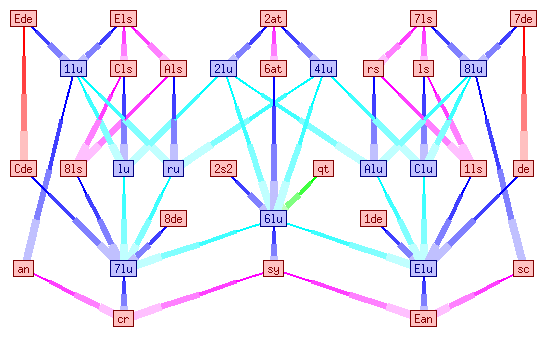}
\caption{Implications about left uniqueness (lu)}
\LABEL{Implications about left uniqueness (lu)}
\end{center}
\end{figure}

\begin{figure}
\begin{center}
\includegraphics[scale=\lpImplIScale]{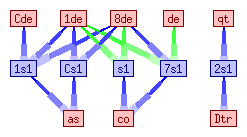}
\caption{Implications about semi-order property 1 (s1)}
\LABEL{Implications about semi-order property 1 (s1)}
\end{center}
\end{figure}

\begin{figure}
\begin{center}
\includegraphics[scale=\lpImplIScale]{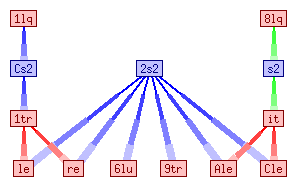}
\caption{Implications about semi-order property 2 (s2)}
\LABEL{Implications about semi-order property 2 (s2)}
\end{center}
\end{figure}

\begin{figure}
\begin{center}
\includegraphics[scale=\lpImplIScale]{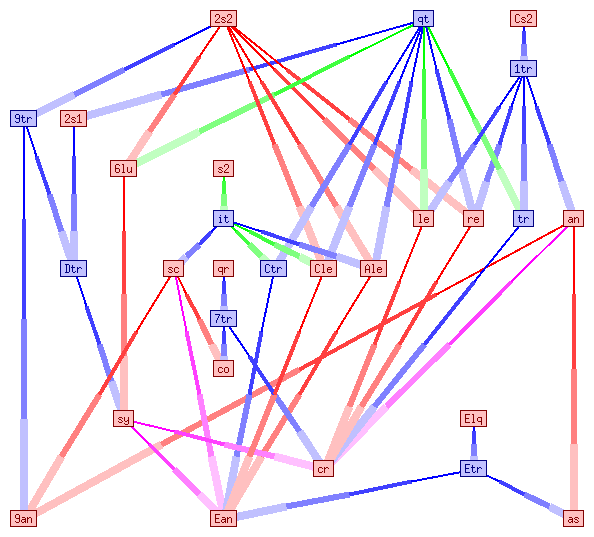}
\caption{Implications about transitivity (tr)}
\LABEL{Implications about transitivity (tr)}
\end{center}
\end{figure}

\begin{figure}
\begin{center}
\includegraphics[scale=\lpImplIScale]{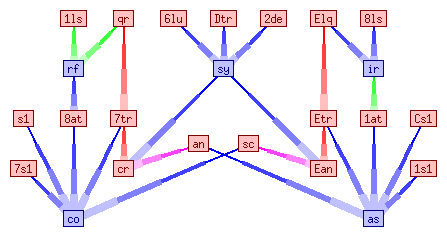}
\caption{Implications about
	symmetry (sy), asymmetry (as), and reflexivity (rf)}
\LABEL{Implications about symmetry (sy), asymmetry (as), and reflexivity (rf)}
\end{center}
\end{figure}

\begin{figure}
\begin{center}
\includegraphics[scale=\lpImplIScale]{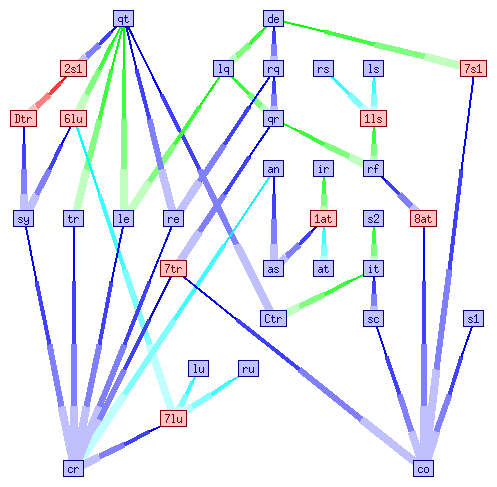}
\caption{Implications about basic properties}
\LABEL{Implications about basic properties}
\end{center}
\end{figure}

\definecolor{bSetII}	{rgb}{0.95,0.95,0.95}
\definecolor{bSetI}	{rgb}{0.90,0.90,0.90}
\definecolor{bSetZ}	{rgb}{0.85,0.85,0.85}
\definecolor{fInco}	{rgb}{0.99,0.00,0.00}

\begin{figure}
\newcommand{\7}[2]{%
	\textcolor{white}{\put(#1,#2){\makebox(0,0)[bl]{\rule{2mm}{2mm}}}}%
}
\newcommand{\8}{\rule[-0.5mm]{0mm}{0mm}}
\newcommand{\9}[7]{%
	\textcolor{#5}{\put(#1,#2){\makebox(0,0)[bl]{\rule{#3mm}{#4mm}}}}%
	\put(#1,#2){\framebox(#3,#4)#6{#7}}%
}
\newcommand{\3}[3]{\9{#1}{#2}{35}{20}{bSetII}{[bl]}{\8~#3}}%
\newcommand{\2}[3]{\9{#1}{#2}{30}{20}{bSetII}{[bl]}{\8~#3}}%
\newcommand{\1}[3]{\9{#1}{#2}{20}{12}{bSetI} {[bl]}{\8~#3}}%
\newcommand{\0}[3]{\9{#1}{#2}{10}{ 5}{bSetZ} {}    {#3}}%
\begin{center}
\begin{picture}(160,165)%
\2{ 85}{130}{Alu,Clu,Cat}%
\0{ 87}{140}{Elu}%
\0{ 103}{140}{Cat}%
\2{130}{130}{1at,1lu}%
\1{132}{135}{1at}%
\0{134}{140}{as}%
\3{ 45}{ 60}{6lu,le,re,Ale,Cle}%
\0{ 47}{ 70}{sy}%
\2{  5}{ 60}{6ls}%
\1{  7}{ 65}{2ls,4ls}%
\2{100}{ 60}{at,lu,ru}%
\0{102}{ 70}{at}%
\2{ 45}{  0}{7at,7lu}%
\0{ 47}{ 10}{7lu}%
\0{ 63}{ 10}{7at}%
\2{130}{  0}{8at,8lu}%
\1{132}{  5}{8at}%
\0{134}{ 10}{co}%
\0{ 65}{120}{9an}%
\0{ 45}{113}{2lu}%
\0{ 45}{120}{4lu}%
\0{ 45}{127}{2at}%
\0{130}{ 40}{rf}%
\0{130}{100}{ir}%
\0{134}{160}{1ls}%
\0{110}{ 10}{8ls}%
\thicklines%
\color{fInco}%
\put( 35, 72){\line( 1, 0){ 12}}
\put( 27, 70){\line( 1, 0){ 18}}
\put( 70,120){\line( 0,-1){ 40}}
\put( 65,122){\line(-1, 0){ 10}}
\put( 65,122){\line(-3, 2){ 10}}
\put( 65,122){\line(-3,-2){ 10}}
\put(120, 12){\line( 1, 0){ 14}}
\put( 75,  5){\line( 1, 0){ 55}}
\put(135, 45){\line( 0, 1){ 55}}
\put(139,145){\line( 0, 1){ 15}}
\put(130,143){\line(-1, 0){ 15}}
\put(150, 17){\line( 0, 1){118}}
\put(113,130){\line( 0,-1){ 50}}
\put( 87,141){\line(-1, 0){ 72}}
\put( 15,141){\line( 0,-1){ 64}}
\put( 87,143){\line(-1, 0){ 87}}
\put( 0,143){\line( 0,-1){132}}
\put( 0, 11){\line( 1, 0){ 47}}
\put( 15, 65){\line( 0,-1){ 52}}
\put( 15, 13){\line( 1, 0){ 32}}
\put( 70,125){\line( 0, 1){ 10}}
\put( 70,135){\line( 1, 0){ 15}}
\put(133, 45){\line( 0, 1){ 27}}
\put(133, 72){\line(-1, 0){ 21}}
\put(139, 45){\line( 1, 6){  9}}
\put(148, 99){\line( 0, 1){ 36}}
\7{142.5}{70.8}%
\7{142.5}{71.8}%
\7{142.5}{72.8}%
\put(148, 17){\line( 0, 1){ 29}}
\put(148, 46){\line(-1, 6){  9}}
\put(130,102){\line(-1, 0){ 10}}
\put(120,102){\line( 0, 1){ 39}}
\put(120,141){\line(-1, 0){ 7}}
\put( 75,121){\line( 1, 0){ 10}}
\put( 85,121){\line( 0,-1){105}}
\put( 85, 16){\line(-1, 0){ 10}}
\put( 75,123){\line( 1, 0){ 28}}
\put(103,123){\line( 0,-1){ 43}}
\7{ 94}{122}%
\put( 95,130){\line( 0,-1){116}}
\put( 95, 14){\line(-1, 0){ 20}}
\put(130, 42){\line(-1, 0){ 25}}
\put(105, 42){\line( 0,-1){ 30}}
\put(105, 12){\line(-1, 0){ 32}}
\7{124}{41}%
\put(125, 60){\line( 0,-1){ 42}}
\put(125, 18){\line( 1, 0){  5}}
\end{picture}
\caption{Inconsistencies between lifted properties}
\LABEL{Inconsistencies between lifted properties}
\end{center}
\end{figure}
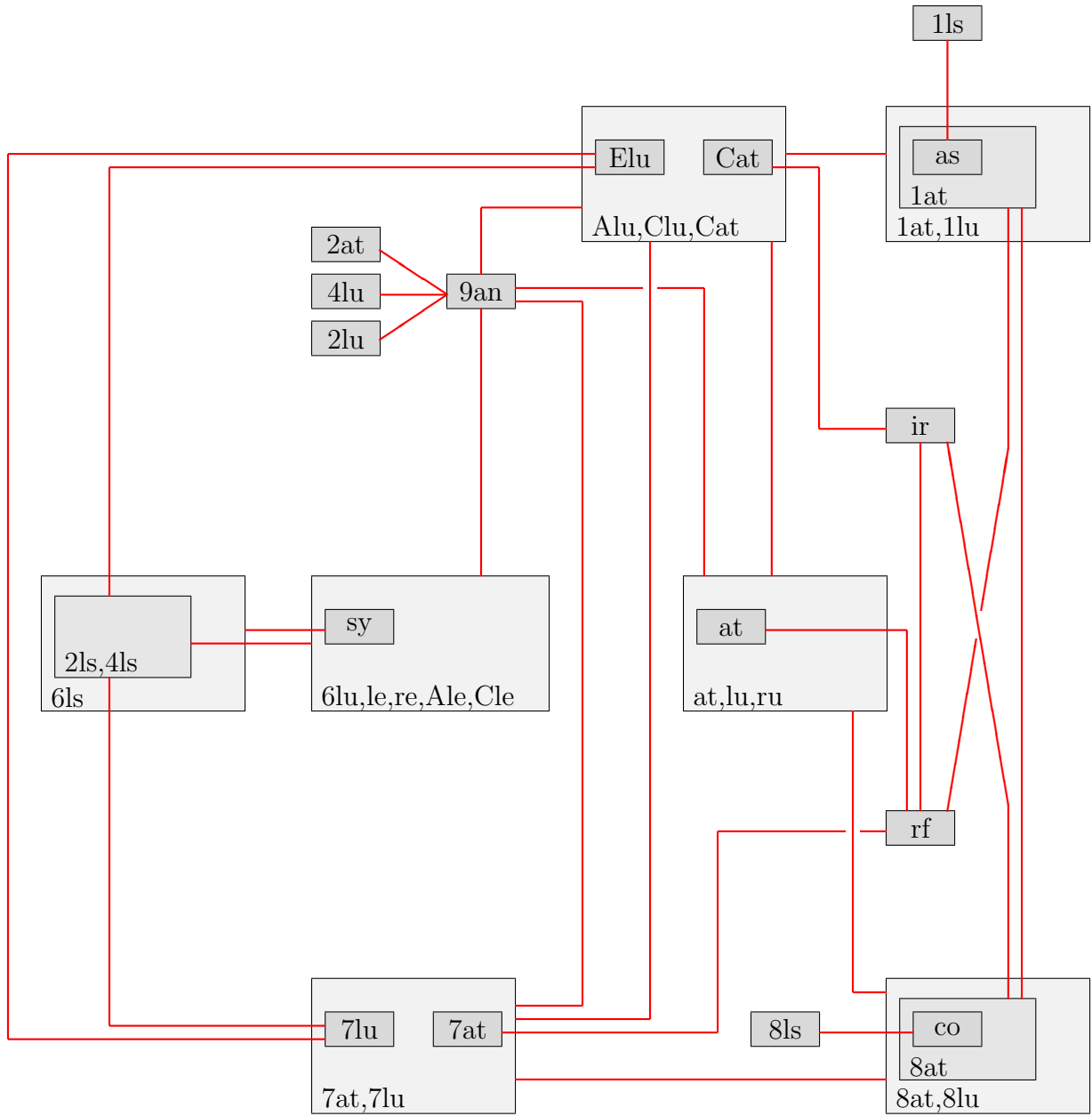

\clearpage

\subsection{Axiom set}
\LABEL{Axiom set}

In order to obtain a human-comprehensible result, we searched for a
small set of ``axioms'' from which all valid $3$-implications
could be derived by trivial inferences in the sense of
Sect.~\REF{Trivial inferences}.
Our approach is described in an abstract way by
Def.~\REF{Abstract derivation} and Lem.~\REF{Axiomatization}.

\DEFINITION{Abstract derivation}{%
\LABEL{Abstract derivation}%
Let $U$ (``universe'') be a set of ``abstract propositions''.
An abstract inference rule $r$ on $U$ is a partial
mapping from finite subsets of $U$ to $U$.
If $r(X)$ is defined for a finite $X \subseteq U$,
we say that $r(X)$ is directly derivable from $X$.

We say that $p$ is derivable from a set $Y$
if an abstract
derivation tree exists
with $p$ at its root and members of $Y$ as its leaves
such that each non-leaf node corresponds to a valid
direct derivation.
The set $Y$ needn't be finite, and not each of its members need to
appear as a leaf.
We will assume a fixed given set of inference rules in the following.

Moreover, let a fixed well-founded total order $<$ on $U$ be given.
We say a derivation tree is ordered if every of its nodes is larger
than all its (direct) descendants; we call its root orderly derivable
from the set of its leaves.
\qed
}

\begin{figure}
\begin{picture}(50,50)%
\put(33,46){\makebox(0,0)[b]{$k$}}%
\textcolor{red}{\put(33,45){\vector(-1,-1){8}}}%
\textcolor{red}{\put(33,45){\vector(1,-1){8}}}%
\textcolor{red}{\put(25,37){\vector(-1,-2){4}}}%
\textcolor{red}{\put(25,37){\vector(1,-2){4}}}%
\textcolor{red}{\put(41,37){\vector(-1,-2){4}}}%
\textcolor{red}{\put(41,37){\vector(1,-2){4}}}%
\put(37,26){\makebox(0,0){$k_i$}}%
\put(33,26){\oval(32,6)}%
\put(15,26){\makebox(0,0)[r]{$K \setminus \set{k}$}}%
\textcolor{green}{\put(37,23){\vector(-1,-1){8}}}%
\textcolor{green}{\put(37,23){\vector(1,-1){8}}}%
\textcolor{green}{\put(29,15){\vector(-1,-2){4}}}%
\textcolor{green}{\put(29,15){\vector(1,-2){4}}}%
\textcolor{green}{\put(45,15){\vector(-1,-2){4}}}%
\textcolor{green}{\put(45,15){\vector(1,-2){4}}}%
\put(37,4){\oval(32,6)}%
\put(19,4){\makebox(0,0)[r]{$A$}}%
\end{picture}
\hfill
\begin{picture}(55,55)%
\color{red}%
\put(0,10){\circle*{1}}%
\put(0.5,10.5){\vector(1,1){9}}%
\put(10,20){\circle{1}}%
\put(10.5,20.5){\vector(1,1){9}}%
\put(20,30){\circle{1}}%
\put(20.5,30.5){\vector(1,1){19}}%
\put(40,50){\circle{1}}%
\put(10,10){\circle*{1}}%
\put(10.5,10.5){\vector(1,1){9}}%
\put(20,20){\circle{1}}%
\put(20.5,20.5){\vector(1,1){9}}%
\put(30,30){\circle{1}}%
\put(10,0){\circle*{1}}%
\put(10.5,0.5){\vector(1,1){9}}%
\put(20,10){\circle{1}}%
\put(20.5,10.5){\vector(1,1){9}}%
\put(30,20){\circle{1}}%
\put(30.5,20.5){\vector(1,1){9}}%
\put(40,30){\circle{1}}%
\put(40.5,30.5){\vector(1,1){9}}%
\put(50,40){\circle{1}}%
\color{blue}%
\put(1,10){\circle*{1}}%
\put(11,0){\circle*{1}}%
\put(11,0.5){\vector(0,1){9}}%
\put(11,10){\circle{1}}%
\put(11,10.5){\vector(0,1){9}}%
\put(11,20){\circle{1}}%
\put(21,10){\circle*{1}}%
\put(21,10.5){\vector(0,1){9}}%
\put(21,20){\circle{1}}%
\put(21,20.5){\vector(0,1){9}}%
\put(21,30){\circle{1}}%
\put(31,20){\circle*{1}}%
\put(31,20.5){\vector(0,1){9}}%
\put(31,30){\circle{1}}%
\put(41,30){\circle*{1}}%
\put(41,30.5){\vector(0,1){19}}%
\put(41,50){\circle{1}}%
\put(51,40){\circle*{1}}%
\color{green}%
\put(10,1){\circle*{1}}%
\put(0,11){\circle*{1}}%
\put(0.5,11){\vector(1,0){9}}%
\put(10,11){\circle{1}}%
\put(10.5,11){\vector(1,0){9}}%
\put(20,11){\circle{1}}%
\put(10,21){\circle*{1}}%
\put(10.5,21){\vector(1,0){9}}%
\put(20,21){\circle{1}}%
\put(20.5,21){\vector(1,0){9}}%
\put(30,21){\circle{1}}%
\put(20,31){\circle*{1}}%
\put(20.5,31){\vector(1,0){9}}%
\put(30,31){\circle{1}}%
\put(30.5,31){\vector(1,0){9}}%
\put(40,31){\circle{1}}%
\put(50,41){\circle*{1}}%
\put(40,51){\circle*{1}}%
\end{picture}
\caption{Proof sketch in Lem.~\REF{Axiomatization} (lf) 
	\hspace*{1cm}
	Order-based axiom set computation (rg)}
\LABEL{Proof sketch in Lem Ax / Order-based axiom set computation}
\end{figure}

\LEMMA{Axiomatization}{%
\LABEL{Axiomatization}%
Let $V$ (``valid'') be a set of propositions;
let $K$ (``kernel'') be a subset of $V$
such that each member of $V$ can be derived from $K$.
Let $A$ (``axioms'')
be the set of propositions $k \in K$ such that $k$ is not orderly
derivable from $K \setminus \set{k}$.
Then:
\begin{enumerate}
\item $A$ is a subset of $K$;
\item Each member of $V$ is derivable already from $A$;
\item No member $a$ of $A$ is orderly derivable from $A \setminus \set{a}$.
\end{enumerate}
}
\PROOF{
We first show, by induction on $(<)$, that each $k \in K$
is orderly derivable from $A$:
If $k$ is not orderly derivable from $K \setminus \set{k}$,
the $k \in A$ by construction,
and we are done (trivial one-node derivation tree).
Else,
let $k$ be orderly derivable from 
$\set{ \:,1n{k_\i} } \subseteq K \setminus \set{k}$
(red tree in the left part of
Fig.~\REF{Proof sketch in Lem Ax / Order-based axiom set computation}).
Then $k_i < k$, and hence by I.H.\ $k_i$ is orderly derivable from $A$
(green tree).
Composing the tree to derive $k$ from $\set{ \:,1n{k_\i} }$
with the trees for all $k_i$ gives an ordered derivation tree for $k$
from $A$.

Now we show the claimed properties:
\begin{enumerate}
\item
	By construction of $A$.
\item
	If $v \in V$,
	then $v$ has a derivation tree with leaves in $K$ by assumption.
	Each of its leaves has an (even ordered) derivation tree from
	$A$.
	Composing the trees appropriately yields a derivation tree for
	$v$ from $A$ (note that it needn't be ordered).
\item
	If $a \in A$ was orderly derivable from
	$A \setminus \set{a}$,
	then it was also orderly derivable from the superset 
	$K \setminus \set{a}$,
	contradicting the construction of $A$.

	Note, however, that some $a \in A$ may still be unorderly
	derivable from larger members of $A$.
	%
\qed
\end{enumerate}
}

Figure~\REF{Proof sketch in Lem Ax / Order-based axiom set computation}
(rg)
shows a geometrical analogy of this approach of order-based axiom set
computation.
A member of the set $V$ is represented by $3$ adjacent circles in
different colors (red, green, blue).
Each color corresponds to an own proof ordering.
A colored arrow indicates an ordered derivation;
for simplicity we assume each inference rule to work on a singleton
set.
A full circle indicates an axiom w.r.t.\ the proof order of its color;
a hollow circle represents a non-axiom.
The proof order shown in red, green, and blue, needs $3$, $6$, and $6$
axioms, respectively.
Note that the set of red axioms can be further reduced by applying,
in turn, e.g.\ the green order to it.

In our implementation, an abstract proposition consists of a valid or
invalid $3$-implication, like
``C-Refl $\land$ 3-Trans $\ra$ 3-ASym''
and
``3-Trans $\land$ 3-Trans $\not\ra$ 3-Sym''
(\QDef{eccj} and \QDef{ilcb}, respectively).
We represented
such a proposition by its base-$27$ code,
using different permutations of character positions.
For example, the implications
``3-Trans $\land$ C-Refl $\ra$ 3-ASym'' and
``3-SemiOrd1 $\land$ 3-Refl $\ra$ C-ASym''
correspond to \QRem{ijrj} and \QRef{faok},
respectively;
so a derivation of the former from the latter is unorderly w.r.t.\ 
the identity permutation,
but orderly when comparing the reverse codes (\QRem{jrji} and
\QRem{koaf}).
Applying \texttt{xor} masks to the characters
before comparing, we obtain more possible orderings.
In order to obtain best-comprehensible axioms, we prepended
symbol encoding the ``human simplicity'' of an implication.\footnote{
	For example, the list
	``Refl $\land$ LfEucl $\ra$ 1-RgSerial'',
	``Refl $\ra$ 1-RgSerial'',
	``Refl $\ra$ RgSerial'',
	``Refl $\ra$ LfSerial''
	is ordered from worst to best simplicity.
}

We used an ad hoc evolutionary algorithm based on
Lem.~\REF{Axiomatization} to minimize the cardinality of the axiom set
as far as we could.
Starting from the set of all $34727$ proven implications,
we applied Lem.~\REF{Axiomatization} in parallel
for a couple of different
orderings.
One of them resulted in a set of $360$ axioms, the minimum number we
could achieve.
Repeating that procedure $5$ times,
we obtained axiom sets comprising $285$, $265$, $256$, $253$, and
$252$ axioms.
After that we manually tried to omit axioms one by one, and ended up
with a set of $124$ axioms.

Subsequently, we manually exchanged a few axioms by equivalent ones
for presentation reasons.
For example, we replaced ``Refl $\ra$ LfQuasiRefl'' (\QRem{dsqt})
by ``Refl $\ra$ QuasiRefl'' (\QRef{dsqv}), to get a nicer Hasse diagram in
Fig.~\REF{Implications about basic properties};
in the presence of the remaining $123$ axioms, both are equivalent.
Similarly, we replaced right by left properties where possible,
e.g.\ ``RgEucl $\ra$ QuasiTrans'' (\QRem{pej0}) by the equivalent
``LfEucl $\ra$ QuasiTrans'' (\QRef{owg0}).





The result is shown in Fig.~\REF{Minimal axiom set};
each valid $3$-implication can be inferred from these axioms
using the trivial inference rules from Sect.~\REF{Trivial inferences}
(Thm.~\REF{$3$-implication axioms}).

In Fig.~\REF{Minimal axiom set}, we marked implications where all
atoms share the same basic property.
We used ``$^*$'' when this can be obtained by applying just the
redundancies from Fig.~\REF{Redundant properties},
and ``$^+$'' when equivalences from
Fig.~\REF{Computed partitions of lifted properties}
are needed in addition;
for example, the top left implication ``IncTrans $\ra$ SemiOrd2''
can be rephrased as ``7-SemiOrd2 $\ra$ SemiOrd2''.

When a minimum cardinality is required for an implication, we noted it
in the same column; for example, the last axiom, \QDef{wfzb}, can be
falsified on a $2$-element set, but is valid for all relation domains
of $\geq 3$ elements.

Moreover, we marked axioms for which some kind of converse also holds.
An axiom $I \land J \ra K$ is marked ``$=$'' 
if also $K \ra I \land J$ holds; 
it is marked ``$\triangleright$'' if also $K \ra J$ holds,
and ``$>$'' if also $I \land K \ra J$ holds.
The axioms could be manually tuned such that the cases
$J \land K \ra I$ and $K \ra I$ don't occur.
Note that the reverse or semi-reverse versions are not part of the
axiomatization; for example, the bottom right axiom 
``ASym $\land$ Dense $\ra$ 2-Dense'' (\QRem{cllc}) is
marked "$>$" since ASym and 2-Dense also implies Dense (\QRem{cijd}),
but the latter can be proved from axioms
\QRem{adsd},
\QRem{cixg},
\QRem{cpbj},
\QRef{dcaa},
\QRef{owht},
\QRem{pgrm}, and
\QRef{xzgd}.

Figure~\REF{Minimal axiom set} shows the $20$
``$2$-implication laws'' first,
then the single ``$2$-inconsistency law'' \QRef{hpr0},
then the $103$ proper ``$3$-implication laws''.
We additionally grouped implications by the basic/lifted distinctions
of their atoms.

As for invalid implications, we minimized the set of implications
obtained from finite (Sect.~\REF{Finite counter-examples})
and infinite (Sect.~\REF{Infinite counter-examples})
counter-examples
by manual ad-hoc removals; we didn't apply Lem.~\REF{Axiomatization}
for them, and we didn't try to minimize the set of computer-generated
counter-examples.
$17$ implications were falsified by examples on an infinite domain,
$133$ on a $5$-element domain (computer-generated),
and $16$ on another finite domain.

\begin{figure}
\renewcommand{\arraystretch}{1.0}
\newcommand{\1}{$\scriptstyle *$}	
\newcommand{\2}{$\scriptscriptstyle +$}	
\newcommand{\3}[1]{$\scriptstyle #1$}	
\newcommand{\7}{$\scriptstyle >$}	
\newcommand{\8}{$\scriptstyle\triangleright$}	
\newcommand{\9}{$\scriptstyle =$}	
\newcommand{\0}{@{}|@{\hspace*{0.1mm}}l@{}c@{}r@{\hspace*{0.1mm}}|rr@{ $\ra$ }r|@{}}
\begin{tabular}[t]{\0}
\hline
\QRef{biuv} & \2  &    && it & s2	\\	
\QRef{dsqv} &     &    && rf & qr	\\	
\QRef{ild0} & \1  &    && tr & qt	\\
\QRef{owht} &     &    && le & lq	\\
\QRef{owg0} &     &    && le & qt	\\
\QRef{xzgd} &     &    && lq & de	\\
\QRef{yqnt} & \1  &    && qr & lq	\\
\cline{1-3}
\QRef{dspx} &     &    && rf & 1-ls	\\	
\QRef{fbao} &     &    && s1 & 1-de	\\
\QRef{fbar} &     &    && s1 & 8-de	\\
\QRef{gknw} &     &    && s2 & 8-lq	\\
\QRem{ilce} & \1  &    && tr & 1-tr	\\
\QRef{jlmq} & \1  &    && de & 7-de	\\
\QRef{owgk} &     &    && le & 6-at	\\
\QRef{xzhs} & \1  &    && lq & 1-lq	\\
\cline{1-3}
\QRef{fkdd} &     &    && 7-s1 & de	\\
\QRef{lmco} &     &    && 1-at & ir	\\
\QRem{pwog} & \2  &    && C-le & it	\\	
\QRem{vpm0} &     &    && 6-lu & qt	\\
\cline{1-3}
\QRef{fkdo} &     &    && 7-s1 & 1-de	\\
\hline
\QRef{hpr0} & \34 &    && 8-lu & $\lnot$at	\\
\hline
\QRef{0udt} &     & \7 & sy & tr & le	\\
\QRef{d0gh} & \2  & \9 & an & sy & cr	\\	
\QRef{dcaa} &     & \8 & an & qt & tr	\\
\QRef{faok} &     & \9 & s1 & rf & co	\\
\QRef{hlmf} &     & \9 & at & de & em	\\
\QRef{ikuv} &     &    & tr & at & s2	\\
\QRef{oskj} &     & \7 & le & an & lu	\\
\QRef{phfb} &     & \7 & re & lq & sy	\\
\QRef{ujry} &     & \8 & lu & ir & at	\\
\QRef{ukaa} &     &    & lu & s1 & tr	\\
\QRef{ulmt} &     & \8 & lu & de & le	\\
\QRef{xtlr} &     &    & lq & s2 & s1	\\
\QRef{xxcn} &     & \9 & lq & rs & rf	\\
\QRef{xzkv} & \1  & \9 & lq & rq & qr	\\
\cline{1-3}
\QRef{0t0p} & \2  & \7 & sy & s1 & 1-s1	\\	
\QRef{0ucf} & \2  & \7 & sy & tr & 7-tr	\\	
\QRef{0wyx} &     & \8 & sy & ls & 1-ls	\\
\QRef{0yen} &     & \8 & sy & lu & 7-lu	\\
\QRem{bldc} & \1  &    & it & tr & D-tr	\\
\QRem{cllc} &     & \7 & as & de & 2-de	\\
\hline
\end{tabular}
\hfill
\begin{tabular}[t]{\0}
\hline
\QRef{cnyy} &     & \8 & as & ls & 2-ls	\\
\QRem{djid} & \1  & \9 & sc & an & 9-an	\\
\QRem{ik0q} &     &    & tr & s1 & 2-s1	\\
\QRem{ikvl} &     & \8 & tr & at & 7-at	\\
\QRem{ot0q} &     &    & le & s1 & 2-s1	\\
\QRem{qtaq} &     &    & ls & s1 & 7-de	\\
\QRef{qudq} &     &    & ls & tr & 7-de	\\
\QRef{uivk} & \34 &    & lu & it & 6-at	\\
\QRef{uklu} & \35 &    & lu & s2 & 2-s2	\\
\QRem{uodr} & \33 &    & lu & rs & 8-de	\\
\QRem{uoed} & \32 &    & lu & rs & 8-ls	\\
\QRem{xt0s} &     &    & lq & s1 & 7-s1	\\
\cline{1-3}
\QRef{0vcy} &     & \8 & sy & 1-at & at	\\
\QRef{0vyd} &     & \8 & sy & 7-de & de	\\
\QRef{0xjz} &     & \8 & sy & 7-ls & ls	\\
\QRef{0xzj} &     & \8 & sy & 1-lu & lu	\\
\QRef{0zet} & \2  & \8 & sy & 1-lq & lq	\\	
\QRem{cpbj} &     & \7 & as & 2-lu & lu	\\
\QRem{cixg} &     & \7 & as & 9-tr & it	\\
\QRem{ckry} &     & \8 & as & 2-at & at	\\
\QRem{ikcr} &     &    & tr & 7-s1 & s1	\\
\QRem{otcr} &     & \7 & le & 7-s1 & s1	\\
\QRem{pgrm} & \33 & \8 & re & 8-lu & sc	\\
\QRem{ulfl} & \33 &    & lu & C-tr & an	\\
\QRem{umxl} &     &    & lu & 7-de & an	\\
\QRem{unrb} &     &    & lu & 1-ls & sy	\\
\QRem{unuj} &     &    & lu & 2-ls & as	\\
\QRem{xtg0} &     &    & lq & C-s1 & qt	\\
\QRem{xxjz} &     & \8 & lq & 7-ls & ls	\\
\QRef{xzog} &     & \8 & lq & 8-lq & it	\\
\QRem{ykcr} &     & \7 & qr & 7-s1 & s1	\\
\cline{1-3}
\QRem{bioh} & \1  &    & it & 1-tr & 9-tr	\\
\QRem{cmyp} &     & \7 & as & 7-de & 6-de	\\
\QRem{cokb} &     & \8 & as & 7-ls & 6-ls	\\
\QRem{cpnn} &     & \8 & as & 6-lu & 7-lu	\\
\QRem{ioas} &     &    & tr & 4-ls & E-de	\\
\QRem{jljo} &     &    & de & D-tr & 1-de	\\
\QRem{otfp} &     & \7 & le & C-s1 & 1-s1	\\
\QRem{owae} &     &    & le & 8-de & C-de	\\
\QRem{uiyr} & \33 &    & lu & 9-tr & 8-de	\\
\QRem{ukfp} &     &    & lu & C-s1 & 1-s1	\\
\hline
\end{tabular}
\hfill
\begin{tabular}[t]{\0}
\hline
\QRem{unae} &     &    & lu & 8-de & C-de	\\
\QRem{uoge} &     &    & lu & 6-ls & C-de	\\
\QRem{vhqx} &     & \8 & ru & 8-lq & A-lq	\\
\QRem{ykfq} &     &    & qr & C-s1 & 2-s1	\\
\cline{1-3}
\QRem{btlg} & \2  & \8 & 9-tr & s2 & it	\\	
\QRem{eqgt} &     & \8 & 1-s1 & lq & le	\\
\QRem{npkj} &     & \7 & 6-de & ru & lu	\\
\QRem{ntuj} &     &    & 7-de & at & as	\\
\QRem{pudt} &     & \7 & C-le & tr & le	\\
\QRem{qbad} &     &    & 1-ls & s1 & de	\\
\QRem{vtlj} & \33 & \7 & 7-lu & s2 & as	\\
\QRem{vxcb} &     &    & 7-lu & rs & sy	\\
\cline{1-3}
\QRef{gtlu} & \1  &    & C-s2 & s2 & 2-s2	\\
\QRem{ituh} &     &    & C-tr & at & 9-tr	\\
\QRem{itui} &     &    & C-tr & at & E-tr	\\
\QRem{j0vb} & \1  & \8 & D-tr & it & C-tr	\\
\QRem{ltvl} & \1  & \9 & 6-at & at & 7-at	\\
\QRem{pxdx} &     & \8 & C-le & rs & 1-ls	\\
\QRef{uydr} & \35 &    & 4-lu & lu & 8-de	\\
\QRem{wcch} &     &    & 8-lu & tr & 9-tr	\\
\QRem{wxdx} &     & \8 & C-lu & rs & 1-ls	\\
\cline{1-3}
\QRem{adsd} & \1  &    & 2-de & 1-de & de	\\
\QRem{etfr} & \1  & \7 & 2-s1 & C-s1 & s1	\\
\QRem{jbcr} &     &    & D-tr & 7-s1 & s1	\\
\QRem{pkcr} &     & \7 & A-le & 7-s1 & s1	\\
\QRem{pxjz} &     & \8 & C-le & 7-ls & ls	\\
\QRem{qavd} &     &    & 1-ls & 1-s1 & de	\\
\QRem{qbgd} &     &    & 1-ls & C-s1 & de	\\
\QRem{rwpd} &     &    & 6-ls & C-le & de	\\
\QRef{vxfj} &     &    & 7-lu & 6-ls & as	\\
\QRem{wgba} & \33 &    & 8-lu & 2-lu & rs	\\
\QRef{xlfr} &     &    & 1-lq & C-tr & s1	\\
\cline{1-3}
\QRem{ewvq} &     &    & 2-s1 & 2-ls & 7-de	\\
\QRem{fjvk} &     &    & 7-s1 & 1-s1 & 6-at	\\
\QRef{g0rw} & \2  &    & 2-s2 & 7-tr & C-s2	\\	
\QRem{nnas} & \1  &    & 6-de & 8-de & E-de	\\
\QRef{voye} & \35 &    & 6-lu & 1-lu & C-de	\\
\QRef{voyr} & \35 &    & 6-lu & 1-lu & 8-de	\\
\QRem{vozd} & \33 &    & 6-lu & 1-lu & 8-ls	\\
\QRem{wfxh} &     &    & 8-lu & 1-lu & 9-tr	\\
\QRef{wfyp} & \35 &    & 8-lu & 1-lu & 6-de	\\
\QRef{wfzb} & \33 &    & 8-lu & 1-lu & 6-ls	\\
\hline
\end{tabular}
\begin{center}
\caption{Minimal axiom set}
\LABEL{Minimal axiom set}
\end{center}
\end{figure}

%

\begin{figure}
\begin{center}
\includegraphics[width=\linewidth]{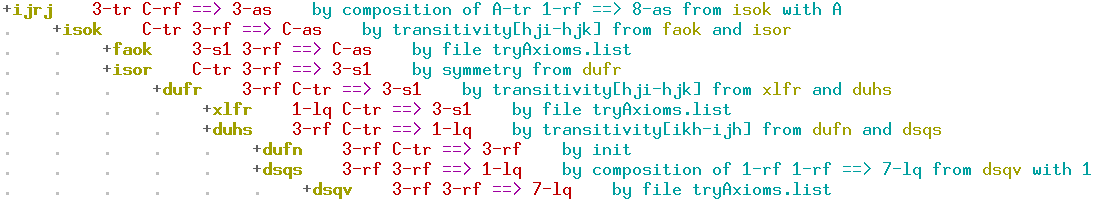}
\caption{Proof tree of ``Trans $\land$ Irrefl $\ra$ ASym''}
\LABEL{Proof tree of ``tr land ir ra as''}
\end{center}
\end{figure}

\begin{figure}
\begin{center}
\includegraphics[width=\linewidth]{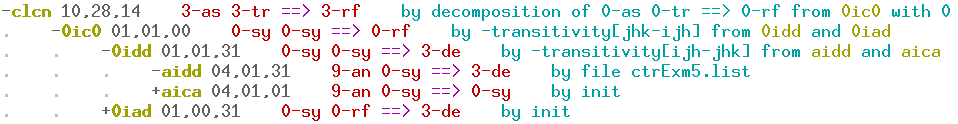}
\caption{Disproof tree of ``ASym $\land$ Trans $\not\ra$ Refl''}
\LABEL{Disproof tree of ``as land tr not ra rf''}
\end{center}
\end{figure}

\subsection{Some proofs of axioms}
\LABEL{Some proofs of axioms}

All axioms from Fig.~\REF{Minimal axiom set}
could be proven by EProver, except
\QRef{g0rw},
\QRef{uivk},
\QRef{uklu},
\QRef{uydr},
\QRef{voye},
\QRef{voyr}, and
\QRef{wfyp}.
In this section, we give manual proofs for these and a few other
axioms.
Moreover, we include references to all proofs of axioms that were
already given in \cite{Burghardt.2018c}.

\LEMMA{2-Implications group 0}{%
\LABEL{2-Implications group 0}%
\begin{enumerate}
\item\LABEL{2-Implications group 0 biuv}
	(\QDef{biuv})	IncTrans $\ra$ SemiOrd2
\item\LABEL{2-Implications group 0 dsqv}%
	(\QDef{dsqv})	Refl $\ra$ 7-LfQuasiRefl
\item\LABEL{2-Implications group 0 ild0}
	(\QDef{ild0})	Trans $\ra$ QuasiTrans
\item\LABEL{2-Implications group 0 owg0}
	(\QDef{owg0})	LfEucl $\ra$  QuasiTrans
\item\LABEL{2-Implications group 0 owht}%
	(\QDef{owht})	LfEucl $\ra$ LfQuasiRefl
\item\LABEL{2-Implications group 0 xzgd}
	(\QDef{xzgd})	LfQuasiRefl $\ra$ Dense
\item\LABEL{2-Implications group 0 yqnt}
	(\QDef{yqnt})	QuasiRefl $\ra$ LfQuasiRefl
\end{enumerate}
}
\PROOF{
\begin{enumerate}
\item
	See \OEF{Lem.34}{31}{}.
\item
	Since Refl implies QuasiRefl by \OEF{Lem.9}{23}{};
	the latter is equivalent to 7-LfQuasiRefl
	by Lem.~\REFF{optn redundant}{13}.
\item
	See \OEF{Lem.18}{27}{}.
\item
	See \OEF{Lem.40}{34}{}.
\item
	See \OEF{Lem.46}{37}{}.
\item
	See \OEF{Lem.48.3}{38}{}.
\item
	By Def.~\REFF{def}{QuasiRefl}.
\qed
\end{enumerate}
}

\LEMMA{2-Implications group 1}{%
\LABEL{2-Implications group 1}%
\begin{enumerate}
\item\LABEL{2-Implications group 1 dspx}%
	(\QDef{dspx})	Refl $\ra$ 1-LfSerial
\item\LABEL{2-Implications group 1 fbao}%
	(\QDef{fbao})	SemiOrd1 $\ra$ 1-Dense
\item\LABEL{2-Implications group 1 fbar}%
	(\QDef{fbar})	SemiOrd1 $\ra$ 8-Dense
\item\LABEL{2-Implications group 1 gknw}%
	(\QDef{gknw})	SemiOrd2 $\ra$ 8-LfQuasiRefl
\item\LABEL{2-Implications group 1 jlmq}%
	(\QDef{jlmq})	Dense $\ra$ 7-Dense
\item\LABEL{2-Implications group 1 owgk}%
	(\QDef{owgk})	LfEucl $\ra$ 6-AntiTrans
\item\LABEL{2-Implications group 1 xzhs}%
	(\QDef{xzhs})	LfQuasiRefl $\ra$ 1-LfQuasiRefl
\end{enumerate}
}
\PROOF{
\begin{enumerate}
\item
	Refl implies LfSerial and RgSerial 
	by \OEF{Lem.54}{40}{ser 1};
	the latter conjunction is equivalent to 1-LfSerial
	by Lem.~\REFF{optn bool}{2}
	and~\REFF{optn redundant}{1}.
\item
	Let $x \ro{R}{1} z$,
	that is, $x \rr{R} z \land z \rr{R} x$.
	If $x \rr{R} x$, then $x \ro{R}{1} x$, 
	and we can choose $x$ as
	intermediate element.
	If $\lnot x \rr{R} x$,
	then $z \rr{R} z$ by semi-order property~1
	applied to $z \rr{R} x$, $x,x$ incomparable, 
	$x \rr{R} z$;
	so we can choose $z$ as intermediate element.
\item
	Let $x \ro{R}{8} z$,
	that is, $\lnot x \rr{R} z \land \lnot z \rr{R} x$,
	we distinguish three cases:
	\begin{itemize}
	\item If $\lnot x \rr{R} x$, then $x \ro{R}{8} x$,
		so $x$ can be used as intermediate element.
	\item if $\lnot z \rr{R} z$, 
		then $z$ can be used, in a similar way.
	\item If $x \rr{R} x \land z \rr{R} z$, 
		then applying semi-order property~1
		yields $x \rr{R} z$, contradicting $x \ro{R}{8} z$.
	\end{itemize}
\item
	Let $x \ro{R}{8} y$,
	assume for contradiction $\lnot x \ro{R}{8} x$.
	The latter means $x \rr{R} x$ by definition.
	Hence $x,y$ must be comparable w.r.t.\ $R$,
	by semi-order-property~2 applied to 
	$x \rr{R} x$ and $x \rr{R} x$.
	But $x \ro{R}{8} y$ means that $x,y$ are incomparable w.r.t.\ $R$.
\item
	Let $x \ro{R}{7} z$,
	then $x \rr{R} z$ or $z \rr{R} x$ by definition.
	In the first case, we have
	$x \rr{R} y \land y \rr{R} z$ for some $y$,
	by density of $R$,
	hence $x \ro{R}{7} y \land y \ro{R}{7} z$
	by Lem.~\REFF{optn bool}{4}.
	The second case is similar.
\item
	Let $x \ro{R}{6} y$ and $y \ro{R}{6} z$, 
	we show $\lnot x \ro{R}{6} z$.
	By definition, we have to distinguish the following cases:
	\begin{itemize}
	\item $x \rr{R} y \land \lnot y \rr{R} x$ 
		and $y \rr{R} z \land \lnot z \rr{R} y$:
		\\
		Then $y \rr{R} y$ by left Euclideanness,
		hence $y \rr{R} x$ by the same property,
		contradicting the case assumption.
	\item $x \rr{R} y \land \lnot y \rr{R} x$ 
		and $\lnot y \rr{R} z \land z \rr{R} y$:
		\\
		Then $x \rr{R} z$ and $z \rr{R} x$ 
		by left Euclideanness,
		that is, $\lnot x \ro{R}{6} z$.
	\item $\lnot x \rr{R} y \land y \rr{R} x$ 
		and $y \rr{R} z \land \lnot z \rr{R} y$:
		\\
		Then $x \rr{R} z$ would imply the contradiction 
		$x \rr{R} y$,
		and $z \rr{R} x$ would imply the contradiction 
		$z \rr{R} y$.
		But $\lnot x \rr{R} z \land \lnot z \rr{R} x$ 
		implies $\lnot x \ro{R}{6} z$.
	\item $\lnot x \rr{R} y \land y \rr{R} x$ 
		and $\lnot y \rr{R} z \land z \rr{R} y$:
		\\
		Then $z \rr{R} x$ would imply the contradiction 
		$y \rr{R} z$,
		and $x \rr{R} z$ would imply 
		$z \rr{R} x$ (since $z \rr{R} z$) which just
		has been show to contradict.
		Again,  $\lnot x \rr{R} z \land \lnot z \rr{R} x$
		implies $\lnot x \ro{R}{6} z$.
	\end{itemize}
\item
	Let $x \ro{R}{1} y$,
	then $x \rr{R} y \land y \rr{R} x$ by definition,
	hence $x \rr{R} x$ by left quasi-reflexivity,
	hence $x \ro{R}{1} x$ by definition.
\qed
\end{enumerate}
}

\LEMMA{2-Implications group 2}{%
\LABEL{2-Implications group 2}%
\begin{enumerate}
\item\LABEL{2-Implications group 2 fkdd}%
	(\QDef{fkdd})	7-SemiOrd1 $\ra$ Dense
\item\LABEL{2-Implications group 2 lmco}%
	(\QDef{lmco})	1-AntiTrans $\ra$ C-Refl
\end{enumerate}
}
\PROOF{
\begin{enumerate}
\item
	Let $x \rr{R} z$,
	we distinguish two cases:
	\begin{itemize}
	\item If $x \rr{R} x$, we can use $x$ as intermediate element.
	\item If $\lnot x \rr{R} x$,
		then applying semi-order property~1
		to $z \ro{R}{7} x$,
		$x,x$ incomparable w.r.t.\ $R^7$,
		$x \ro{R}{7} z$
		yields $z \ro{R}{7} z$,
		that is, $z \rr{R} z$,
		and we can use $z$ as intermediate element.
	\end{itemize}
\item
	Assume for contradiction $\lnot x \ro{R}{C} x$ for some $x$.
	Then $x \rr{R} x$,
	hence $x \ro{R}{1} x \land x \ro{R}{1} x$ by definition,
	hence $\lnot x \ro{R}{1} x$ by anti-transitivity.
\qed
\end{enumerate}
}

\LEMMA{2-Implications group 3}{%
\LABEL{2-Implications group 3}%
\begin{enumerate}
\item\LABEL{2-Implications group 3 fkdo}%
	(\QDef{fkdo})	7-SemiOrd1 $\ra$ 1-Dense
\end{enumerate}
}
\PROOF{
\begin{enumerate}
\item
	The proof is similar to that 
	of~\REFF{2-Implications group 2}{fkdd}:
	Let $x \ro{R}{1} z$,
	that is, $x \rr{R} z \land z \rr{R} x$,
	then $x \ro{R}{7} z$ and $z \ro{R}{7} x$.
	We distinguish two cases:
	\begin{itemize}
	\item If $x \rr{R} x$,
		then $x \ro{R}{1} x$,
		and we can use $x$ as intermediate element.
	\item If $\lnot x \rr{R} x$,
		then applying semi-order property~1
		to $z \ro{R}{7} x$,
		$x,x$ incomparable w.r.t.\ $R^7$,
		$x \ro{R}{7} z$
		yields $z \ro{R}{7} z$,
		that is, $z \rr{R} z$,
		that is $z \ro{R}{1} z$,
		and we can use $z$ as intermediate element.
\qed
	\end{itemize}
\end{enumerate}
}

\LEMMA{2-Contradictions group 0}{%
\LABEL{2-Contradictions group 0}%
\begin{enumerate}
\item\LABEL{2-Contradictions group 0 hpr0}
	(\QDef{hpr0})
	No relation on a set of $\geq 4$ elements can be both
	AntiTrans and 8-LfUnique.
\end{enumerate}
}
\PROOF{
\begin{enumerate}
\item
	Assume for contradiction that $R$ is AntiTrans and $\ro{R}{8}$ is
	LfUnique.
	By \OEF{Lem.22}{28}{},
	AntiTrans implies Irrefl.
	Let $w,x,y,z$ be $4$ distinct elements.
	By 8-LfUnique, since each element is incomparable to itself,
	it must be comparable to the other three.
	Hence, in the directed graph corresponding to $R$,
	the vertice of $w$ must have two incoming edges or two 
	outgoing ones,
	leading w.l.o.g.\ to $x$ and $y$ such that, moreover, w.l.o.g.\
	$x \rr{R} y$.
	In the incoming case, we have 
	$x \rr{R} y$ and
	$y \rr{R} w$, but
	$x \rr{R} w$.
	In the outgoing case, we have
	$w \rr{R} x$ and
	$x \rr{R} y$, but
	$w \rr{R} y$.
	Both contradict AntiTrans.
\qed
\end{enumerate}
}

\LEMMA{3-Implications group 0}{%
\LABEL{3-Implications group 0}%
\begin{enumerate}
\item\LABEL{3-Implications group 0 0udt}%
	(\QDef{0udt})	Sym $\land$ Trans $\ra$ LfEucl
\item\LABEL{3-Implications group 0 d0gh}
	(\QDef{d0gh})
	AntiSym $\land$ Sym $\ra$ CoRefl
\item\LABEL{3-Implications group 0 dcaa}
	(\QDef{dcaa})
	AntiSym $\land$ QuasiTrans $\ra$ Trans
\item\LABEL{3-Implications group 0 faok}
	(\QDef{faok})
	SemiOrd1 $\land$ Refl $\ra$ Connex
\item\LABEL{3-Implications group 0 hlmf}
	(\QDef{hlmf})
	AntiTrans $\land$ Dense $\ra$ Empty
\item\LABEL{3-Implications group 0 ikuv}
	(\QDef{ikuv})
	Trans $\land$ AntiTrans $\ra$ SemiOrd2
\item\LABEL{3-Implications group 0 oskj}
	(\QDef{oskj})
	LfEucl $\land$ AntiSym $\ra$ LfUnique
\item\LABEL{3-Implications group 0 phfb}
	(\QDef{phfb})
	RgEucl $\land$ LfQuasiRefl $\ra$ Sym
\item\LABEL{3-Implications group 0 ujry}
	(\QDef{ujry})
	LfUnique $\land$ Irrefl $\ra$ AntiTrans
\item\LABEL{3-Implications group 0 ukaa}
	(\QDef{ukaa})
	LfUnique $\land$ SemiOrd1 $\ra$ Trans
\item\LABEL{3-Implications group 0 ulmt}
	(\QDef{ulmt})
	LfUnique $\land$ Dense $\ra$ LfEucl
\item\LABEL{3-Implications group 0 xtlr}
	(\QDef{xtlr})
	LfQuasiRefl $\land$ SemiOrd2 $\ra$ SemiOrd1
\item\LABEL{3-Implications group 0 xxcn}
	(\QDef{xxcn})
	LfQuasiRefl $\land$ RgSerial $\ra$ Refl
\item\LABEL{3-Implications group 0 xzkv}
	(\QDef{xzkv})
	LfQuasiRefl $\land$ RgQuasiRefl $\ra$ QuasiRefl
\end{enumerate}
}
\PROOF{
\begin{enumerate}
\item
	See \OEF{Lem.36}{33}{}.
\item
	See \OEF{Lem.7.7}{22}{}.
\item
	See \OEF{Lem.19}{27}{}.
\item
	See \OEF{Lem.66}{45}{}.
\item
	See \OEF{Lem.49}{39}{}.
\item
	See \OEF{Lem.24}{28}{}.
\item
	See \OEF{Lem.45}{37}{}.
\item
	See \OEF{Lem.37}{33}{}.
\item
	Assume for contradiction $x \rr{R} y$ and $y \rr{R} z$,
	but $x \rr{R} z$.
	then $x=y$ by LfUnique,
	contradicting Irrefl.
\item
	See \OEF{Lem.62.2}{43}{}.
\item
	See \OEF{Lem.47.1}{37}{}.
\item
	See \OEF{Lem.73}{47}{}.
\item
	See \OEF{Lem.55.2}{40}{}.
\item
	By Def.~\REFF{def}{QuasiRefl}.
\qed
\end{enumerate}
}

In Lem.~\REFF{3-Implications group 0}{ikuv},
note that a relation $R$ is both Trans and AntiTrans iff
$R$ is vacuously transitive, i.e.\ 
$\lnot \exists x,y,z. \;\; x \rr{R} y \land y \rr{R} z$.
Such a relation also is
ASym,
1-, 2-, 6-, 7-AntiTrans,
1-, 8-, C-, E-Dense,
1-, 8-, A-, C-, E-LfQuasiRefl,
8-, A-, C-, E-LfSerial,
1-LfUnique,
1-, C-SemiOrd1,
2-, 3-, C-SemiOrd2,
and
1-, 2-, E-Trans.

As a side remark,
a relation is both LfQuasiRefl and C-LfQuasiRefl
iff
it is ``right-constant'', i.e.\ it satisfies
$\forall x, y_1, y_2 \in X. x \rr{R} y_1 \lra x \rr{R} y_2$.
Such a relation also satisfies
2- and 6-AntiTrans,
1-, 3-, 7-, 8-, C- E-Dense,
3-, C-LfEucl,
1-, 8-LfQuasiRefl,
1-, 2-, 3-, 7-, C-SemiOrd1,
2-, 3-, C-SemiOrd2,
and
1-, 2-, 3-, 8-, 9-, C-, D-Trans.
Dually, $R$ is 5-LfQuasiRefl and A-LfQuasiRefl 
iff 
it is left-constant (defined similarly).
A relation is LfQuasiRefl and A-LfQuasiRefl
iff
it is totally constant, i.e.\ it satisfies
$\forall x_1, x_2, y_1, y_2 \in X. \;\;
x_1 \rr{R} y_1 \lra x_2 \rr{R} y_2$.

It seems promising to investigate more classes of relations
that are characterizable by conjunctions of lifted properties.

\LEMMA{3-Implications group 1}{%
\LABEL{3-Implications group 1}%
\begin{enumerate}
\item\LABEL{3-Implications group 1 0t0p}%
	(\QDef{0t0p})
	Sym $\land$ SemiOrd1 $\ra$ 1-SemiOrd1
\item\LABEL{3-Implications group 1 0ucf}%
	(\QDef{0ucf})
	Trans $\land$ Sym $\ra$ 7-Trans
\item\LABEL{3-Implications group 1 0wyx}%
	(\QDef{0wyx})
	Sym $\land$ LfSerial $\ra$ 1-LfSerial
\item\LABEL{3-Implications group 1 0yen}%
	(\QDef{0yen})
	Sym $\land$ LfUnique $\ra$ 7-LfUnique
\item\LABEL{3-Implications group 1 cnyy}
	(\QDef{cnyy})
	ASym $\land$ LfSerial $\ra$ 2-LfSerial
\item\LABEL{3-Implications group 1 qudq}%
	(\QDef{qudq})
	LfSerial $\land$ Trans $\ra$ 7-Dense
\item\LABEL{3-Implications group 1 uivk}%
	(\QDef{uivk})
	On a set of $\geq 4$ elements,
	3-LfUnique $\land$ 8-Trans $\ra$ 6-AntiTrans
\item\LABEL{3-Implications group 1 uklu}%
	(\QDef{uklu})
	On a set of $\geq 5$ elements,
	3-LfUnique $\land$ 3-SemiOrd2 $\ra$ 2-SemiOrd2
\end{enumerate}
}
\PROOF{
\begin{enumerate}
\item\LABEL{proof 0t0p}
	If $R$ is symmetric, then it coincides with
	its symmetric kernel $\ro{R}{1}$;
	hence, if the former is SemiOrd1, so is the latter.
\item
	Similar to~\REF{proof 0t0p}.
\item
	Similar to~\REF{proof 0t0p}.
\item
	Similar to~\REF{proof 0t0p}.
\item
	Given $y$, we find some $x$ with $x \rr{R} y$ by LfSerial,
	this implies $\lnot y \rr{R} x$ by ASym,
	hence $x \ro{R}{2} y$.
\item
	Let $x \ro{R}{7} y$ hold; w.l.o.g.\ let $x \rr{R} y$ hold.
	By left seriality of $R$, we find some $x'$ such that $x' R x$;
	by transitivity, this implies $x' R y$.
	Therefore, $x \ro{R}{7} x'$ and $x' \ro{R}{7} y$.
\item
	Assume for contradiction
	$x \ro{R}{6} y$ and $y \ro{R}{6} z$, but $x \ro{R}{6} z$.
	That is, we have exactly one of $xRy$ and $yRx$,
	and similar for $y,z$ and for $x,z$.
	Of the $8$ cases, only $2$ satisfy LfUnique,
	viz.\ those corresponding to directed cycles.
	By 8-Trans,
	any fourth element $w$ must be comparable to
	(w.l.o.g.) $x$ and $y$;
	this is impossible due to LfUnique.
\item
	Assume for contradiction $x \ro{R}{2} y$ and $y \ro{R}{2} z$,
	let $w \neq w'$ be both distinct from $x,y,z$.
	From 3-SemiOrd2,
	we obtain that $w$ is related (w.r.t.\ $R$) to one of $x,y,z$.
	By 3-LfUnique, 
	neither $w \rr{R} z$ nor $w \rr{R} y$ can hold,
	so the only case that could violate 2-SemiOrd2 is
	$x \rr{R} w \land w \rr{R} x$.
	By the same argument, we obtain $x R w' \land w' Rx$.
	However, $w \rr{R} x \land w' R x$ violates 3-LfUnique.
\qed
\end{enumerate}
}

\LEMMA{3-Implications group 2}{%
\LABEL{3-Implications group 2}%
\begin{enumerate}
\item\LABEL{3-Implications group 2 0vcy}%
	(\QDef{0vcy})
	Sym $\land$ 1-AntiTrans $\ra$ AntiTrans
\item\LABEL{3-Implications group 2 0vyd}%
	(\QDef{0vyd})
	Sym $\land$ 7-Dense $\ra$ Dense
\item\LABEL{3-Implications group 2 0xjz}%
	(\QDef{0xjz})
	Sym $\land$ 7-LfSerial $\ra$ LfSerial
\item\LABEL{3-Implications group 2 0xzj}%
	(\QDef{0xzj})
	Sym $\land$ 1-LfUnique $\ra$ LfUnique
\item\LABEL{3-Implications group 2 0zet}%
	(\QDef{0zet})
	Sym $\land$ 1-LfQuasiRefl $\ra$ LfQuasiRefl
\item\LABEL{3-Implications group 2 xzog}%
	(\QDef{xzog})
	LfQuasiRefl $\land$ 8-LfQuasiRefl $\ra$ 8-Trans
\end{enumerate}
}
\PROOF{
\begin{enumerate}
\item\LABEL{proof 0vcy}
	If $R$ is symmetric, then its symmetric kernel $\ro{R}{1}$
	coincides with $R$ itself;
	hence, if the former is anti-transitive, so is the latter.
\item
	Similar to~\REF{proof 0vcy}, using the symmetric closure
	instead of the kernel.
\item
	Similar to~\REF{proof 0vcy}.
\item
	Similar to~\REF{proof 0vcy}.
\item
	Similar to~\REF{proof 0vcy}.
\item
	Let $x,y$ and $y,z$ be incomparable w.r.t.\ $R$.
	Then $x,x$ and $y,y$ are incomparable w.r.t.\ $R$ by 8-LfQuasiRefl.
	Assume for contradiction $x \rr{R} z \lor z \rr{R} x$.
	In the first case, 
	we have the contradiction $x \rr{R} x$ by 3-LfQuasiRefl.
	In the second case, we have $z \rr{R} z$ by 3-LfQuasiRefl;
	but $y,y$ and $y,z$ incomparable should imply $z,z$ incomparable.
\qed
\end{enumerate}
}



\LEMMA{lu, de}{
\LABEL{lu, de}%
Let ${\it Sym} = \set{0, 1, 6, 7, 8, 9, E, F}$
denote the set of all unary operations that
are guaranteed to yield always a symmetric relation 
by Lem.~\REFF{optn univ}{17}.
Let $p, q, r$ be unary operations such that
\begin{enumerate}
\item $(p \land q) = 0$,
\item $(p \lor q) \supseteq (\lnot r)$, and
\item $(\lnot r) \in {\it Sym}$, or $p,q \in {\it Sym}$.
\end{enumerate}
Then $p$-LfUnique and $q$-LfUnique implies $r$-Dense.
}
\PROOF{
First, we prove the case $(\lnot r) \in {\it Sym}$.
Assume for contradiction that $x \ro{R}{r} z$,
but (using symmetry of $\ro{R}{\lnot r}$)
that $y \ro{R}{\lnot r} x \lor y \ro{R}{\lnot r} z$ for all $y$.
Let $p' = (p \land \lnot r)$ and $q' = (q \land \lnot r)$,
then $p'$ and $q'$ are disjoint like $p$ and $q$,
and $(\lnot r) = (p' \lor q')$,
so we can make the following case
distinction for all $y$:
\begin{enumerate}
\item\LABEL{lu, de 1} If $y \ro{R}{p'} x$, 
	then $y \ro{R}{p} x$ by Lem.~\REFF{optn bool}{4};
\item\LABEL{lu, de 2} if $y \ro{R}{q'} x$, 
	then similarly $y \ro{R}{q} x$;
\item\LABEL{lu, de 3} if $y \ro{R}{p'} z$, then $y \ro{R}{p} z$;
\item\LABEL{lu, de 4} if $y \ro{R}{q'} z$, then $y \ro{R}{q} z$.
\end{enumerate}
Choosing $5$ distinct elements $\:,15{y_\i}$, one of the cases must
appear twice.
Double appearance of case~\REF{lu, de 1} or~\REF{lu, de 3}
contradicts $p$-LfUnique,
double appearance of case~\REF{lu, de 2} or~\REF{lu, de 4}
contradicts $q$-LfUnique.

For the case $p,q \in {\it Sym}$, we use a similar reasoning:
Assume for contradiction that $x \ro{R}{r} z$,
but $x \ro{R}{\lnot r} y \lor y \ro{R}{\lnot r} z$ for all $y$.
Defining $p', q'$ as above, we get these cases:
\begin{enumerate}
\item\LABEL{lu, de 5} If $x \ro{R}{p'} y$, 
	then $x \ro{R}{p} y$ by Lem.~\REFF{optn bool}{4},
	hence $y \ro{R}{p} x$ by symmetry of $\ro{R}{p}$;
\item\LABEL{lu, de 6} if $y \ro{R}{q'} x$, 
	then similarly $y \ro{R}{q} x$;
\item\LABEL{lu, de 7} if $y \ro{R}{p'} z$, then $y \ro{R}{p} z$;
\item\LABEL{lu, de 8} if $y \ro{R}{q'} z$, then $y \ro{R}{q} z$.
\end{enumerate}
Again, double appearance of~\REF{lu, de 5} or~\REF{lu, de 7}
contradicts $p$-LfUnique,
and double appearance of case~\REF{lu, de 6} or~\REF{lu, de 8} 
contradicts $q$-LfUnique.
\qed
}

\LEMMA{3-Implications group 5}{%
\LABEL{3-Implications group 5}%
\begin{enumerate}
\item\LABEL{3-Implications group 5 gtlu}%
	(\QDef{gtlu})
	C-SemiOrd2 $\land$ SemiOrd2 $\ra$ 2-SemiOrd2
\item\LABEL{3-Implications group 5 uydr}%
	(\QDef{uydr})
	On a domain of $\geq 5$ elements,
	3-LfUnique $\land$ 4-LfUnique $\ra$ 8-Dense
\end{enumerate}
}
\PROOF{
\begin{enumerate}
\item
	The proof extensively uses boolean connectives on operations
	(Lem.~\REF{optn bool}) to shorten the proof presentation.
	Let $x_1 \ro{R}{2} x_2$ and $x_2 \ro{R}{2} x_3$,
	and let $w$ be given; 
	we have to show $w \ro{R}{6} x_i$
	for some $i \in \set{1,2,3}$.
	Assume for contradiction there is no such $i$,
	that is, $w \ro{R}{9} x_i$ for all $i$..
	We distinguish two cases:
	\begin{itemize}
	\item If $w \rr{R} w$,
		\\
		then $w \rr{R} w$ and $w \rr{R} w$
		implies $w \ro{R}{7} y$ for each $y$, 
		by SemiOrd2.
		Hence $w \rr{R}{1} x_i$
		for each $i$ by assumption.
		But $x_3 \ro{R}{C} x_2$ and $x_2 \ro{R}{C} x_1$
		requires $w \ro{R}{E} x_i$
		for some $i$, by C-SemiOrd2,
		which is a contradiction.
	\item If $\lnot w \rr{R} w$,
		\\
		then $w \ro{R}{C} w$ and $w \ro{R}{C} w$ 
		implies $w \ro{R}{E} y$
		for each $y$, by C-SemiOrd2.
		Hence, $w \rr{R}{8} x_i$
		for each $i$ by assumption.
		But this contradicts SemiOrd2.
	\end{itemize}
\item
	Follows from Lem.~\REF{lu, de}
	with $p=3$, $q=4$, and $r=8$,
	since $(\lnot r) = 7 \in {\it Sym}$.
\qed
\end{enumerate}
}

\LEMMA{3-Implications group 6}{%
\LABEL{3-Implications group 6}%
\begin{enumerate}
\item\LABEL{3-Implications group 6 vxfj}
	(\QDef{vxfj})
	7-LfUnique $\land$ 6-LfSerial $\ra$ ASym
\item\LABEL{3-Implications group 6 xlfr}
	(\QDef{xlfr})
	1-LfQuasiRefl $\land$ C-Trans $\ra$ SemiOrd1
\end{enumerate}
}
\PROOF{
\begin{enumerate}
\item
	Let $x \rr{R} y$,
	assume for contradiction $y \rr{R} x$,
	i.e.\ $x \ro{R}{1} y$.
	By 6-LfSerial, we obtain some $w$ such that
	$w \ro{R}{6} x$.
	Since $x \ro{R}{1} y$ implies $x \ro{R}{7} y$,
	and similarly $w \ro{R}{6} x$ implies $x \ro{R}{7} w$,
	we get $w=y$ by 7-LfUnique.
	But $x \ro{R}{1} y$ and
	$x \ro{R}{6} y$ (using symmetry of $\ro{R}{6}$)
	implies $x \ro{R}{0} y$ by Lem.~\REFF{optn bool}{2},
	a contradiction.
\item
	Let $w \rr{R} x$ and $x \ro{R}{8} y$ and $y \rr{R} z$.
	Assume for contradiction $\lnot w \rr{R} z$.
	Then also $z \rr{R} y$,
	since else,
	we had $\lnot w \rr{R} y$ by C-Trans,
	and $\lnot w \rr{R} x$ by C-Trans again.
	Applying 1-LfQuasiRefl, we obtain $y \ro{R}{1} y$.
	But $\lnot x \rr{R} y$ and $\lnot y \rr{R} x$ 
	imply the contradiction
	$\lnot y \rr{R} y$ by C-Trans.
\qed
\end{enumerate}
}

\LEMMA{3-Implications group 7}{%
\LABEL{3-Implications group 7}%
\begin{enumerate}
\item\LABEL{3-Implications group 7 ewvq}
	(\QDef{ewvq})
	2-SemiOrd1 $\land$ 2-LfSerial $\ra$ 7-Dense
\item\LABEL{3-Implications group 7 g0rw}%
	(\QDef{g0rw})
	2-SemiOrd2 $\land$ 7-Trans $\ra$ C-SemiOrd2
\item\LABEL{3-Implications group 7 voye}%
	(\QDef{voye})
	On a domain of $\geq 5$ elements,
	1-LfUnique $\land$ 6-LfUnique $\ra$ C-Dense
\item\LABEL{3-Implications group 7 voyr}%
	(\QDef{voyr})
	On a domain of $\geq 5$ elements,
	1-LfUnique $\land$ 6-LfUnique $\ra$ 8-Dense
\item\LABEL{3-Implications group 7 wfyp}%
	(\QDef{wfyp})
	On a domain of $\geq 5$ elements,
	1-LfUnique $\land$ 8-LfUnique $\ra$ 6-Dense
\end{enumerate}
}
\PROOF{
\begin{enumerate}
\item
	Let $x \ro{R}{7} y$, w.l.o.g.\ let $x \rr{R} y$.
	By 2-LfSerial, we obtain $x', y'$
	with $x' \ro{R}{2} x$ and $y' \ro{R}{2} y$.
	We have two cases:
	\begin{itemize}
	\item If $x,y'$ are comparable w.r.t.\ $R$,
		then $x \ro{R}{7} y'$, 
		that is, $y'$ is an intermediate element,
		and we are done.
	\item Else, we apply 2-SemiOrd1 to
		$x' \ro{R}{2} x$, 
		$x,y'$ incomparable w.r.t.\ $\ro{R}{2}$, 
		$y' \ro{R}{2} y$
		to obtain $x' \ro{R}{2} y$,
		hence $x' \ro{R}{7} y$,
		hence $x'$ being an intermediate element.
	\end{itemize}
\item
	Assume for contradiction
	$\lnot x \rr{R} y \land \lnot y \rr{R} z$,
	but 
	$w \rr{R} x \land x \rr{R} w 
	\land w \rr{R} y \land y \rr{R} w 
	\land w \rr{R} z \land z \rr{R} w$
	for some $x,y,z,w$.
	By 7-Trans, we get $x \rr{R} y \lor y \rr{R} x$ 
	and $y \rr{R} z \lor z \rr{R} y$,
	that is, $y \rr{R} x$ and $z \rr{R} y$.
	That is, $z \ro{R}{2} y \land y \ro{R}{2} x$,
	but $w$ is incomparable (w.r.t.\ $R^2$) to $x$, $y$, and $2$,
	contradicting 2-SemiOrd2.
\item
	Follows from Lem.~\REF{lu, de}
	with $p=1$, $q=6$, and $r=C$,
	since $p,q \in {\it Sym}$.
\item
	Follows from Lem.~\REF{lu, de}
	with $p=1$, $q=6$, and $r=8$.
	since $(\lnot r) = 7 \in {\it Sym}$.
\item
	Follows from Lem.~\REF{lu, de}
	with $p=1$, $q=8$, and $r=6$,
	since $(\lnot r) = 9 \in {\it Sym}$.
\qed
\end{enumerate}
}

\LEMMA{cjdj}{%
If $R$ is ASym, then $R^1$ is CoRefl,
LfEucl, LfUnique, Sym, AntiTrans, ASym, AntiSym, Trans, SemiOrd1, Dense
(\QDef{cjdj}).
}
\PROOF{
From the assumption follows that $R^1$ is empty,
hence has all claimed properties
by \OEF{Exm.~74}{48}{}.
\qed
}

\THEOREM{$3$-implication axioms}{%
\LABEL{$3$-implication axioms}%
Recall that a ``$3$-implication'' is a formula of the form
$\forall R. \; {\it lprop}_1(R) \land {\it lprop}_2(R)
        \ra {\it lprop}_3(R)$,
where ${\it lprop}_i$ are unnegated lifted properties.
Consider the inference system from Def.~\REF{Trivial inference rules}
and~\REF{Trivial initialization rules};
inferences are understood to apply the equivalences from
Thm.~\REF{optn equiv} as needed.

A $3$-implication
is valid (w.r.t.\ a relation domain of $\geq 7$ elements)
iff it
can be inferred from the set of axioms shown in 
Fig.~\REF{Minimal axiom set}.
}
\PROOF{
``$\La$'':
All axioms have been proven to be valid, by EProver, or manually
(Sect.~\REF{Some proofs of axioms}), 
or both.
The inference rules of Def.~\REF{Trivial inference rules}
and~\REF{Trivial initialization rules} are obviously sound.

``$\Ra$'':
All $3$-implications that couldn't be inferred were proven to be
invalid, either by a computer-generated counter-example on a domain of
$5$ elements, or a manual counter-example w.r.t.\ a finite
(Sect.~\REF{Finite counter-examples}) 
or infinite (Sect.~\REF{Infinite counter-examples}) domain.
\qed
}

\clearpage

\section{References}

\bibliographystyle{elsarticle-num-names}
\bibliography{lit}

\end{document}